\documentclass{article}

\usepackage{amsfonts,amsmath}
\usepackage{amssymb,latexsym}
\usepackage[dvips]{epsfig}
\usepackage{amscd}

\title{ \hspace{2.3in} {\it {\normalsize To my teacher Igor Frenkel}} \\
 \hspace{1in}  \\
{ \bf A categorification of the Jones polynomial}
\footnotetext{This paper will appear in Duke Mathematical Journal}}
\author{Mikhail Khovanov}  

\newtheorem{prop}{Proposition}
\newtheorem{theorem}{Theorem}
\newtheorem{lemma}{Lemma}
\newtheorem{corollary}{Corollary}
\newtheorem{definition}{Definition}
\newtheorem{conjecture}{Conjecture} 
\newcommand{\oplusop}[1]{{\mathop{\oplus}\limits_{#1}}}

\begin{document}
\maketitle
\baselineskip 12pt

\tableofcontents
\bigbreak

\def\C{\mathbb C}
\def\R{\mathbb R}
\def\N{\mathbb N}
\def\Z{\mathbb Z}
\def\Q{\mathbb Q}
\def\F{\mathbb F}
\def\S{\mathbb S}
\def\Zq{\Z[q,q^{-1}]}
\def\l{\lbrace}
\def\r{\rbrace}
\def\o{\otimes}
\def\D{\Delta}
\def\O{\mathcal{O}}
\def\lra{\longrightarrow}

\def\cH{{\cal H}} 
\def\cA{{\cal A}}
\def\cC{{\cal C}}
\def\cV{{\cal V}}
\def\whc{\widehat{\chi}}
\def\wta{\widetilde{\L}}

\def\drawing#1{\begin{center} \epsfig{file=#1} \end{center}}

\def\mo{\mathbf{1}}

\def\mc{\mathcal} 
\def\mf{\mathfrak}

\def\CB{\mathcal B}
\def\MM{\mathcal M} 
\def\M1{\mc{M}_1}
\def\I{\mathcal I} 
\def\L{\mathcal L} 

\def\C{\overline{C}}

\def\I{\mc{I}}
\def\J{\mc{J}}

\def\Hom{\textrm{Hom}}
\def\Kom{\mbox{Kom}}
\def\Rm{R{\mbox{-mod}_0}}
\def\Rmb{R{\mbox{-mod}}}
\def\KR{\Kom(R\mbox{-mod})}

\newcommand{\ddeg}{{\mathrm{deg}}}
\newcommand{\ddim}{{\mathrm{dim}}}
\newcommand{\Tor}{{\mathrm{Tor}}}
\newcommand{\Mor}{{\mathrm{Mor}}}
\newcommand{\Tot}{{\mathrm{Tot}}}
\newcommand{\wt}{{\mathrm{wt}}}
\newcommand{\hl}{{\mathrm{hl}}}
\newcommand{\cm}{{\mathrm{cm}}}
\newcommand{\cl}{{\mathrm{cl}}}
\newcommand{\kker}{{\mathrm{ker }}}

\section{Introduction}

During the past 15 years many new 
 structures have arisen 
in the topology of low-dimensional manifolds: the Jones and HOMFLY 
polynomials of links, Witten-Reshetikhin-Turaev invariants of 3-manifolds, 
Floer homology groups of homology 3-spheres, Donaldson and Seiberg-Witten 
invariants of 4-manifolds. These invariants of 3- and 4-manifolds 
naturally split into two groups. Members of the first group are 
 combinatorially defined invariants of knots and 3-manifolds, 
such as various link polynomials, finite type invariants and
quantum invariants of 3-manifolds. 
Floer and Seiberg-Witten 
homology groups of 3-manifolds and Donaldson-Seiberg-Witten invariants 
of 4-manifolds constitute the second group. While invariants from the first 
group have a combinatorial description and in each instance can 
be computed algorithmically, invariants from the second group 
are understood through moduli spaces of solutions of 
suitable differential-geometric equations and the 
infinite-dimensional Morse theory and have evaded
all attempts at a finite combinatorial definition. These invariants 
have been computed for many 3- and 4-manifolds, yet the methods of 
computation use some extra structure on these manifolds, 
such as Seifert fibering or complex structure and the problem of 
finding an algorithmic construction of these invariants remains open. 

  It is probably due to this striking difference in the 
  origins and computational complexity that so far 
  not many  direct relations have been found
  between invariants from different groups. 
  The most notable connection is the 
  Casson invariant of homology 3-spheres [AM], which is equal to the 
  Euler characteristic of Floer homology [F]. Yet the Casson 
  invariant is
  computable and intimately related to the Alexander polynomial of 
  knots and (see [M]) to Witten-Reshetikhin-Turaev invariants, which are 
  examples of invariants from the first group. 
  A similar relation has recently been discovered between 
  Seiberg-Witten invariants and Milnor torsion of 3-manifolds ([MT]). 
  In summary, Euler characteristics of 
  Floer and Seiberg-Witten homology groups bear an algorithmic 
  description, while no such procedure is known for 
  finding the groups themselves. 

  A speculative question now comes to mind:
  quantum invariants 
  of knots and 3-manifolds tend to have good integrality properties. 
  What if these invariants can be interpreted as Euler characteristics 
  of some homology theories of 3-manifolds? 

  Our results suggest that such an interpretation 
  exists for the Jones polynomial of links in 3-space ([Jo]). We give an 
  algorithmic procedure that to a generic plane projection $D$ of 
  an oriented link $L$ in $\R^3$ associates cohomology groups 
  $\cH^{i,j}(D)$ that depend on two integers $i,j.$ 
  If two diagrams $D_1$ and $D_2$ of the same link $L$ 
  are related by a Reidemeister move, a canonical 
  isomorphism of groups $\cH^{i,j}(D_1)$ and $\cH^{i,j}(D_2)$ is constructed. 
  Thus, isomorphism classes of these groups are invariants of the link 
  $L.$ These groups are finitely generated and may have non-trivial 
  torsion. Tensoring these groups with $\Q$
    we get a two-parameter family 
  $\{ \mbox{dim}_{\Q}(\cH^{i,j}(D)\o \Q)\}_{i,j\in \Z}$ 
  of integer-valued link invariants. 
  
  From our construction of groups $\cH^{i,j}(D)$ we immediately 
  conclude that the graded Euler characteristic 
  \begin{equation}  
  \sum_{i,j}(-1)^i q^j \mbox{dim}_{\Q}(\cH^{i,j}(D)\o \Q) 
  \end{equation} 
   is equal, 
  up to a simple change of variables, to the Jones polynomial of $L,$
  multiplied by $q+q^{-1}.$   

  We conjecture that not just the isomorphism classes of $\cH^{i,j}(D)$
  but the groups themselves are invariants of links. We will 
  consider this conjecture in a subsequent paper. 
  
   To define cohomology groups $\cH^{i,j}(D)$ we start with the Kauffman 
   state sum model [Ka] for the Jones polynomial and then, roughly 
   speaking, turn all integers 
   into complexes of abelian groups.  
  In the Kauffman model a link is projected generically onto 
  the plane so that the projection has a finite number of double 
  transversal intersections. There are two ways to ``smooth'' the 
  projection near the double point, i.e., erase the intersection 
  of the projection with a small neighbourhood of the double point 
  and connect the four resulting ends by a pair of simple, nonintersecting 
  arcs: 

  \drawing{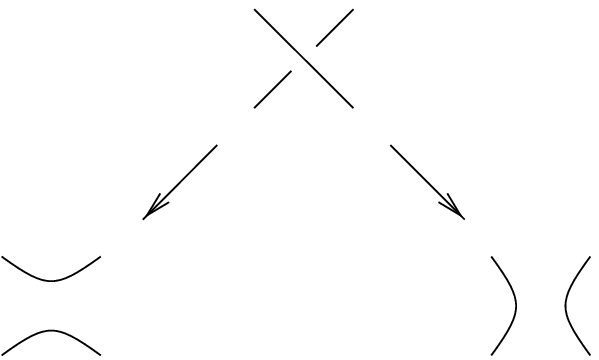}

  A diagram $D$ with $n$ double points admits $2^n$ resolutions of 
  these double points. Each of the resulting diagrams is a collection 
  of disjoint simple closed curves on the plane. In [Ka] 
  Kauffman associates Laurent polynomial 
  $(-q-q^{-1})^k$ to a collection of $k$ simple curves and 
  then forms a weighted sum of these numbers over all $2^n$ resolutions. 
  After normalization, Kauffman obtains the Jones polynomial of the 
  link $L.$ 
  The principal constant in this construction is $-q-q^{-1}, $ the 
  number associated to a simple closed curve. 
  
  In our approach 
  $q+q^{-1}$ becomes a certain module $A$ over the base ring $\Z[c].$ 
  In detail, we work over the graded ring $\Z[c]$ of polynomials in $c,$ 
  where $c$ has degree $2,$ and we define  $A$ to be a free $\Z[c]$-module 
  of rank two with generators in degrees $1$ and $-1.$ 
  This is the object we associate to a simple closed curve in the plane. 

  Given a diagram $D$, to each resolution of all double points 
  of $D$ we associate the graded $\Z[c]$-module 
  $A^{\o k},$ where $k$ is the number of curves in the 
  resolution. Then we glue these modules over all $2^n$ resolutions into 
  a complex $C(D)$ of graded $\Z[c]$-modules. The gluing maps come 
  from commutative algebra and cocommutative coalgebra structure on 
  $A.$ When two diagrams $D_1$ and 
  $D_2$ are related by a Reidemeister move, we construct a 
  quasi-isomorphism between the complexes $C(D_1)$ and $C(D_2).$ 
  
The 
cohomology groups $H^i(D)$ of the complex $C(D)$ are graded $R$-modules   
and we prove that isomorphism classes of $H^i(D)$ do not depend of 
the choice of a diagram of the link. We then look at some elementary 
properties of these groups and introduce several 
cousins of $H^i(D).$

  {\it Outline of the paper.} In Section~\ref{section-prelim} 
  we define an algebra $A$ over the ring $R=\Z[c]$ 
  and use $A$ to construct a 2-dimensional topological 
  quantum field theory. In our case this 
  topological quantum field theory is a functor from
  the category of two-dimensional cobordisms between one-dimensional 
  manifolds to  the category of graded $\Z[c]$-modules.  
  In Section~\ref{Kauffman-bracket} we review the Kauffman
  state sum model [Ka] for the Jones polynomial of oriented links. 

  In Section~\ref{diagrams-to-comcubes}, we  associate 
  a complex of $\Z[c]$-modules to a plane diagram of a link. 
  As an intermediate step, to a diagram $D$  
  we associate a commutative cube $V_D$ of $\Z[c]$-modules and maps between 
  them, i.e., we consider an $n$-dimensional cube with its edges 
  standardly oriented and, given a plane projection with $n$ double points
  of a link, 
  to each vertex of the cube we associate a $\Z[c]$-module and to each 
  oriented edge a map of modules so that all square facets of this diagram 
  are commutative squares. 
  This is done in Section~\ref{diagrams-to-comcubes}.  
  In the same section we pass from commutative cubes to 
  complexes of $\Z[c]$-modules and to a diagram $D$ 
  we associate a complex $C(D)$ of graded $\Z[c]$-modules. 

  Earlier, in Section~\ref{section-Cubes}, we review the 
  notions of a commutative cube and a map between commutative cubes. 

In Section~\ref{Reidemeister-moves} we review Reidemeister moves. 
In Section~\ref{section-transformations}, 
which is the technical core of the paper, to a Reidemeister move 
between diagrams $D_1$ and $D_2$ 
we associate a quasi-isomorphism between the complexes $C(D_1)$ and 
$C(D_2).$ These isomorphisms seem to be canonical. We conjecture 
that the quasi-isomorphisms are coherent, which would naturally 
associate cohomology groups to links. Our quasi-isomorphism result 
shows that the isomorphism classes of the cohomology groups are 
invariants, but not necessarily that the groups are functorial under 
link isotopy. 
  
We define $H^i(D)$ to the be $i$-th cohomology group of 
the complex $C(D).$ These cohomology groups are graded $\Z[c]$-modules, 
and the isomorphism class of each $H^i(D)$ is a link invariant. 
If we split these groups into the direct sum of their graded components, 
\begin{equation} 
H^i(D)= \oplusop{j\in \Z}H^{i,j}(D),
\end{equation} 
we get a two-parameter family of ``abelian group valued'' link  
invariants. These results are 
stated earlier, at the end of  Section~\ref{diagrams-to-comcubes}, 
as Theorems~\ref{closed-i} and \ref{closed-i-and-j}.  
For a diagram $D,$ 
the groups $H^{i,j}(D)$ are trivial for $j\ll 0.$ Moreover, 
for each $j$ only finitely many of the groups $H^{i,j}(D)$ are non-zero. 
  Consequently, the graded Euler characteristic of $C(D),$ defined as 
  \begin{equation} 
  \widehat{\chi}(C(D))= 
  \sum_{i,j\in \Z}(-1)^i q^j \mbox{dim}_{\Q}(H^{i,j}(D)\o \Q), 
  \end{equation} 
 is well-defined as a Laurent series in $q.$ Since our construction 
  of $C(D)$ lifts Kauffman's construction of the Jones polynomial, 
  it is not surprising that the graded Euler characteric of $C(D)$ is related 
  to the Jones polynomial. 
  Namely, $\widehat{\chi}(C(D)),$ multiplied by 
  $\frac{1-q^2}{q+ q^{-1}},$ is equal, after a simple change of 
  variables, to the Jones polynomial of $L.$

  In Section~\ref{c-is-zero} 
   we explain how a version of our construction when the base
  algebra $\Z[c]$ is reduced to $\Z$ by taking $c=0$ 
  produces graded cohomology groups 
  $\cH^{i,j}(D).$ 
    The complex which is used to define $\cH^{i,j}(D)$ is given by 
   tensoring $C(D)$ with $\Z$ over $\Z[c].$ 
  As before, the isomorphism classes of these groups 
  are invariants of links. These groups are ``smaller'' than the 
  groups $H^{i,j}(D).$ In particular, for each $D,$ these groups are non-zero 
  for only finitely many pairs $(i,j)$ of integers. 
As with the groups groups $H^{i,j}(D),$ the graded Euler characteristic 
  \begin{equation} 
  \sum_{i,j}(-1)^iq^j \mbox{dim}_{\Q}(\cH^{i,j}(D)\o \Q),
  \end{equation} 
divided by $q+q^{-1}$, is equal to the 
Jones polynomial of the link represented by the diagram $D.$  
In Section~\ref{spectral-sequence} 
we exhibit a spectral sequence whose $E_1$-term 
is made of $\cH(D)$ and which converges
to $H(D).$ Apparently, $H(D)$ is
a kind of $S^1$-equivariant version of the groups $\cH^{i,j}(D).$ 
  In Section~\ref{crossing-number} we use these cohomology groups to
  reprove a result of Thistlethwaite on the crossing number of 
  adequate links. 
  
  If a link in $\R^3$ has cohomology groups, 
  then cobordisms between links, i.e.,
   surfaces embedded in $\R^3\times [0,1],$ 
   should provide maps between the associated 
  groups.   A surface embedded in the 4-space can be visualized as a sequence 
  of plane projections of its 3-dimensional sections (see [CS]). 
  Given such a presentation $J$ of a compact oriented surface $S$ 
  properly embedded in $\R^3\times [0,1]$ with the boundary of $S$ 
  being the union of two links $L_0\subset \R^3\times \{ 0\} $ and 
  $L_1 \subset \R^3\times \{ 1\},$   we explain 
  in Section~\ref{four-d-cobordisms} 
  how to associate to $J$ a map 
  of cohomology groups 
  \begin{equation} 
  \theta_J: H^{i,j}(D_0)\lra H^{i, j + \chi(S)}(D_1), \hspace{0.4in} 
  i,j\in \Z,
  \end{equation} 
  $\chi(S)$ being the Euler characteristic of the surface $S$ 
  and $D_0$ and $D_1$ -- diagrams of $L_0$ and $L_1$ induced by $J.$  
  We conjecture that, up to an overall minus sign, 
   this map does not depend on the choice of $J,$ in other words, 
  $\pm \theta_J$ behaves invariantly  under isotopies of $S.$ 

  If this conjecture is true, we get a 4-dimensional topological 
  quantum field theory, restricted to links in $\R^3$ and 
  $\R^3\times [0,1]$-cobordisms between them. Because the theory 
  has a combinatorial definition, all cohomology groups and maps between 
  them will be algorithmically computable. If successful, this program 
  will realize the Jones polynomial as the Euler characteristic of 
  a cohomology theory of link cobordisms. 

  Section~\ref{section-invariants-from-modules} presents 
  mild variations on cohomology groups $H^i(D)$ and $\cH^{i,j}(D).$  
  There we switch from links to $(1,1)$-tangles.
   We consider the category $A\mbox{-mod}_0$ of graded $A$-modules
   and grading-preserving homomorphisms between them. 
   Given a plane diagram $D$ of a $(1,1)$-tangle 
   $L$ and a graded $A$-module $M,$ in Section~\ref{section-loong-links} 
  we define cohomology groups $H^i(D,M)$ which are graded $A$-modules. 
  The arguments of 
  Sections~\ref{section-diagrams}-\ref{section-transformations} 
  go through without a single alteration and show that isomorphism 
  classes of $H^i(D,M)$ do not depend on the choice of $D$ and 
  are invariants of the underlying $(1,1)$-tangle $L$. In fact, to every 
  $(1,1)$-tangle and an integer $i$ 
  we associate an isomorphism class of functors from 
  the category of graded $A$-modules to itself.  
  
  {\it Motivations for this work and its 
  relations to representation theory.}   
   What is the representation-theoretical meaning of the cohomology groups 
  $H^{i,j}(D)?$
  The Jones polynomial of links is encoded in the finite-dimensional 
  representation theory of the quantum group $U_q(\mf{sl}_2).$  
  It was shown in [FK] and [K] that the integrality and positivity 
  properties of the Penrose-Kauffman $q$-spin networks 
  calculus, 
  of which the Jones polynomial is a special instance, are related 
  to Lusztig canonical bases in tensor products of finite-dimensional 
  $U_q(\mf{sl}_2)$-representations. Lusztig's theory [L], among other 
  things, says that various structure coefficients of quantum groups 
  can be obtained as dimensions of cohomology groups 
  of sheaves on quiver varieties. 
  This suggests a ``categorification'' of quantum groups and their 
  representations, i.e. that there exist certain categories 
  and 2-categories whose Grothendieck groups produce quantum groups 
  and their representations. 

Louis Crane and Igor Frenkel [CF] conjectured  
that quantum  $\mf{sl}_2$ invariants of 3-manifolds 
can be lifted to a 4-dimensional topological quantum field theory 
  via canonical bases of Lusztig. They also introduced a notion of 
  Hopf category and associated to it 4-dimensional invariants. 
  Representations of a Hopf category form a 2-category, 
  and a relation between 2-categories and invariants of 2-knots 
  in $\R^4$ was established in [Fs].  
 
In joint work with Bernstein and Frenkel [BFK], 
  we propose a categorification of 
  the representation theory of  $U_q(\mf{sl}_2)$ 
  via categories of highest weight representations for 
  Lie algebras $\mf{gl}_n$ for all natural $n.$ 
  This approach can be viewed as an algebraic counterpart of Lusztig's 
  original geometric approach to canonical bases. Motivated by 
  the geometric constructions of [BLM] and [GrL], we obtain a 
  categorification of the Temperley-Lieb algebra and Schur quotients of 
  $U(\mf{sl}_2)$ via projective and Zuckerman functors. 
   We consider categories $\O_n$ 
  that are direct sums of certain singular blocks  of the category 
  $\O$ for $\mf{gl}_n.$ Given a tangle $L$ in the 3-space with $n$ 
  bottom and $m$ top ends and a plane projection $P$ of $L,$ 
  we associate to $P$ a functor between derived categories $D^b(\O_n)$ 
  and $D^b(\O_m).$ Properties of these functors suggest
  that their isomorphism classes, up to shifts in the derived 
  category, are invariants of tangles. When the tangle is a link $L,$ 
  we expect to get cohomology groups $\mathbb{H}^i(L)$ 
  as invariants of links. 
  These groups will be a special case of the cohomology groups 
  constructed in this paper: conjecturally
  \begin{equation}
  \mathbb{H}^i(L)= \oplusop{j} (\cH^{i,j}(L)\o \mathbb C) .
  \end{equation}  

  {\it Acknowledgements.}  This work was started during  a
  visit to the Institut des Hautes Etudes Scientifiques
  and finished at the Institute for Advanced 
 Study. I am grateful to these institutions for creating a wonderful 
 working atmosphere. During my stay at the Institute for Advanced 
  Study I was supported by grant DMS 9304580  from the NSF.

  I am indebted to Joseph Bernstein for interesting discussions, to 
  Greg Kuperberg for reading and correcting the first version of the 
   manuscript and explaining to me a natural way to hide minus signs, 
  and to Oliver Dasbach and Arkady Vaintrob for pointing out that 
  Corollary~\ref{due-to-Th} was proved by Thistlethwaite [T]. 

  On numerous occasions 
  Igor Frenkel, who was my supervisor at Yale University, 
  advised  me to look for a lift of the Penrose-Kauffman 
  quantum spin networks calculus to a calculus of surfaces in $\S^4.$ 
  This work can be seen as a partial
  answer to his questions. It is a pleasure to dedicate this paper 
  to my teacher Igor Frenkel.

\section{Preliminaries}
 \label{section-prelim}

\subsection{ The ring $R$} 
\label{here-is-R}

Let $R=\Z[c]$ denote the ring of polynomials with integral coefficients. 
Introduce a $\Z$-grading on $R$ by 
\begin{equation} 
\label{grading-1}
\ddeg(1)=0,  \hspace{0.3in} \ddeg (c) =2. 
\end{equation} 

Denote by $\Rm$ the abelian category of 
graded $R$-modules. Denote the $i$-th graded component of an 
object $M$ of $\Rm$ by $M_i.$ Morphisms in 
the category $\Rm$ are grading-preservings homomorphisms of modules. 
For $n\in \Z$ denote by $\{ n \}$ the automorphism of $\Rm$ 
given by shifting the grading down by $n.$ Thus for a graded $R$-module 
$N=\oplus_i N_i,$ the shifted module $N\{ n \} $ has 
graded components $N\{ n\} _i = N_{i+n}.$ 

In this paper we will sometimes consider graded, rather than just 
grading-preserving, maps. A map $\alpha: M\to N$ 
of graded $R$-modules is called graded of degree $i$ if 
$\alpha(M_j)\subset N_{i+j}$ for all $j\in \Z.$ 

Let $\Rmb$ be the category of graded $R$-modules and graded maps 
between them. This category has the same objects as the category 
$\Rm,$ but more morphisms. It is not an abelian category. 

A graded map $\alpha$ 
is a morphism in the category $\Rm$ if and only if the degree of $\alpha$ 
is $0.$ At the end we will favor grading-preserving maps 
and when at some point we look at a graded map $\alpha: M\to N$ of 
degree $i,$ later we will make it grading-preserving by appropriately 
shifting the degree of one of the modules. For example, $\alpha$ gives 
rise to a grading-preserving map $M\to N\{ i \},$ also denoted $\alpha.$ 

Let $M$ be a finitely-generated graded $R$-module. As an abelian group, 
$M$ is the direct sum of its graded components: $M= \oplusop{j\in \Z} 
M_j,$ where each $M_j$ is a finitely-generated abelian group. Define 
{\it the graded Euler characteristic } $\whc(M)$ of $M$ by 
\begin{equation} 
\whc(M)= \sum_{j\in \Z} \ddim_{\Q}(M_i\o_{\Z} \Q) q^j 
\end{equation} 
Since
\begin{equation} 
\whc(R)= 1 + q^2+q^4 + \dots = \frac{1}{1-q^2}, 
\end{equation} 
$\whc(M)$ is not, in general, a Laurent polynomial in $q$, but 
an element of the Laurent series ring.   
Moreover, for any $M$ as above, 
there are Laurent polynomials $a,b\in \Z[q,q^{-1}]$ 
such that 
\begin{equation} 
\whc(M) = a + \frac{b}{1-q^2}
\end{equation}

\subsection{ The algebra $A$}
\label{here-is-A} 

Let $A$ be a  free graded $R$-module of rank $2$ spanned by $\mo$ 
and $X$ with 
\begin{equation}
\label{grading-2}
 \ddeg(\mo) = 1, \hspace{0.3in}  \ddeg (X) = -1 
\end{equation} 
We equip $A$ with a commutative algebra structure with  
the unit $\mo $ and multiplication 
\begin{equation}
\mo  X=X \mo  =X, X^2=0.
\end{equation}
We denote by $\iota$ the unit map $R\to A$ which sends $1$ to $\mo.$ 
This map is a graded map of graded $R$-modules and it increases 
 the degree by $1.$

We equip $A$ with a coalgebra structure with a  
coassociative cocommutative comultiplication 
\begin{eqnarray} 
& \Delta(\mo ) & = \mo \o X + X \o \mo + c X \o X \\
& \Delta(X) & = X\o X 
\end{eqnarray} 
and a counit 
\begin{equation}
\epsilon(\mo ) =-c   , \hspace{0.5in} \epsilon(X) = 1 .
\end{equation}

$A,$ equipped with these structures, is not a Hopf algebra. 
Instead, the identity 
\begin{equation}
\label{m-and-delta}
\Delta \circ m = (m \o 
\mbox{Id} ) \circ (\mbox{Id} \o \Delta)
\end{equation}
holds.

Grading $\ddeg,$ given by (\ref{grading-1}),(\ref{grading-2}), 
induces a grading, also denoted $\ddeg,$
 on tensor powers of $A$ by 
\begin{equation} 
\ddeg(a_1\o \dots \o a_n) = \ddeg (a_1)+ \dots +\ddeg (a_n)
\mbox{ for } a_1,\dots , a_n \in A
\end{equation} 
Here and further on all tensor products are taken over the ring $R$
unless specified otherwise.  

We next describe the effect of the structure maps $\iota, m, 
\epsilon, \Delta$ on the gradings.
We say that a map $f$ between 
two graded $R$-modules  $V=\oplus V_n$ and $W=\oplus W_n$ 
has degree $k$ if $f(x)\in W_{n+k}$ whenever $x\in V_n.$ 

\begin{prop} Each of the structure maps $\iota, m, \epsilon, 
\Delta$ is graded relative to the grading $\ddeg.$ Namely, 
\begin{equation}
\ddeg(\iota)=1,\hspace{0.3in} \ddeg(m)=  -1 , \hspace{0.3in}
\ddeg(\epsilon)=1, \hspace{0.3in} \ddeg(\Delta)=-1.
\end{equation} 
\end{prop}

\begin{prop} 
\label{decomp-one} 
We have an $R$-module decomposition 
\begin{equation}
A\o A = (A\o \mo) \oplus \Delta (A)
\end{equation}
 which respects the grading $\ddeg.$  
\end{prop}

\subsection{Algebra $A$ and (1+1)-dimensional cobordisms} 
\label{A-and-surfaces} 

Consider the surfaces $S_2^1,S_1^2,S_0^1,S_1^0,S_2^2$ and $S_1^1$,  
depicted below 

\drawing{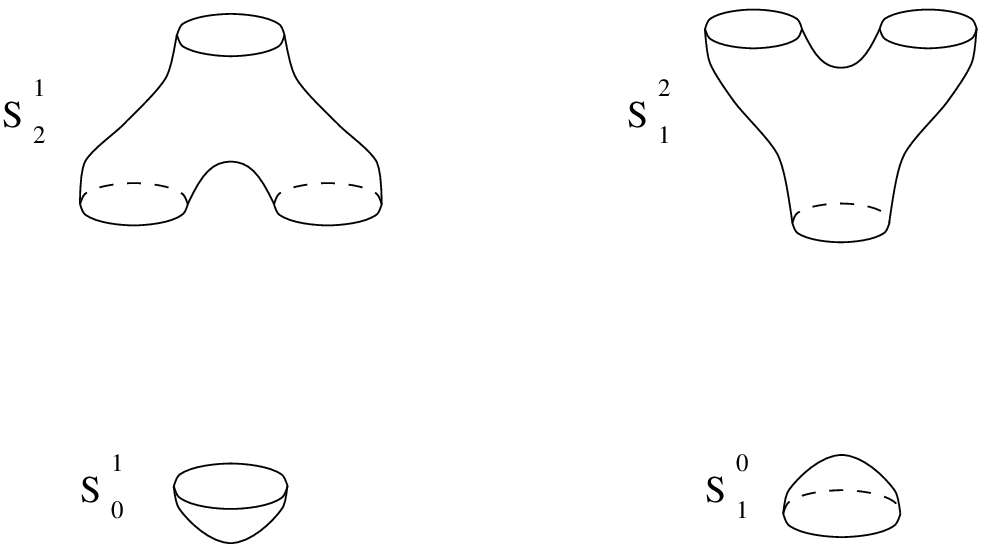}

\vspace{0.3in} 

\drawing{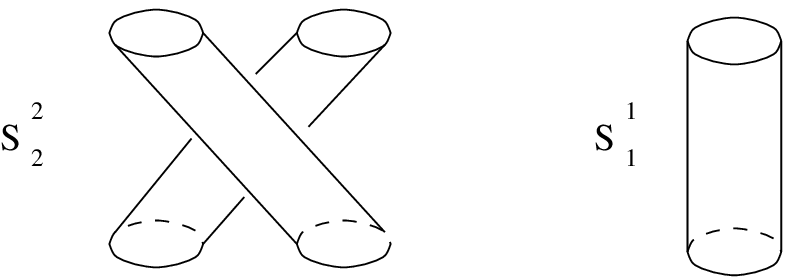} 

Each of these surfaces $S_a^b$ defines a cobordism from a union of $a$ 
circles  to a union of $b$ circles. 
We denote by $\MM$ the category whose objects are closed one-dimensional 
manifolds and morphisms are two-dimensional cobordisms between these
manifolds generated by the above cobordisms.  
Specifically, objects of $\MM$ are enumerated by 
nonnegative integers  $\mbox{Ob}(\MM)=\l \overline{n} | n\in \Z_+\r .$
A morphism between $\overline{n}$ and $\overline{m}$ is a compact
oriented surface $S$ with boundary being the union of $n+m$ circles. 
The 
boundary circles are split into two sets $\partial_0 S$ and $\partial_1 S$ 
with $\partial_0 S$ containing $n$ and $\partial_1 S$ containing 
$m$ circles. An ordering of elements of each of these two sets is fixed.  
The surface $S$ is presented as a concatenation of disjoint unions of 
elementary surfaces, depicted above.
Morphisms are composed in the usual way by gluing boundary 
circles. Two morphisms are equal if the surfaces $S, T$ 
representing these morphisms 
are diffeomorphic via a diffeomorphism that extends the identification 
$\partial_0S \cong \partial_0 T,\partial_1S\cong \partial_1T$ 
of their boundaries.  
$\MM$ is a monoidal category with tensor product of morphisms defined
by taking the disjoint union of surfaces.

Let us construct a monoidal functor from $\MM $ to the category $\Rm$ 
of graded $R$-modules and graded module maps. Assign graded $R$-module  
$A^{\o n}$ to the object $\overline{n}$  and to the 
elementary surfaces $S_2^1,S_1^2,S_0^1,S_1^0,S_2^2,S_1^1$ 
assign morphisms $m, \Delta, \iota, \epsilon, \mbox{Perm},\mbox{Id}$

\begin{equation}
\label{functor-F-defined} 
\begin{array}{lll}  
  F(S_2^1) =m, &  
  F(S_1^2) = \Delta, & 
  F(S_0^1)= \iota, \hspace{0.1in}  \\
  F(S_1^0)=\epsilon, & 
  F(S_2^2)=\mbox{Perm} ,  & 
      F(S_1^1)= \mbox{Id},  \\
\end{array} 
\end{equation}
where $\mbox{Perm}: A\o A\to A\o A$ is the permutation map, 
$\mbox{Perm}(u\o v) = v\o u,$ and $\mbox{Id}$ is the identity 
 map $\mbox{Id}: A \to A.$  

To check that $F$ is well-defined one must verify that for any two 
ways to glue an arbitrary surface $S$ 
in $\Mor(\MM)$ from copies of these six 
elementary surfaces, the two maps of $R$-modules, defined by these 
two decompositions of $S,$ coincide. 
This follows from the commutative algebra and cocommutative coalgebra 
axioms of $A$ and the identity (\ref{m-and-delta}). 

\emph{Remark:} 
Suppose that a surface $S\in \mc M$ contains a 
punctured genus two surface as a subsurface. Then $F(S)$ is 
the zero map. Indeed, we only need to check this 
when $S$ has genus two and one boundary component. Then that $F(S)=0$ 
follows from $m\circ \Delta \circ m \circ \Delta =0.$

Maps $F(S_2^1),F(S_1^2),F(S_0^1),F(S_1^0),F(S_2^2)$ and $F(S_1^1)$  
between tensor powers of $A$ are graded relative to $\ddeg$ with degrees 
\begin{eqnarray*}
\ddeg (F(S_2^1)) & = & \ddeg(F(S_1^2)) =-1, \\ 
\ddeg (F(S_0^1)) & = & \ddeg (F(S_1^0))=1, \\
\ddeg (F(S_2^2)) & = & \ddeg (F(S_1^1)=0   
\end{eqnarray*} 
From this we deduce 

\begin{prop} 
\label{degree-Euler-characteristic} For a surface $S\in \Mor(\MM),$  
the degree of the map $F(S)$ of graded $R$-modules 
 is equal to the Euler characteristic of $S.$
\end{prop}

\subsection{Kauffman bracket} 
\label{Kauffman-bracket} 

In this section we review the Kauffman bracket 
and its relation to the Jones polynomial following Kauffman [Ka].  
Fix an orientation of the 3-space $\R^3$.
A plane projection $D$ of an oriented link $L$ in $\R^3$ is 
called generic if it has no triple intersections, no tangencies and 
 no cusps. In this paper by a plane projection we will mean a 
 generic plane projection. Given a plane projection $D,$ 
we assign a Laurent polynomial $<D>\in \Zq$ to $D$ by the following rules 

\begin{enumerate}
\item A simple closed loop evaluates to $q+q^{-1}:$ 

\drawing{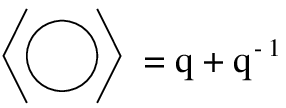}

\item Each over and undercrossing is a linear combination of two
simple resolutions of this crossing: 

   \drawing{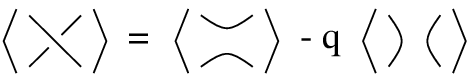} 

\item $< D_1 \bigsqcup D_2> = < D_1> <D_2> $ where 
$< D_1\bigsqcup D_2>$ stands for the disjoint union 
of the diagrams $D_1$ and $D_2.$ 
\end{enumerate} 

From these rules we deduce that 

\drawing{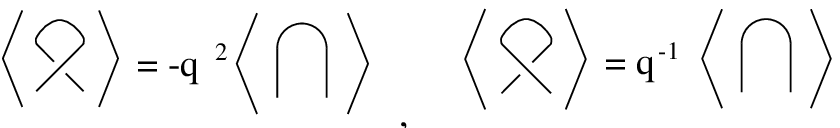}

\drawing{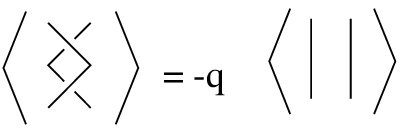}

\drawing{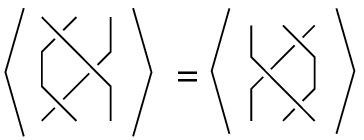}

Curves of the diagram $D$ inherit orientations from that of $L.$ 
Let $x(D)$ be the number of double points in the diagram $D$ that 
look like   
\drawing{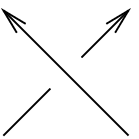}
and $y(D)$ the number of double points that look like 
\drawing{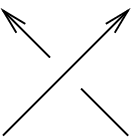}

Then the quantity 
\begin{equation}
\label{KD-defined}  
K(D)=(-1)^{x(D)} q^{y(D)-2x(D)} <D>
\end{equation} 
 does not depend on the choice of a diagram $D$ of the oriented link $L$ 
and is an invariant of $L.$ We denote this invariant by $K(L).$ 
Up to a simple normalization, $K(L)$ 
is the Kauffman bracket of link $L$ and equal to the 
Jones polynomial of $L.$  
The Kauffman bracket, $f[L],$ as defined in [Ka], is a Laurent 
polynomial in an indeterminate $A$ (this $A$ has no relation to the 
algebra $A$ in Section~\ref{here-is-A} of this paper). 
One easily sees that setting our $q$ to $-A^{-2}$ 
and dividing by $(-A^2-A^{-2})$ we get $f[L]$: 
\begin{equation} 
K(L)_{(q=-A^{-2})}= (-A^2-A^{-2})f[L]
\end{equation} 
In this paper we will call $K(L)$ {\it the scaled Kauffman bracket}. 

Let $L_1, L_2$ and $L_3$ be three oriented links that differ as 
shown below. 
\drawing{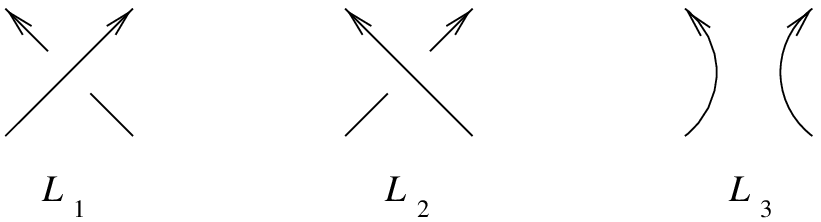}
 The rules for computing the Kauffman bracket imply 
\begin{equation} 
q^{-2}K(L_1) - q^2K(L_2) = (q^{-1} -q )K(L_3)
\end{equation} 
Moreover, $K(L)=q+q^{-1}$ if $L$ is the unknot. 

The Jones polynomial $V(L)$ of an oriented link $L$ is determined 
by two properties: 

\begin{enumerate} 
\item The Jones polynomial of the unknot is $1.$ 
\item For oriented links $L_1,L_2,L_3$ as above
\begin{equation} 
t^{-1}V(L_1) - t V(L_2) = (\sqrt{t}-\frac{1}{\sqrt{t}}) V(L_3)
\end{equation} 
\end{enumerate} 

Therefore, the scaled Kauffman bracket and the Jones polynomial 
are related by 
\begin{equation} 
V(L)_{\sqrt{t}= -q}= \frac{K(L)}{q+q^{-1}}
\end{equation}

\section{Cubes}
\label{section-Cubes} 

\subsection{Complexes of $R$-modules} 
\label{notations-complexes} 

Denote by $\Kom(\CB)$ the category of complexes of 
an abelian category  $\CB.$ An object $N$ of $\Kom(\CB)$ is a collection 
of objects $N^i\in \CB, i\in \Z$ together with morphisms $d^i: N^i\lra 
N^{i+1}, i\in \Z$ such that  $d^{i+1}d^i =0.$ 
A morphism $f: M \to N$ of complexes 
is a collection of morphisms $f^i: M^i \to N^i$ 
such that $f^{i+1}d^i = d^i f^i, i\in \Z.$ 
A morphism $f: M \to N$ is called a quasi-isomorphism if the induced 
map of the cohomology groups $H^i(f): H^i(M) \to H^i(N) $ is an isomorphism
for all $i\in \Z.$

For $n\in \Z$ denote by $[n]$ the automorphism of $\Kom(\CB)$ that 
is defined on objects by $N[n]^i = N^{i+n}, d[n]^i= (-1)^{n} d^{i+n}$
and continued to morphisms in the obvious way.  

The cone of a morphism $f: M \to N$ of complexes is a complex $C(f)$ with 
\begin{eqnarray} 
\label{cone-of-morphism} 
C(f)^i=M[1]^i\oplus N^i, \hspace{0.1in} 
 d_{C(f)}(m^{i+1},n^i)= (-d_M m^{i+1}, 
 f(m^{i+1})+ d_Nn^i).
\end{eqnarray}

The automorphism $\{ n \}$ of shifting the grading down by $n,$ introduced 
  in Section~\ref{here-is-R}, can be naturally extended to an automorphism 
of the category 
$\Kom(\Rm)$ of complexes of graded $R$-modules. This automorphism of 
$\Kom(\Rm)$ will also be denoted $\{ n \}.$ 

To a complex $M$ of graded $R$-modules we associate a graded $R$-module  
$\oplusop{i\in \Z} M^i.$ Each $M^i$ is a graded 
$R$-module, $M^i=\oplusop{j\in \Z} M^i_j$ and thus 
$\oplusop{i\in \Z}M^i$ is a bigraded $R$-module when 
we extend our usual grading of $R$ to a bigrading 
with $c\in R$ having degree $(0,2).$

From this viewpoint the differential $d_M$ of a complex $M$ 
is a homogeneous map of degree $(1,0)$ of bigraded $R$-modules. 

\subsection{Commutative cubes} 
\label{com-cubes} 

  Let $\I$ be a finite set. Denote by $|\I|$ the cardinality of $\I$ 
  and by $r(\I)$ the set 
  of all pairs $(\L,a)$ where $\L$ is a subset of $\I$ and $a$ an element of 
  $\I$ that does not belong to $\L.$ 
  To simplify notation we will often 

  (a) denote a  one-element set $\{ a\}$ by $a,$ 
 
  (b) denote a finite set $\{ a,b, \dots , d\}$ by $ab\dots d,$ 

  (c) denote the disjoint union $\L_1 \sqcup \L_2$ of two sets 
  $\L_1, \L_2$ by $\L_1 \L_2,$ in particular, we denote 
  by $\L a$ the disjoint union of a set $\L$ and a one-element set $\{ a\},$ 
  similarly, $\L  ab$ means $\L \sqcup \{ a\} \sqcup \{ b\},$ etc.  

\begin{definition} 
\label{com-cubes-definition} Let $\I$ be a finite set and $\CB$ a 
category. A commutative $\I$-cube $V$ over $\CB$ is a
collection of objects $V(\L)\in Ob(\CB)$ for each subset $\L$ of $\I,$ 
morphisms 
\begin{equation} 
\xi^{V}_a(\L): V(\L) \lra V(\L  a)
\end{equation} 
for each $(\L,a)\in r(\L),$ 
such that for each triple $(\L,a,b),$ where $\L$ is a subset of $\I$ and 
 $a,b,a\not= b$ are two elements of $\I$ that do not lie in $\L,$ 
 there is an equality of morphisms 
\begin{equation} 
\label{comcube-equation}
\xi^{V}_b(\L  a)
\xi^{V}_a(\L) = 
\xi^{V}_a(\L  b)
\xi^{V}_b(\L), 
\end{equation}  
i.e., the following diagram is commutative 
\[
\begin{CD}
\label{comcube-diagram}
V(\L) @>{\xi^V_a(\L)}>> V(\L  a) \\
@VV{\xi^V_b(\L)}V          @VV{\xi^V_b(\L   a )}V \\
V(\L   b)  @>{\xi_a^V(\L  b)}>>      V(\L   a   b) 
\end{CD}
\]

\end{definition} 

We will call a commutative $\I$-cube an \emph{$\I$-cube} or, 
sometimes, a \emph{cube} when it is clear what $\I$ is.  
Maps $\xi_a^V$ are called the structure maps of $V.$

{\it Example } If $\I$ is the empty set, an $\I$-cube is an object 
 in $\CB.$ If $\I$ consists of one element, an $\I$-cube is a morphism 
 in $\CB.$ If $\I$ consists of two elements, $\I= \{ a,b\},$ an $\I$-cube 
  is a commutative square of objects and morphisms in $\CB$: 

\[
\begin{CD}
V(\emptyset)     @>>>       V(a)   \\
@VVV                   @VVV \\
V(b)    @>>>       V(a  b)  
\end{CD}
\]

In general, an $\I$-cube can be visualized in the following manner. 
Let $n$ be the cardinality of $\I.$ We take an $n$-dimensional cube 
in standard position in 
the Euclidean $n$-dimensional space, i.e. each vertex has coordinates 
$(a_1, \dots , a_n)$ where  $a_i \in \{ 0,1\}$ 
 and orient each edge in the direction of the vertex with the bigger 
sum of the coordinates. Then the edges of 
any 2-dimensional facet of this cube are oriented 
as shown on the diagram below.

\drawing{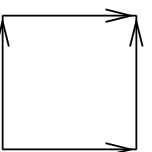} 

Choose a bijection between elements of $\I$ and coordinates of 
 $\R^n.$ This bijection defines a bijection between vertices of the 
 $n$-cube and subsets of $\I,$ with the $(0, \dots, 0)$ vertex associated to 
 the empty set. Oriented edges of the $n$-cube correspond to 
 pairs $(\L,a)\in r(\I).$ 

 Given an $\I$-cube $V,$ put object $V(\L)$ into the vertex 
 associated to the set $\L$ and assign morphism $\xi^V_a(\L)$ to the 
 arrow going from the vertex associated to $\L$ to the vertex associated to 
 $\L  a.$  Equation (\ref{comcube-equation}) is equivalent to 
the commutativity of diagrams in all 2-dimensional faces of the $n$-cube.

Given two $\I$-cubes $V,W$ over a category $\CB,$ an $\I$-cube 
map $\psi:V\lra W$ is a collection of maps 
$$\psi(\L): V(\L)\lra W(\L), \hspace{0.3in}\mbox{ for all } \L\subset \I$$
that make diagrams 
\begin{equation} 
\label{comcube-map} 
\begin{CD}
V(\L) @>{\psi(\L)}>> W(\L) \\
@VV{\xi^V_a(\L)}V          @VV{\xi^W_a(\L)}V \\
V(\L  a)  @>{\psi(\L  a)}>>      W(\L  a) 
\end{CD}
\end{equation} 
commutative for all $(\L,a)\in r(\I).$
The map $\psi$ is called an isomorphism if $\psi(\L)$ is an isomorphism 
for all $\L\subset \I.$ The map $\psi$ of $\I$-cubes over an 
abelian category $\CB$ is called injective/surjective 
if $\psi(\L)$ is injective/surjective for all $\L\subset \I.$ 

The class of $\I$-cubes over an abelian category $\CB$ and maps between 
$\I$-cubes constitute an abelian category in the obvious way. 
In particular, direct sums of $\I$-cubes are defined. 

 For a finite set $\I$ and $a\in \I,$ let $\J$ be the 
 complement, $\I= \J \sqcup \{ a\}.$ Given an $\I$-cube 
  $V,$ let $V_a(\ast 0), V_a(\ast 1)$ be $\J$-cubes defined as 
  follows: 
  \begin{equation} 
  V_a(\ast 0) (\L) = V(\L), \hspace{0.2in} 
  V_a(\ast 1) (\L) = V(\L   a), \hspace{0.2in} 
  \mbox{ for } \L\subset \J
  \end{equation} 
  and the structure maps of $V_a(\ast 0), V_a(\ast 1)$ are determined by 
  the structure maps $\xi_b^V, b\in \I_1$ of $V$ in the obvious fashion. 
  Sometimes we will write $V(\ast 0)$ for $V_a(\ast 0),$ etc. 
The structure map $\xi_a^V$ of $V$ defines an $\I_1$-cube map 
$\xi_a^V: V_a(\ast 0)\lra V_a(\ast 1).$ 
This provides a one-to-one correspondence between $\I$-cubes and 
maps of $\I_1$-cubes.

We say that a  map $\psi: V\to W$ of $\I$-cubes over the category 
$R$-mod of graded $R$-modules and graded maps is graded of degree $i$ if 
the map $\psi(\L): V(\L)\to W(\L)$ has degree $i$ for all $\L\subset\I.$

For a cube $V$ over $\Rm$ denote by $V\{ i\}$ the cube $V$ with 
the grading shifted by $i$:  
\begin{equation} 
V\{ i\} (\L) = V(\L) \{ i\} \mbox{ for all } \L \subset \I
\end{equation} 
and the structure maps being appropriate shifts of the structure maps of $V.$ 
A degree $i$ map $\psi: V\to W$ of $\I$-cubes over $\Rm$ induces a 
grading-preserving map $V\to W\{ i \},$ also denoted $\psi.$ 

\subsection{skew-commutative cubes} 
\label{section-cubes} 

We next define skew-commutative $\I$-cubes over an additive 
 category $\CB.$ A skew-commutative  $\I$-cube is almost the 
same as a commutative 
 $\I$-cube but now we require that for every square facet 
of the cube the associated diagram of objects and morphisms of $\CB$ 
anticommutes.

\begin{definition}
\label{cubes-definition} 
 Let $\I$ be a finite set and $\CB$ an additive category. 
A skew-commutative $\I$-cube $V$ over $\CB$ is a
collection of objects $V(\L)\in Ob(\CB)$ for $\L\subset\I$ and 
morphisms 
\begin{equation*}
\xi^V_a(\L): V(\L) \lra V(\L  a).
\end{equation*}
such that for each triple $(\L,a,b),$ where $\L$ is a subset of $\I$ and 
  $a,b,a\not= b$ are two elements of $\I$ that do not lie in $\L,$
 there is an equality 
\begin{equation*}
\xi^{V}_b(\L  a)
\xi^{V}_a(\L) +  
\xi^{V}_a(\L  b)
\xi^{V}_b(\L)=0.  
\end{equation*}

\end{definition} 

We will call a skew-commutative $\I$-cube over $\CB$ a 
\emph{skew $\I$-cube} or, without specifying $\I$, a \emph{skew cube}. 

Given $\I$-cubes or skew $\I$-cubes $V$ and $W$
over $\Rm$, their tensor product is defined 
to be an $\I$-cube (if $V$ and $W$ are both cubes or both skew cubes) 
or a skew $\I$-cube (if one of $V, W$ is a cube and the other is a
skew cube), denoted $V\o W,$ given by 
\begin{eqnarray*} 
& & (V\o W)(\L) = V(\L)\o W(\L), \hspace{0.3in} \L \subset\I, \\
& & \xi_a^{V\o W}(\L)=\xi_a^{V}(\L)\o \xi_a^{W}(\L),
\hspace{0.3in} (\L,a)\in r(\I),  
\end{eqnarray*} 
where, recall, the tensor products are taken over $R.$ 

For a finite set $\L$ denote by $o(\L)$ the set of complete orderings 
  or elements of $\L.$ For $x,y\in o(\L)$ let $p(x,y)$ be the parity 
 function, $p(x,y)=0$ if $y$ can be obtained by from $x$ via an even number 
  of transpositions of two neighboring elements in the ordering, 
 otherwise, $p(x,y)=1.$ 
 To a finite set $\L$ associate a graded $R$-module $E(\L)$ defined 
  as the quotient of the graded $R$-module, freely generated by elements 
  $x$ for all $x\in o(\L),$ by relations $x = (-1)^{p(x,y)}y $ for 
  all pairs $x,y\in o(\L).$ Module $E(\L)$ is a free graded $R$-module 
  of rank $1.$ 
  For $a\not\in \L$ there is a canonical isomorphism of graded $R$-modules 
  $E(\L) \lra E(\L   a)$ induced by the map  $o(L)\to o(L  a)$
  that takes $x\in o(L)$ to $xa\in o(L  a).$ Moreover, for 
  $a,b, a\not= b,$ the diagram below anticommutes 

   \begin{equation} 
  \begin{CD}
  E(\L) @>>> E(\L  a) \\
  @VVV          @VVV \\
  E(\L  b)  @>>>      E(\L  a  b) 
  \end{CD}
  \end{equation} 
 
  Denote by $E_{\I}$ the skew $\I$-cube with $E_{\I}(\L)= E(\L)$ for 
  $\L \subset \I$ and the structure map $E_{\I}(\L)\to E_{\I}(\L   a)$ 
  being canonical isomorphism $E(\L) \to E(\L   a).$ 

  We will use $E_{\I}$ to pass 
  from $\I$-cubes over $\Rm$ to skew $\I$-cubes 
  over $\Rm$ by tensoring an $\I$-cube with $E_{\I}.$

\subsection{skew-commutative cubes and complexes}
\label{section-cubes-complexes} 

Let $V$ be a skew $\I$-cube over an abelian category $\CB.$ 
To $V$ we associate a complex $\C(V)=(\C^i(V),d^i), i\in \Z$ 
of objects of $\CB$ by 
\begin{equation} 
\C^i(V)=\oplusop{\L\subset \I, |\L|=i} V(\L)
\end{equation} 
The differential $d^i:\C^i(V)\to \C^{i+1}(V)$ is given on an element 
$x\in V(\L), |\L|=i$ by 
\begin{equation} 
d^i (x)=\sum_{a\in \I\setminus \L} \xi^V_a(\L)x.
\end{equation} 

\emph{Examples: } 
\begin{enumerate} 
\item  If $|\I|=1, \I= \{ a\},$ 
\[ \C^i(V) = 
\left\{ 
\begin{array}{ll}
V(\emptyset)  & \mbox{ if $i=0 $ } \\
V(a)  & \mbox{ if $i=1 $ } \\
0     & \mbox{ otherwise }
\end{array} 
\right.  \] 
The differential $d^0= \xi_a^V(\emptyset)$ and  $d^i=0$ if $i\not= 0,$
so $\C(V)$ is the complex 
\begin{equation} 
 \cdots \lra 0 \lra V(\emptyset) \stackrel{\xi^V_a(\emptyset)}{\lra} V(\I) 
\lra 0 \lra \cdots 
\end{equation} 

\item If $\I$ contains two elements, say, $\I= \{ a,b\},$ then   
\[ \C^i(V) = 
\left\{ 
\begin{array}{ll}
V(\emptyset)  & \mbox{ if $i=0 $ } \\
V(a)\oplus V(b)  & \mbox{ if $i=1 $ } \\
V(a  b)              & \mbox{ if $i=2 $ } \\
0     & \mbox{ otherwise }
\end{array} 
\right.  \] 
and differentials 
\begin{eqnarray*} 
& d^0: & V(\emptyset) \lra V(a)\oplus V(b) \\
& d^0=  & \xi_b^V(\emptyset)+ \xi_a^V(\emptyset)  \\
& d^1: & V(b) \oplus V(a) \lra V(a  b) \\
& d^1= & (\xi_a^V(b) , \xi_b^V(a) )
\end{eqnarray*}
\end{enumerate}

\begin{prop}
\label{isomorphism-of-cubes}  
Let $V$ be a skew 
 $\I$-cube over an abelian category $\CB$ and suppose 
that for some $a\in \I$ and any $\L\subset \I\setminus\{ a\} $ the map 
$\xi^V_a: V(\L) \to V(\L  a)$ is an isomorphism. 
Then the complex $\C(V)$ is acyclic.  
\end{prop} 

\emph{Proof: }The complex $\C(V)$ is isomorphic to the cone of 
the identity map of the complex $\C(V_a(\ast 1))[-1]$ and, therefore, 
acyclic. $\square$ 

Every map of $\I$-cubes $\phi: V\lra W$ over $\Rm$ induces 
a map of complexes 
\begin{equation} 
\C(\phi) : \C(V\o E_{\I}) \lra \C(W\o E_{\I}).
\end{equation} 
If $\phi$ is an isomorphism of commutative cubes, 
$\C(\phi)$ is an isomorphism of complexes. 

\begin{prop} 
Let $V$ be an $\I$-cube over $\Rm$ and 
suppose that for some $a\in \I$ the structure map 
\begin{equation} 
\xi^V_a: V_a( \ast 0)\lra V_a(\ast 1) 
\end{equation} 
is an isomorphism. Then the complex $\C(V\o E_{\I})$ is acyclic. 
\end{prop} 

\emph{Proof: } Immediate from Proposition~\ref{isomorphism-of-cubes}.
$\square$  

The following proposition and its corollary are obvious. 

\begin{prop} We have a canonical splitting of complexes
\begin{equation} 
\C(V\oplus W) = \C(V) \oplus \C(W) 
\end{equation} 
where $V$ and $W$ are skew-commutative 
$\I$-cubes over an abelian category and 
$V\oplus W$ is the direct sum of $V$ and $W$. 
\end{prop} 

\begin{corollary} 
We have a canonical splitting of complexes 
\begin{equation} 
\C((V\oplus W)\o E_{\I}) = \C(V\o E_{\I}) \oplus \C(W\o E_{\I}) 
\end{equation} 
where $V$ and $W$ are $\I$-cubes over $\Rm.$
\end{corollary}

\section{Diagrams} 
\label{section-diagrams} 

\subsection{Reidemeister moves} 
\label{Reidemeister-moves} 

Given a link $L$ in $\R^3$ we can take its generic projection on the 
plane. A generic projection is the one
without triple points and double tangencies. An isotopy class of 
such projections is called a \emph{plane diagram} of $L,$
or, simply, a \emph{diagram}. 
The following four types of transformations of plane diagrams 
preserve the isotopy type of the associated link. 

I. Addition/removal  of a left-twisted curl: 

\drawing{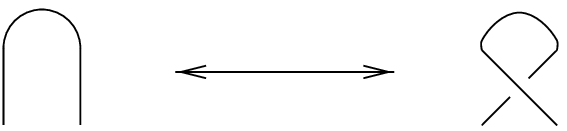}

II. Addition/removal of a right-twisted curl: 

\drawing{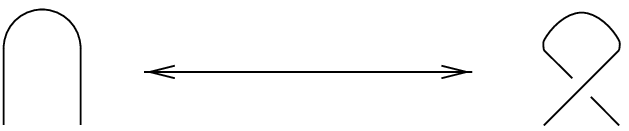} 

III. Tangency move: 

\drawing{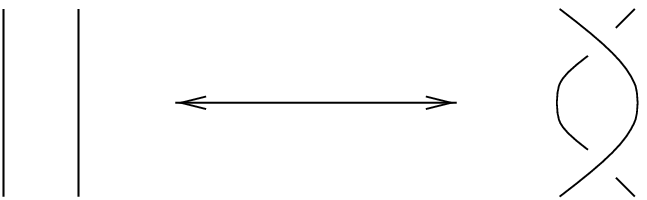} 

IV. Triple point move: 

\drawing{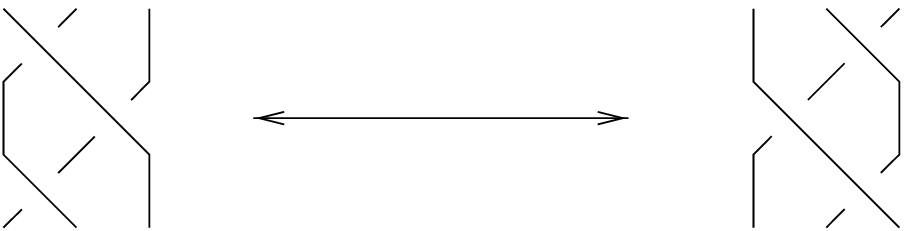} 

\begin{prop}
\label{five-moves}
 If plane diagrams $D_1$ and $D_2$ represent 
isotopic oriented links, these diagrams can be connected by a 
chain of moves I-IV. 
\end{prop} 

$\square$ 

\subsection{Constructing cubes and complexes from plane diagrams} 
\label{diagrams-to-comcubes} 

Fix a plane diagram $D$ with $n$ double points of an oriented link $L.$
 Denote by $\I$ the set of double points of $D.$ 
To $D$ we will associate an $\I$-cube $V_D$
over the category $\Rm$ of graded $R$-modules.  
This cube will not depend on the orientation of components of $L.$ 

Given a double point of a diagram $D$, it can be resolved
in two possible ways: 

\drawing{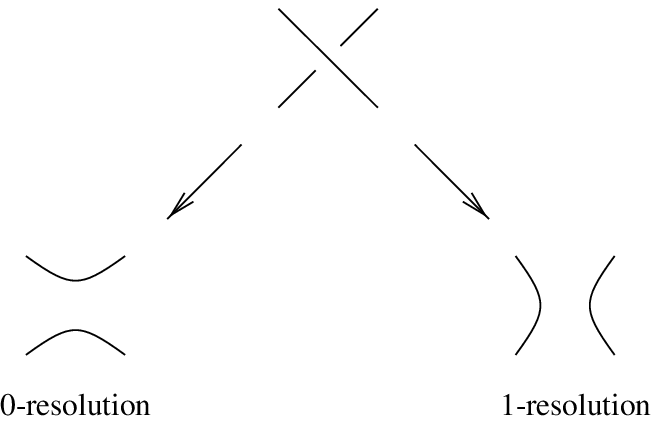} 

Let us call the resolution on the left the $0$-resolution, and the one 
on the right the $1$-resolution. A {\it resolution} of $D$ is a 
 resolution of each double point of $D.$ Thus, $D$ admits $2^n$ 
 resolutions. There is a one-to-one correspondence between 
 resolutions of $D$  and subsets $\L$ of the set $\I$ of double 
 points. Namely, to $\L \subset \I$  we associate a resolution, 
 denoted $D(\L),$ by taking $1$-resolution of each double point that 
 belongs to $\L$ and $0$-resolution if the double point does not lie in 
 $\L.$ 

 A resolution of a diagram $D$ is always a collection of simple 
 disjoint curves on the plane and is thus a $1$-manifold embedded in the 
 plane. Now the functor $F$ from $(1+1)$ cobordisms to 
 $R$-modules (see Section~\ref{A-and-surfaces}) comes into play.
 To a union of $k$ circles it assigns the $k$-th tensor
 power of $A.$ The functor $F$, applied to the diagram $D(\L),$ considered 
 as a one-dimensional manifold, produces a graded $R$-module $A^{\o k}$ 
 where $k$ is the number of components of $D(\L).$ We raise the 
grading of $A^{\o k}$ by $|\L|,$ the cardinality of $\L,$  and assign 
the $R$-module $F(D(a))\{- |\L|\}$ to the vertex $V_D(\L)$ 
of the cube $V_D:$ 
\begin{equation} 
\label{here-is-VD} 
V_D(\L) = F(D(\L))\{- |\L|\}
\end{equation} 
(Recall from Section~\ref{notations-complexes}
 that the automorphism $\{ 1\}$ of the category $\Rm$ lowers 
the grading by $1$.) 
Let us now define maps between vertices of $V_D.$ 
Choose $(\L,a)\in r(\I).$ 
We want to have a map 
\begin{equation} 
\xi_a^{V_D}(\L): V_D(\L)\lra V_D(\L  a).
\end{equation}  
The diagrams $D(\L)$ and $D(\L  a)$ differ only in the neighborhood 
of the double point $a$ of $D,$ as an example below (for $n=1$, so 
that $D$ has one double point, $\I= \{ a\}$) demonstrates
($D$ is the leftmost diagram, $D(\emptyset)$ is the diagram in the center 
 and $D(a)$ is depicted on the right):

\drawing{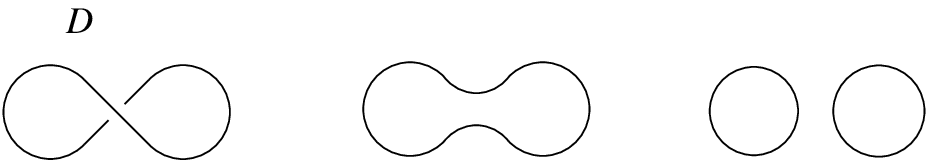}

Take the direct product 
of the plane $\R^2$ and the interval $[0,1].$ We identify the 
diagram $D(\L)$ (resp. $D(\L  a)$) 
with a one-dimensional submanifold of $\R^2\times \{ 0 \}$ 
(resp $\R^2\times \{ 1\}$.) 
We can choose a small neighbourhood $U$ of $a$ such that 
$D(\L)$ and $D(\L  a)$ coincide outside $U$ and inside 
they look as follows: 

\drawing{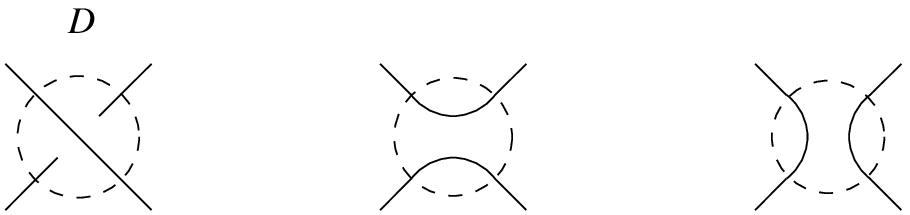} 

The boundary of $U$ is depicted by a dashed circle, the central 
 picture shows the intersection of $D(\L)$ and $U,$ the 
  rightmost picture shows the intersection of $D(\L   a)$ and $U.$  
Let $S$ be a surface properly embedded in $\R^2\times [0,1]$ such 
that 

\begin{enumerate} 
\item The boundary of $S$ is the union of the diagrams $D(\L)$ and 
$D(\L  a).$
\item Outside of $U\times [0,1]$ surface $S$ is the direct product 
of $D(\L)\cap (\R^2\setminus U)$ and the interval $[0,1].$ 
\item The connected component of $S$ that has a nonempty intersection 
with $U\times [0,1]$ is homeomorpic to the two-sphere with 
three holes.  
\item The projection $S\lra [0,1]$ onto the second component of 
the product $\R^2\times [0,1]$ has only one critical point -- the 
saddle point that lies inside $U\times [0,1].$  
\end{enumerate}

Example: Let 
$D$ be a diagram with two double points, $\I=\{ a,b\},$ 
 diagram $D(a),$ respectively $D(a b),$ 
consists of 2, respectively 3 simple curves:

\drawing{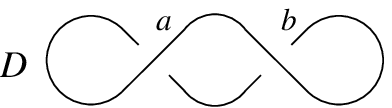}

\drawing{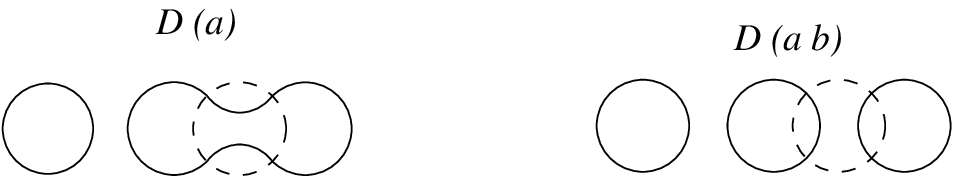} 

The boundary of the neighbourhood $U$ of the double point $a$ is 
depicted by the dashed circle on the diagram above. 
Then the surface $S$ looks like 

\drawing{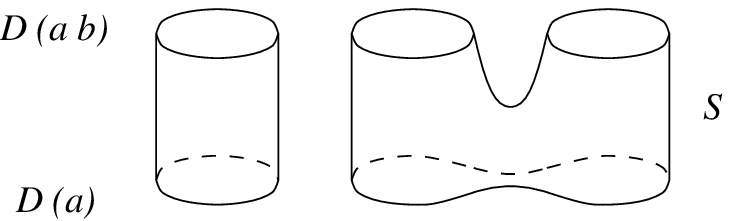} 

Recall that earlier we defined $V_D(\L)$ to be $F(D(\L))$ for 
$\L\subset \I,$ with the degree raised by $|\L|.$ Now define the map 
$$\xi_a^{V_D}(\L): V_D(\L)\lra V_D(\L  a)$$
to be given by 
$$F(S): F(D(\L)) \lra F(D(\L  a)).$$
Note that the degree of $F(S)$ is equal to $-1,$ the Euler characteristic 
of the surface $S$ 
(Proposition~\ref{degree-Euler-characteristic}). 
But $|\L  a|= |\L|+1,$ 
so, 
with degrees shifted: 
\begin{eqnarray}
V_D(\L) & = & F(D(\L))\{- |\L| \} \\
V_D(\L  a) & = & F(D(\L  a))\{- |\L|-1 \} 
\end{eqnarray} 
  and the map $\xi_a^{V_D}(\L)$ is a grading-preserving map 
of graded $R$-modules.  

\begin{prop} $V_D$, defined in this way, is an $\I$-cube over the 
category $\Rm$ of graded $R$-modules and grading-preserving maps.  
\end{prop}

The proof consists of verifying  
commutativity relations (\ref{comcube-equation}) for maps 
$\xi_a^{V_D}(\L).$  They follow immediately from the functoriality of $F.$ 
 
$\square$

\emph{Example:} Let $D$ be the diagram 

\drawing{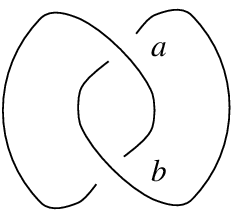} 

The four resolutions of this diagram are given below 

\drawing{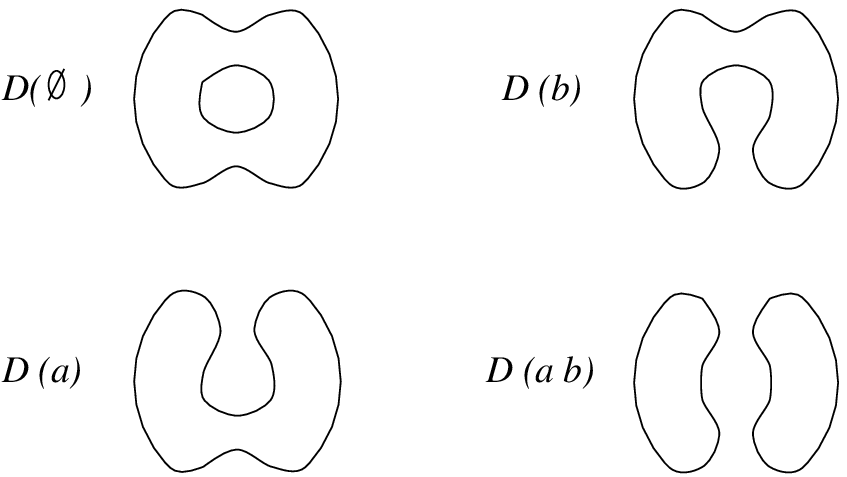} 

Applying the functor $F$, we get 
\begin{eqnarray*} 
& & F(D(\emptyset))=A^{\o 2}, \hspace{0.3in} 
F(D(a))= A, \\ 
& & F(D(b))= A, \hspace{0.3in}
F(D(a  b))= A^{\o 2}.  
\end{eqnarray*} 
The cube $V_D$ has the form 
\[
\begin{CD}
A^{\o 2} @>{m}>>   A\{ -1 \} \\
@VV{m}V        @VV{\Delta}V  \\
A\{ -1 \}     @>{\Delta}>>  A^{\o 2}\{ -2 \}  
\end{CD} 
\]

$\square$

 Let us now go back to our construction. So far, to a plane diagram 
 diagram $D$ with the set $\I$ of double points 
 we associated a  $\I$-cube  
 $V_D$ over the category $\Rm$ of graded $R$-modules. 
 We would like to build a complex of graded $R$-modules out of $V_D.$ 
 We know how to build a complex from a skew-commutative  $\I$-cube (see 
 Section~\ref{section-cubes-complexes}). To 
 make a skew-commutative 
  $\I$-cube out of an $\I$-cube $V_D$ we put minus sign in front 
of some structure maps $\xi^{V_D}$ of $V_D$ so that for any commutative 
square 
of $V_D$ an odd number out of the four maps constituting the square change 
signs. A more intrinsic way to do this is to tensor
$V_D$ with the skew-commutative $\I$-cube $E_{\I},$ defined at the end of 
 Section~\ref{section-cubes}. 

To the skew-commutative 
$\I$-cube $V_D\o E_{\I}$ there is associated the  complex 
$\C(V_D\o E_{\I})$ of graded $R$-modules (see  
Section~\ref{section-cubes-complexes}). Denote this 
complex by $\C(D):$ 
\begin{equation}
\label{bar-CD-defined}  
\C(D) \stackrel{\mbox{\scriptsize def}}{=} \C(V_D\o E_{\I}) 
\end{equation}

Thus, $\C(D)$ is a complex of graded $R$-modules 
and grading-preserving homomorphisms. It does not depend on 
the orientations of the components of the link $L.$ 

Recall (Section~\ref{notations-complexes}) 
that the category $\Kom(\Rm)$ has two commuting 
automorphisms: $[1],$  which shifts a complex 
one term to the left; and $\{ 1\},$ which  
lowers the grading of each component of the complex by $1.$ 

To the  diagram $D$ of link $L$ we associated 
(Section~\ref{Kauffman-bracket}) two 
numbers, $x(D)$ and $y(D).$ 
Define a complex $C(D)$ by 
\begin{equation}
\label{CD-defined}  
C(D)=\C(D)[x(D)]\{2 x(D)-y(D)\}
\end{equation}

Define $H^i(D)$ as the $i$-th cohomology group of $C(D).$ It is 
a finitely-generated graded $R$-module.  

\begin{theorem} 
\label{closed-i} 
If $D$ is a plane diagram of an oriented link $L$, then for each $i\in \Z,$ 
the isomorphism class of the graded $R$-modules $H^i(D)$ is an invariant 
of $L.$   
\end{theorem} 

The proof of this theorem 
occupies Section~\ref{section-transformations}, together 
with some preliminary material contained in 
Section~\ref{surfaces-comcube-morphisms}.

Define $H^{i,j}(D)$ as the $i$-th cohomology group of the degree $j$ 
subcomplex of $C(D).$ Thus, $H^{i,j}(D)$ is the graded component of 
$H^i(D)$ of degree $j$ and we have a decomposition of abelian groups 
\begin{equation}
H^i(D)= 
\oplusop{j\in \Z}H^{i,j}(D).
\end{equation}

We denote by $H^i(L)$ the isomorphism class of $H^i(D)$ in the 
category of graded $R$-modules. For an oriented link $L$ only finitely 
many of $H^i(L)$ are non-zero as $i$ varies over all integers. 

\begin{corollary}
\label{closed-i-and-j} 
 If $D$ a plane diagram of an oriented link $L$, then for each $i,j\in \Z,$ 
the isomorphism class of the abelian group $H^{i,j}(D)$
is an  invariant of $L.$  
\end{corollary} 

We next show that the Kauffman bracket is equal to a 
suitable Euler characteristic of these cohomology groups.

\begin{prop} For an oriented link $L,$ 
\begin{equation} 
\label{first-characteristic} 
K(L)= (1-q^2)\sum_{i\in \Z} (-1)^{i}\whc(H^i(D))
\end{equation}
where $K(L)$ is the scaled Kauffman bracket, defined in 
Section~\ref{Kauffman-bracket}, $\whc$ is the Euler characteristic 
  (see Section~\ref{here-is-R}) and $D$ is any diagram of $L.$ 
\end{prop} 

{\emph Proof:} First notice that 
$\whc(M\{ n\} ) = q^{-n} \whc(M)$ for a finitely-generated graded $R$-module 
$M.$ Given a bounded complex 
\begin{equation} 
M:\hspace{0.3in} \dots \to M^i \to M^{i+1}\to \dots 
\end{equation} 
 of finitely-generated graded $R$-modules, define 
\begin{equation} 
\whc(M) = \sum_{i\in \Z} (-1)^i \whc(M^i)
\end{equation} 
Since
\begin{equation} 
\whc(C(D))= \sum_{i\in \Z} (-1)^i \whc(H^i(D)), 
\end{equation} 
it is enough to prove 
\begin{equation} 
K(L) = (1-q^2)\whc(C(D))
\end{equation} 
For three diagrams $D_1,D_2$ and $D_3$ that differ as shown below
\drawing{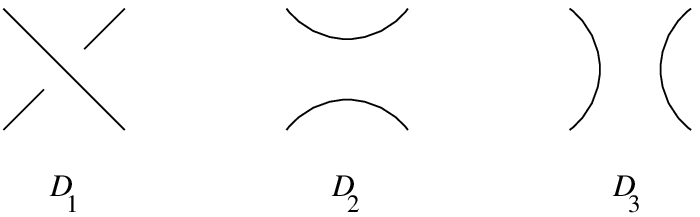} 
the complex $\C(D_1)[1]$ is isomorphic, up to a shift,  to 
the cone of a map of complexes $\C(D_2)\to \C(D_3)\{ -1 \} .$ 
Therefore, 
\begin{equation} 
\whc(\C(D_1)) = \whc(\C(D_2)) - \whc(\C(D_3)\{ -1\} ) = 
  \whc(\C(D_2)) - q \whc(\C(D_3))
\end{equation} 
On the other hand, for diagrams $D_1,D_2,D_3$ as above, we have  
\begin{equation} 
<D_1> =  <D_2> - q < D_3>
\end{equation}  
(see Section~\ref{Kauffman-bracket}, where $<D>$ is defined). 
If the diagram $D$ is a disjoint union of $k$ 
 simple plane curves then  
\begin{equation} 
\whc(\C(D))= \whc(A^{\o k})=
(q+q^{-1})^k\whc(R)= \frac{(q+q^{-1})^k}{1-q^2}
\end{equation} 
and $<D> = (q+q^{-1})^k.$ Therefore, for any diagram $D$ 
\begin{equation} 
<D> = (1-q^2) \whc(\C(D)). 
\end{equation} 
Since 
\begin{equation} 
\whc(C(D))= \whc(\C(D))[x(D)]\{ 2x(D)-y(D)\} = 
(-1)^{x(D)}q^{y(D)-2x(D)} \whc(\C(D)),  
\end{equation} 
and in view of (\ref{KD-defined}), proposition follows. $\square$

\subsection{Surfaces and cube morphisms} 
\label{surfaces-comcube-morphisms} 

Let $U$ be a closed disk in the plane $\R^2$ and $\dot{U}$ the 
interior of $U$ so that $U= \partial U \cup \dot{U}.$ Let $T'$ be 
a tangle in $(\R^2\setminus \dot{U})\times [0,1]$ with $m$ points 
(where $m$ is even) on the boundary $\partial U \times [0,1]$ and 
$T$ a generic projection of $T'$ on $\R^2 \setminus \dot{U}.$ 
The intersection of $T$ with $\partial U$ consists of $m$ points. 
Denote them by $p_1,\dots ,p_m$ (see an example on the diagram below, 
$\partial U$ is shown by a dashed circle). 

\drawing{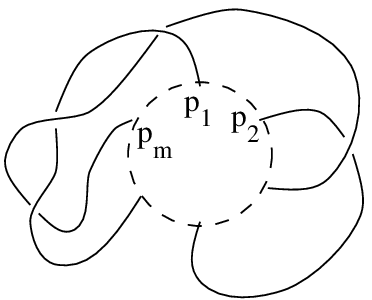}

Let $\I$ be the set of double points of $T.$ Pick two systems $Q_0$ and $Q_1$ 
of $\frac{m}{2}$ simple disjoints arcs in $U$ with ends in points 
$p_1,\dots p_m$: 

\drawing{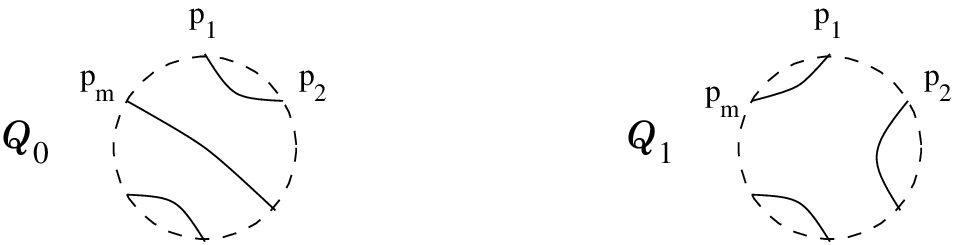} 

Then $Q_0\cup T$ and $Q_1\cup T$ (here and further on we denote them 
by $P_0$ and $P_1$ respectively) can be considered as two plane 
diagrams of links in $\R^3:$ 

\drawing{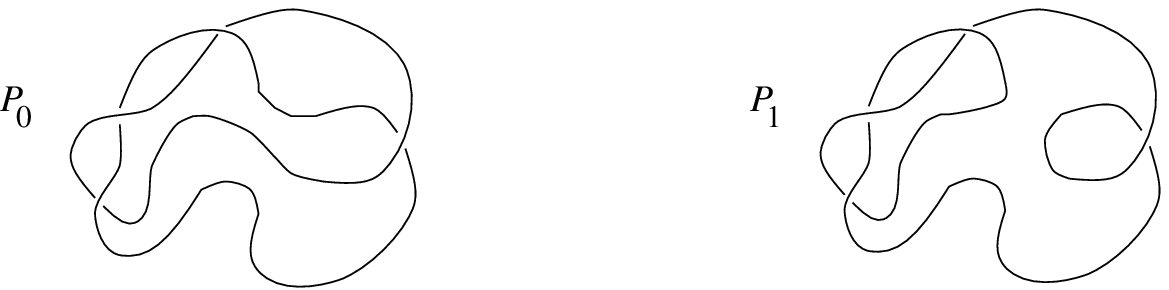} 

To $P_0$ and $P_1$ there are associated $\I$-cubes 
$V_{P_0}$ and $V_{P_1}.$ 

Let $S$ be a compact oriented surface in $U\times [0,1]$ such that 
the boundary of $S$ is the union of $Q_0\times \l 0\r , 
Q_1\times \l 1\r $ and $(p_1\cup \dots \cup p_m)\times [0,1].$ 
To $S$ we associate an $\I$-cube map
$$\psi_S: V_{P_0} \lra V_{P_1}$$
as follows. For each $\L\subset \I$ we must construct a map 
\begin{equation}
\psi_{S,\L}: V_{P_0}(\L)\{-|\L|\} \lra V_{P_1}(\L)\{ -|\L| \}
\end{equation}
and check the commutativity of diagrams (\ref{comcube-map}). 

To $\L$ there is associated a resolution $T(\L)$ of double points 
of $T.$ Thus $T(\L)$ is a collection of simple closed curves and arcs in 
$\R^2 \setminus \mbox{\.{U}}$ with ends in $p_1,\dots ,p_m.$ 
Then (by (\ref{here-is-VD})
\begin{eqnarray}
V_{P_0}(\L) & = & F(T(\L)\cup Q_0)\{ - |\L|\}  \\
V_{P_1}(\L) & = & F(T(\L)\cup Q_1)\{ -|\L|\} 
\end{eqnarray} 
where $F$ is the functor described in Section~\ref{A-and-surfaces} 
( $T(\L)\cup Q_0$ and $ T(\L)\cup Q_1$ 
are collections of simple closed curves on the plane,  so that 
we can apply functor $F$ to them). 

Let $S'$ be a surface in $\R^2\times [0,1]$ which is $S$ inside 
$U\times [0,1]$ and $T(\L)\times [0,1]$ outside $U\times [0,1].$ 
Map $F(S')$ 
\begin{equation}
F(S'): F(T(\L)\cup Q_0)  \lra F(T(\L) \cup Q_1) 
\end{equation} 
is a graded map of $R$-modules of degree 
$\chi(S')= \chi(S)-\frac{m}{2}.$ 
Define $\psi_{S,\L}$ as this map, shifted by $|\L|$: 
\begin{equation} 
\psi_{S,\L}=F(S')\{ -|\L| \} : 
 V_{P_0}(\L)\{ -|\L|\}  \lra V_{P_1}(\L) \{ -\L|\} 
\end{equation} 
The commutativity condition (\ref{comcube-map}) is immediate. 
We sum up our result as 

\begin{prop} The map 
\begin{equation} 
\psi_S: V_{P_0} \lra V_{P_1}
\end{equation}
is a degree $\chi(S)-\frac{m}{2}$ map of $\I$-cubes. 
\end{prop} 

Everything in this section extends to the case when the  
diagrams $Q_0$ and $Q_1$ are allowed to have simple closed circles in 
addition to $\frac{m}{2}$ simple disjoint acts 
joining points $p_1,\dots p_m.$ For instance, $Q_0$ 
may look like 

\drawing{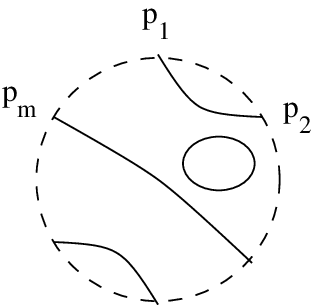}

In this more general case to each compact oriented surface 
$S$ in $U\times [0,1]$ such that 
the boundary of $S$ is the union of $Q_0\times \{ 0\} , 
Q_1\times \{ 1\} $ and $(p_1\cup \dots \cup p_m)\times [0,1],$ 
in exactly the same fashion as before, we associate an $\I$-cube map
\begin{equation} 
\psi_S: V_{P_0} \lra V_{P_1}
\end{equation} 
This map is a graded map of cubes over $\Rm$ of degree equal to 
the Euler characteristic of $S$ minus $\frac{m}{2}.$ 

Tensoring the map $\psi_S$ with the identity map of the 
skew-commutative $n$-cube $E_{\I}$ and 
passing to associated complexes, we obtain a map of complexes 
of graded $R$-modules 
\begin{equation} 
\psi'_S: \C(P_0)\lra \C(P_1)
\end{equation}
In general this map is not a morphism in the category $\Kom(\Rm)$
of complexes of graded $R$-modules and grading-preserving homomorphism,  
as it shifts the grading by 
 $\chi(S)-\frac{m}{2},$ but $\psi'_S$ becomes a morphism in 
$\Kom(\Rm)$ when the grading of $\C(P_0)$ or $\C(P_1)$ is 
appropriately shifted.

\section{Transformations}
\label{section-transformations} 

In this section we associate a quasi-isomorphism 
of complexes of graded $R$-modules $C(D) \lra C(D')$
to a Reidemeister move between two plane diagrams 
$D$ and $D'$ of an oriented link $L.$

\subsection{ Left-twisted curl}
\label{section-left-curl} 

 Let $D$ be a plane diagram with $n-1$ double points and 
 let $D_1$ be a diagram constructed from $D$ by adding 
 a left-twisted curl. Denote by $\I'$ the set of double points of $D_1,$ 
 by $a$ the double point in the curl and by $\I$ the set of double points 
 of $D.$ There is a natural bijection of sets $\I \to \I'\setminus \{ a\},$ 
 coming from identifying a double point of $D$ with the corresponding 
 double point of $D_1.$ We will use this bijection to identify the 
  two sets $\I$ and $\I'\setminus \{ a\}.$ 

  The crossing  $a$ of $D_1$ can be resolved in two ways. 
 The 0-resolution of $a$ is a diagram $D_2$ which is a disjoint 
union of $D$ and a circle. The 1-resolution 
is a diagram isotopic to $D$ and we will identify this 
diagram with $D.$  

\drawing{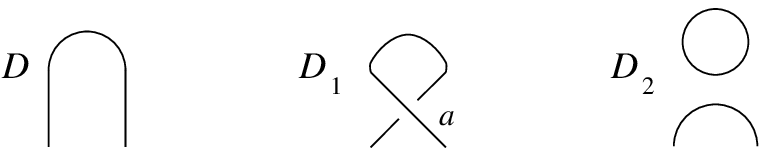}

In this section we will define a quasi-isomorphism of the complexes 
$C(D)$ and $C(D_1).$ This quasi-isomorphism arises from a 
splitting of the $\I'$-cube $V_{D_1}$ as a direct sum of two 
cubes, $V_{D_1}= V'\oplus V''.$ This splitting will induce 
 a decomposition of the complex $C(D_1)$ into a direct sum of an 
 acyclic complex and a complex isomorphic to  $C(D).$ 

Recall that $V_D,V_{D_1}$ and $V_{D_2}$ are the cubes associated with 
the diagrams $D,D_1$ and $D_2$ respectively. $V_{D_1}$ has index set 
$\I',$ while $V_D$ and $V_{D_2}$ are $\I$-cubes.

From the decomposition of $D_2$ as a union of $D$ and a simple 
circle we get a canonical isomorphism of cubes 
\begin{equation} 
V_{D_2} = V_D \o A
\end{equation} 
where $V_D\o A$ is the $\I$-cube obtained from $V_D$ by tensoring 
graded $R$-modules $V_D(\L), \L\subset \I$ with 
$A$ and tensoring the structure maps $\xi^{V_D}_a(\L)$ with the identity 
map of $A$. 

Let $U\subset \R^2$ be a small neighborhood of $a$  
that contains the curl: 

\drawing{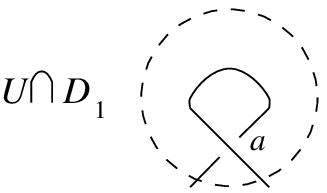} 

The picture above depicts how the  diagram 
$D_1$ looks inside $U.$ The boundary of $U$ is shown by a dashed 
circular line.  
Intersections of $U$ with diagrams $D$ and $D_2$ are depicted below 

\drawing{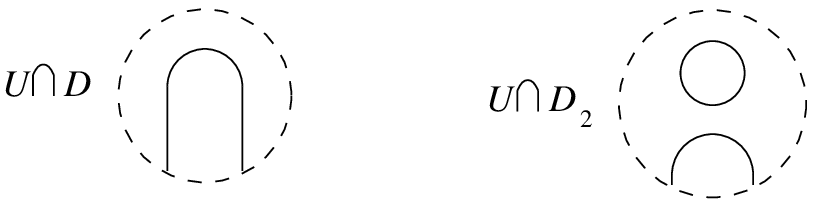}

Outside of $U$ diagrams $D, D_1$ and $D_2$ coincide. It is explained 
in Section~\ref{surfaces-comcube-morphisms} how surfaces in $U\times [0,1]$, 
satisfying certain conditions, give rise to cube maps. Using this 
construction we now define three cube maps between cubes $V_D$ and 
$V_{D_2}:$ 
\begin{eqnarray} 
& m_a: & V_{D_2} \lra V_D \\
& \Delta_a: & V_D \lra V_{D_2} \\
& \iota_a: & V_D \lra V_{D_2} 
\end{eqnarray} 

The map $m_a$ is associated to the following surface:

\drawing{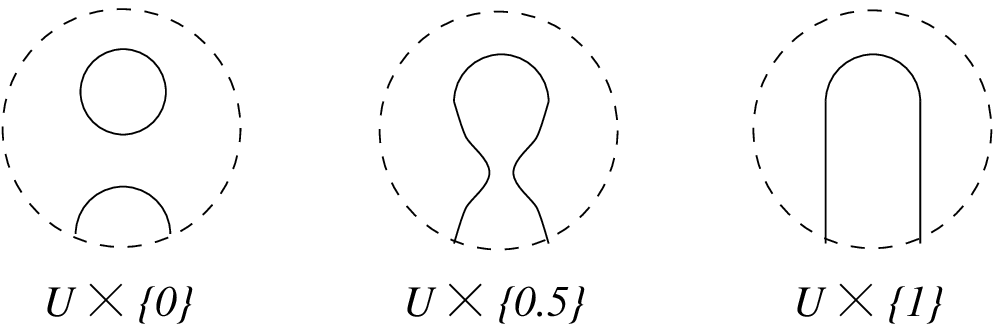}

Here and further on we depict surfaces embedded in $U\times [0,1]$ 
by a sequence of their cross-sections $U\times \{ t \}, t\in [0,1],$ 
the leftmost one being the 
intersection of the surface with $U\times \l 0 \r$, the rightmost 
being the intersection with $U \times \l 1\r.$ For such a surface 
$S\in U\times [0,1]$ we will call the projection $S\to [0,1]$ 
{\it the height function} of $S$. These surfaces will 
have only nondegenerate critical points relative to the height 
function. We depict enough sections of $S$ 
to make it obvious what surface we are considering, sometimes adding extra 
information, i.e., that the above surface has one saddle point and no other 
critical points relative to the height function. 

The intersections 
$S\cap U\times \{ 0 \}, S\cap U\times \{ 1 \}$ of the surface $S$ depicted 
above with the boundary disks $U\times \{0 \} , U\times \{ 1 \}$ 
are isomorphic to the intersections $D_2\cap (U\times \{ 0 \})$, 
respectively $D\cap (U\times \{ 1 \}).$ Thus, $S$ defines 
a map $m_a$ from the cube $V_{D_2}$ to $V_D.$ 

The cube map $\Delta_a$ is associated to the surface

\drawing{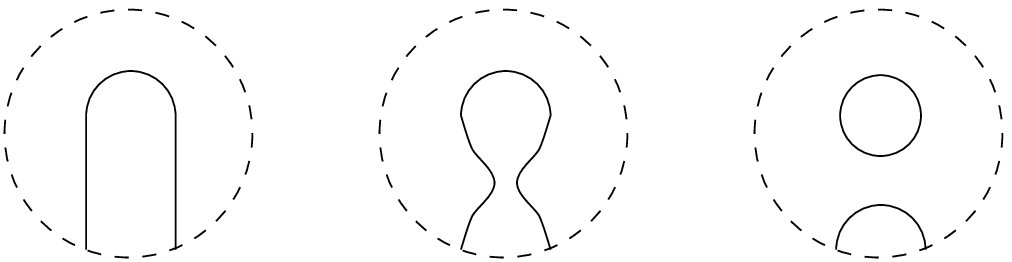}

 This surface has one saddle point and no other critical points relative 
to the height function. 

 $\iota_a$ is associated to 

\drawing{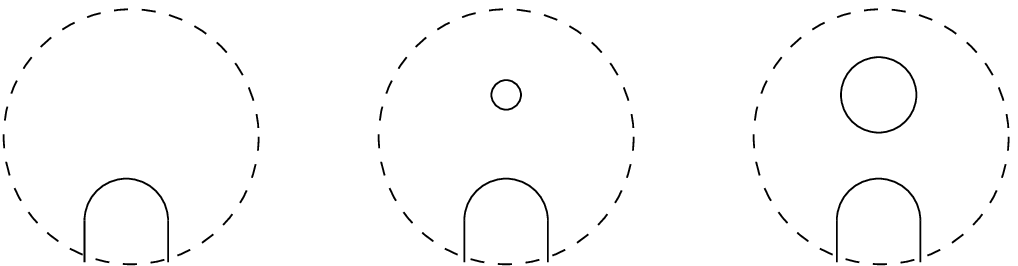}

  The only critical point of the height function is a local minimum.

The cube maps $m_a,\Delta_a,\iota_a$ are graded maps and change the 
grading by $-1,-1,1$ respectively. So let's keep in mind that 
$m_a,\Delta_a,\iota_a$ become grading-preserving if we appropriately 
shift gradings of our cubes, for example,

\begin{eqnarray} 
& m_a: & V_{D_2} \lra V_D \{ -1 \}  \\
& \Delta_a: & V_D \lra V_{D_2} \{ -1 \}  \\
& \iota_a: & V_D \lra V_{D_2} \{ 1 \}
\end{eqnarray} 

are grading-preserving maps of cubes over $\Rm.$ 

The composition $ m_a\iota_a$ is equal to  the identity map from $V_D$ 
to itself. Denote by $\jmath_a$ the map 
\begin{equation} 
\jmath_a\stackrel{\mbox{\scriptsize def}}{=}
\Delta_a-\iota_a m_a \Delta_a: V_D\lra V_{D_2}
\end{equation}  

The map $\jmath_a$ is a graded map of degree $-1.$ 

\begin{prop}
\label{comcube-splitting-iota-jmath} 
 The $\I$-cube
$V_{D_2}$ splits as a direct sum: 
\begin{equation} 
\label{jmath-splitting} 
V_{D_2} = \iota_a (V_D) \oplus \jmath_a(V_D).
\end{equation} 
\end{prop} 

\emph{Proof: } It is enough to consider the case when $D$ is a single 
circle. Then $\I'=\{ a\}, \I=\emptyset, V_D=A$ and 
 $\iota_a(V_D)= \mo \o A.$ But 

\begin{eqnarray*}
& \jmath_a \mo= & (\Delta_a-\iota_a m_a\Delta_a)\mo = 
X\o \mo - \mo \o X + c X\o X \\
& \jmath_a X = & (\Delta_a-\iota_a m_a\Delta_a)X = X \o X
\end{eqnarray*} 
and, thus, $A\o A$ is a direct sum of $\mo \o A$ and 
the $R$-submodule spanned by $\jmath_a\mo$ and $\jmath_a X.$

$\square$

Note that 
\begin{equation} 
m_a\jmath_a=m_a(\Delta_a-\iota_a m_a\Delta_a)=0
\end{equation} 
because $m_a\iota_a =Id.$ 

The $\I'$-cube $V_{D_1}$ contains $V_D$ and $V_{D_2}$ as subcubes
of codimension $1.$ 
Namely, we have canonical isomorphisms 
\begin{eqnarray} 
V_{D_1}(\ast 0) & \cong & V_{D_2} \label{D1-ast-0}\\
V_{D_1}(\ast 1) & \cong & V_D \{ -1 \} 
\end{eqnarray} 

Recall from Section~\ref{com-cubes} that $V_{D_1}(\ast 0)$ denotes 
  the $\I'\setminus \{ a\}$-cube (i.e. $\I$-cube) with 
  $V_{D_1}(\ast 0)(\L)= V_{D_1}(\L)$ for $\L\subset \I,$ etc. 

Under these isomorphisms the structure map $\xi_a^{V_{D_1}}$ 
(denoted below by $\xi_a$) for 
the $\I'$-cube $V_{D_1}$ 
\begin{equation} 
\xi_a: V_{D_1}(\ast 0)\lra V_{D_1}(\ast 1)
\end{equation} 
is equal to the map $m_a$ of $\I$-cubes, i.e., the following 
diagram is commutative 

\[ 
\begin{CD}
 V_{D_1}(\ast 0) @>{\xi_a}>> V_{D_1}(\ast 1) \\
@VV{\cong}V                  @VV{\cong}V  \\
V_{D_2}         @>{m_a}>>          V_D\{ -1 \} 
\end{CD}
\]

Using the splitting (\ref{jmath-splitting}) 
of $V_{D_2}$ we can decompose the $\I'$-cube 
$V_{D_1}$ as a direct sum of two $\I'$-cubes as follows: 
\begin{equation} 
\label{V-split} 
V_{D_1}= V'\oplus V''
\end{equation} 
where 
\begin{eqnarray} 
V'(\ast 0) & = & \jmath_a(V_D) \label{V-ast-0}\\
V'(\ast 1) & = & 0  \label{V-ast-1}\\
V''(\ast 0) & = & \iota_a(V_D) \\
V''(\ast 1) & = & V_{D_1}(\ast 1) 
\end{eqnarray} 
Some explanation: in the formula (\ref{V-ast-0}) $\jmath_a(V_D)$ 
is a subcube of $V_{D_2}$ 
and, due to (\ref{D1-ast-0}), $\jmath_a(V_D)$ sits inside $V_{D_1}$ 
as a subcube of codimension 1.  
Equation (\ref{V-ast-1}) means that $V'(\ast 1)(\L)=0$ 
for all $\L\subset \I'.$ Thus, $V'(\L)= \jmath_a(V_D(\L))\subset V_{D_1}(\L)$
 for $\L\subset \I',$ if $\L$ does not contain $a.$ If $\L$ contains $a,$ 
 $V'(\L)=0.$

Tensoring (\ref{V-split}) with $E_{\I'}$ we get a splitting of 
 skew-commutative $\I'$-cubes 
\begin{equation} 
V_{D_1}\o E_{\I'}= (V'\o E_{\I'}) \oplus (V''\o E_{\I'}) 
\end{equation} 
This induces a splitting of complexes associated to these 
skew-commutative $\I'$-cubes 
\begin{equation} 
\C(V_{D_1}\o E_{\I'})= \C(V'\o E_{\I'}) \oplus \C(V''\o E_{\I'}) 
\end{equation} 

\begin{prop} 
The complex $\C(V''\o E_{\I'})$ is acyclic.
\end{prop} 

\emph{Proof: }The complex $\C(V''\o E_{\I'})$ is isomorphic to the cone of 
the identity map of the complex $\C(V_D \o E_{\I})[-1]\{ -1 \} .$

$\square$

\begin{prop} 
The complexes $\C(V'\o E_{\I'})$ and $\C(V_D\o E_{\I})\{ 1 \}$ 
are isomorphic. 
\end{prop} 

\emph{Proof: }
We have a chain of isomorphisms of complexes 
\begin{eqnarray*} 
\C(V'\o E_{\I'}) & = & \C(V'(\ast 0)\o E_{\I}) \\
             & = & \C(V_D\{ 1 \} \o E_{\I}) \\
             & = & \C(V_D\o E_{\I}) \{ 1 \} 
\end{eqnarray*}

\begin{corollary} The complexes $\C(D_1)$ and $\C(D)\{ 1\}$ are 
quasiisomorphic.  
\end{corollary} 

{\emph Proof: } We have 
\begin{eqnarray*} 
\C(D_1) & = & \C(V_{D_1}\o E_{\I'}) \\
        & = & \C(V'\o E_{\I'}) \oplus \C(V''\o E_{\I'}) \\
        & = & \C(V_D\o E_{\I})\{ 1 \}\oplus \C(V''\o E_{\I'}) \\
        & = & \C(D)\{ 1\} \oplus (\mbox{Acyclic complex}) 
\end{eqnarray*} 

$\square$ 

Note that $x(D_1)=x(D)$ and $y(D_1)=y(D)+1.$ By 
(\ref{CD-defined}) 
\begin{equation} 
C(D)    =  \C(D)[x(D)]\{ 2x(D)-y(D)\} 
\end{equation}
and 
\begin{eqnarray*}   
C(D_1) & = & \C(D_1)[x(D_1)]\{ 2x(D_1) - y(D_1) \} \\
       & = & \C(D_1)[x(D)] \{ 2x(D) -y(D) -1\}.    
\end{eqnarray*} 
Therefore, complexes $C(D)$ and $C(D_1)$ are quasiisomorphic. 
Q.E.D. 

\subsection{Right-twisted curl}
\label{section-right-curl} 

Let $D$ be a diagram with $n-1$ double points and 
let $D_1$ be a diagram constructed from $D$ by adding 
a right-twisted curl. Denote by $a$ the new crossing that appears 
 in the curl. Let $\I$ be the set of crossings of $D$ and $\I'$ the 
  set of crossings of $D_1.$ We have a natural bijection of sets
 $\I \to \I'\setminus \{a\}$ and use it to identify these two sets. 
 
 Crossing $a$ can be resolved in two ways. 
$0$-resolution gives a diagram, isotopic to $D$ and canonically identified 
with $D$. $1$-resolution produces a diagram, denoted $D_2,$ 
which is a disjoint union of $D$ and a simple circle. 
 
\drawing{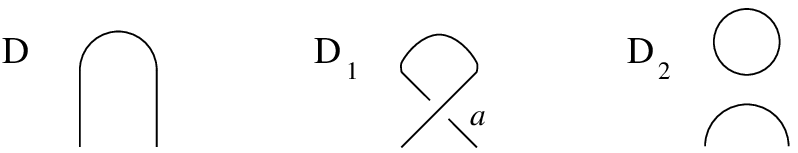}

Note that diagrams $D$ and $D_2$ are the same as diagrams 
$D$ and $D_2$ from Section~\ref{section-left-curl} and 
we will be using cube maps $m_a,\Delta_a,\iota_a$ defined 
in that section. Also define a map
\begin{equation} 
\epsilon_a: V_{D_2} \lra V_D 
\end{equation} 
where $\epsilon_a$ is associated to the surface 
\drawing{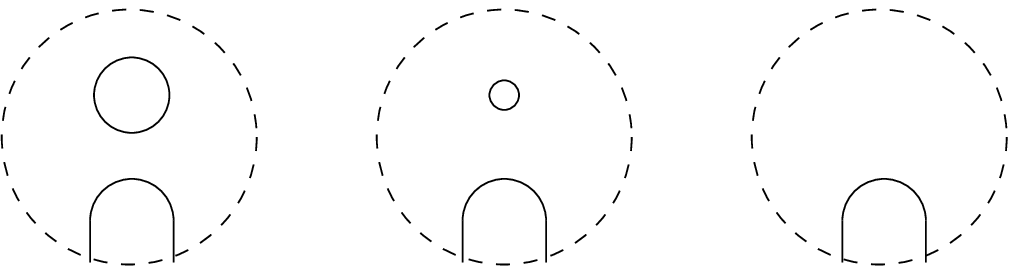}
This surface has one critical point relative to the height function 
and it is a local maximum. 
The cube map $\epsilon_a$ changes the grading by $1$ and becomes 
grading preserving after an appropriate shift: 
\begin{equation} 
\epsilon_a: V_{D_2} \lra V_D\{ 1 \}  
\end{equation} 

Let $\aleph$ be the map 
\begin{equation} 
\aleph= \iota_a - c \iota_a m_a \Delta_a: V_D\lra V_{D_2}
\end{equation} 
$\aleph$ is graded of degree $1.$ 

\begin{prop}
\label{comcube-splitting-aleph-delta}  We have a cube splitting
\begin{equation} 
\label{aleph-splitting} 
V_{D_2} = \aleph (V_D) \oplus \Delta_a(V_D)
\end{equation} 
\end{prop}

\emph{Proof: } It suffices to check this when $D$ is a simple circle. 
Then 
\begin{eqnarray*}
& \aleph (\mo) & = \mo \o \mo - 2c \mo \o X \\
& \aleph (X)   & = \mo \o X 
\end{eqnarray*} 
The $R$-submodule of $A\o A$ generated by these two vectors 
complements $\Delta(A)$ and there is direct sum decomposition of 
 $R$-modules 
 \begin{equation*} 
 A\o A = R\cdot \aleph(\mo)\oplus R\cdot \aleph(X) \oplus  \Delta(A)
 \end{equation*}
$\square$

Denote by $\wp$ the cube map 
\begin{equation} 
\wp=m_a-m_a\Delta_a \epsilon_a : V_{D_2} \lra V_D.
\end{equation} 
Note that $\wp$ is a graded map of degree $-1.$ 

\begin{lemma} \label{two-equalities} We have equalities 
\begin{eqnarray} 
\wp \Delta_a  & = & 0 \label{wp-delta}\\
\wp \aleph    & = & {{\mathrm{Id}}}(V_D)  \label{wp-aleph} 
\end{eqnarray} 
\end{lemma} 
\emph{Proof: } 
 Map $\wp \Delta_a: V_D \lra V_D$ is the zero map because 
\begin{equation} 
\wp \Delta_a= m_a \Delta_a -m_a\Delta_a\epsilon_a \Delta_a =
m_a\Delta_a -m_a\Delta_a = 0
\end{equation} 
(the second equality uses that $\epsilon_a\Delta_a =\mbox{Id}.$) 

The equality  (\ref{wp-aleph}) is checked similarly:  
\begin{eqnarray*} 
\wp \aleph 
& = & (m_a-m_a\Delta_a\epsilon_a)(\iota_a -c \iota_a m_a \Delta_a)  \\
& = & m_a\iota_a - c m_a\iota_a m_a\Delta_a -m_a\Delta_a \epsilon_a
\iota_a + c m_a\Delta_a\epsilon_a \iota_a m_a \Delta_a  \\
& = & {{\mathrm{Id}}}  - c m_a\Delta_a + c m_a\Delta_a - c^2 
m_a \Delta_a m_a\Delta_a  \\
& = & {{\mathrm{Id}}} - c^2 m_a\Delta_a m_a\Delta_a  \\ 
& = & {{\mathrm{Id}}} 
\end{eqnarray*}
The third equality in the computation above follows from 
the identities 
\begin{equation} 
m_a\iota_a ={{\mathrm{Id}}} , \hspace{0.3in} \epsilon_a \iota_a = -c.  
\end{equation} 
The fifth equality is implied by 
$m_a\Delta_a m_a\Delta_a =0.$ This identity follows from 
the nilpotence property $m\Delta m \Delta = 0$ of the structure 
maps $m$ and $\Delta$ of $A.$

Using the splitting (\ref{aleph-splitting}) 
of $V_{D_2}$ and Lemma~\ref{two-equalities}, 
we can decompose the $\I'$-cube $V_{D_1}$ 
 as a direct sum of two $\I'$-cubes as follows: 
\begin{equation} 
\label{V-split-again} 
V_{D_1}= V'\oplus V''
\end{equation} 
where 
\begin{eqnarray} 
V'(\ast 0) & = & 0 \label{V-prime-zero} \\
V'(\ast 1) & = & \aleph(V_D)\{ -1\}  \subset V_{D_2}\{ -1\}
= V_{D_1}(\ast 1) \\
V''(\ast 0) & = & V_D = V_{D_1}(\ast 0) \\
V''(\ast 1) & = & \Delta_a(V_D)\{ -1\}  \subset V_{D_2}\{ -1\} 
= V_{D_1}(\ast 1)
\end{eqnarray} 
Tensoring (\ref{V-split-again}) with $E_{\I'}$ we get a splitting 
 of skew-commutative  $\I'$-cubes
\begin{equation} 
V_{D_1}\o E_{\I'}= (V'\o E_{\I'}) \oplus (V''\o E_{\I'}) 
\end{equation} 
This induces a splitting of complexes associated to these 
skew $\I'$-cubes 
\begin{equation} 
\C(V_{D_1}\o E_{\I'})= \C(V'\o E_{\I'}) \oplus \C(V''\o E_{\I'}) 
\end{equation} 

\begin{prop} 
The complex $\C(V''\o E_{\I'})$ is acyclic.
\end{prop} 
\emph{Proof: }The complex $\C(V''\o E_{\I'})$ is isomorphic to the cone of 
the identity map of the complex $\C(V_D \o E_{\I})[-1].$

$\square$

\begin{prop} 
The complexes $\C(V'\o E_{\I'})$ and $\C(D)[-1]\{ -2 \}$ 
are isomorphic. 
\end{prop} 

\emph{Proof: }
We have a chain of isomorphisms of complexes 
\begin{eqnarray*} 
\C(V'\o E_{\I'}) & = & \C(V'(\ast 1)\o E_{\I})[-1] \\
             & = & \C(V_D\{ -2 \} \o E_{\I})[-1] \\
             & = & \C(V_D\o E_{\I}) [-1]\{ -2 \} \\
            & = & \C(D) [-1]\{ -2 \} 
\end{eqnarray*} 
The first isomorphism here follows from (\ref{V-prime-zero}) and is 
obtained by fixing an  isomorphism between  
skew-commutative $\I$-cubes $E_{\I'}(\ast 1)$ 
and $E_{\I}.$  The second isomorphism comes from an isomorphism 
$V'(\ast 1) = V_D\{ -2 \},$ induced by $\aleph.$   

$\square$

\begin{corollary} The complexes $\C(D_1)$ and $\C(D)[-1]\{ -2 \}$ are 
quasiisomorphic.  
\end{corollary} 

{\emph Proof: } We have 
\begin{eqnarray*} 
\C(D_1) & = & \C(V_{D_1}\o E_{\I'}) \\
        & = & \C(V'\o E_{\I'}) \oplus \C(V''\o E_{\I'}) \\
        & = & \C(D)[-1]\{-2\}\oplus \C(V''\o E_{\I'}) \\
        & = & \C(D)[-1]\{-2\}\oplus (\mbox{Acyclic complex}) 
\end{eqnarray*} 
$\square$ 

Note that $x(D_1)=x(D)+1$ and $y(D_1)=y(D).$ By 
(\ref{CD-defined}) 
\begin{equation} 
C(D)   =  \C(D)[x(D)]\{ 2x(D)-y(D)\} 
\end{equation} 
and 
\begin{eqnarray*}
C(D_1) & = & \C(D_1)[x(D_1)]\{ 2x(D_1) - y(D_1) \} \\
       & = & \C(D_1)[x(D)+1] \{ 2x(D) -y(D) +2\} 
\end{eqnarray*} 
Therefore, complexes $C(D)$ and $C(D_1)$ are quasiisomorphic. 
Q.E.D.

\subsection{The tangency move} 
\label{section-tangency-move} 

Let $D$ and $D_1$ be two diagrams that differ as depicted below

\drawing{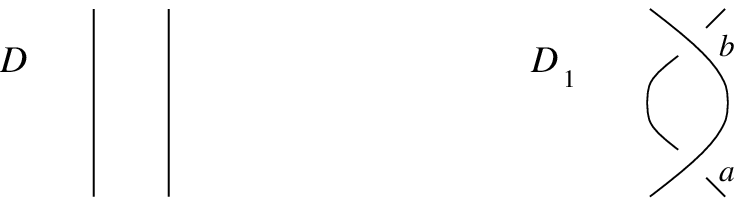} 

In this section we will construct a quasi-isomorphism of 
complexes $C(D)$ and $C(D_1).$ 

We assume that $D$ has $n-2$ double points. Consequently, $D_1$ 
has $n$ double points. Let $\I'$ be the set of double points of $D_1,$ 
let $\I$ be $\I'\setminus\{ a,b\}$ where $a$ and $b$ are double points 
 of $D_1$ depicted above. We identify $\I$ with the double points set 
 of $D.$ 

Denote by $d$ the differential 
of the complex $\C(D_1).$  
Consider diagrams  $D_1(\ast 00), D_1(\ast 01), D_1(\ast 10), D_1(\ast 11)$ 
obtained by resolving double points $a$ and $b$ of $D_1:$ 

\drawing{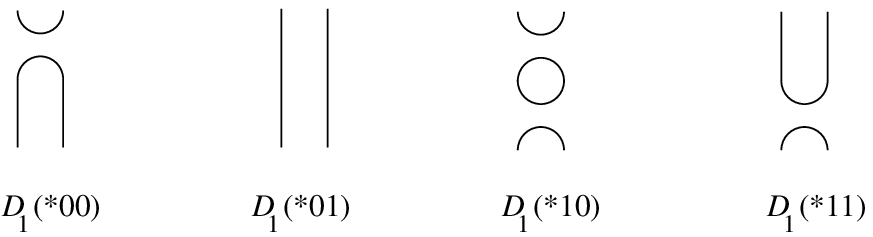}

(e.g., $D_1(\ast 01)$ is constructed from $D_1$ by taking $0$-resolution 
 of $a$ and $1$-resolution of $b,$ etc.) 
Each of these four diagrams has $\I$ as the set of its double points. 

To a diagram $D_1(\ast uv)$ where $u,v\in \{ 0,1\}$ there is 
associated the complex $\C(D_1(\ast uv))$ of graded $R$-modules. 
 Denote by $d_{uv}$ the differential in 
this complex:
\begin{equation} 
d_{ab}: \C(D_1(\ast uv)) \lra  \C(D_1(\ast uv)).
\end{equation} 
We denote by $d_{uv}^{(i)}$ 
the differential in shifted complexes 
\begin{equation} 
d_{uv}^{(i)}: \C(D_1(\ast uv))[i]\{ i \} \lra  
         \C(D_1(\ast uv))[i ]\{ i \} \mbox{ for }i\in \Z. 
\end{equation} 

The commutative $\I$-cube $V_{D_1}$ can be viewed as a commutative square 
of $\I'$-cubes 
\[
\begin{CD}
\label{CD-this}
V_{D_1(\ast 00)}   @>{\phi_2}>>    V_{D_1(\ast 01)}\{ -1\}   \\
@V{\phi_1}VV                   @VV{\phi_3}V \\
V_{D_1(\ast 10)}\{ -1 \} @>{\phi_4}>>  V_{D_1(\ast 11)}\{ -2 \}    
\end{CD}
\]
where $\phi_i, 1\le i\le 4$ denote the corresponding cube maps. 
Recall that these cube maps are associated to certain elementary 
surfaces (see Sections~\ref{diagrams-to-comcubes}, 
\ref{surfaces-comcube-morphisms}) that have one saddle point relative 
to the height function and no other critical points. For example, 
$\phi_1$ is associated to the surface 

\drawing{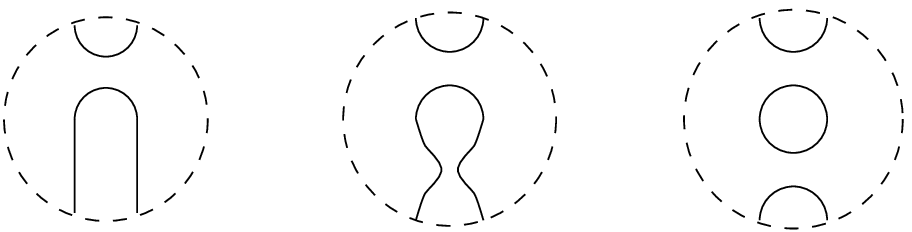} 

The maps $\phi_i$ induce maps $\psi_i$ between complexes: 
\begin{eqnarray*} 
\psi_1 & : & \C(D_1(\ast 00)) \lra \C(D_1(\ast 10))\{ -1 \} \\
\psi_2 & : & \C(D_1(\ast 00)) \lra \C(D_1(\ast 01))\{ -1 \} \\
\psi_3 & : & \C(D_1(\ast 01))[-1]\{ -1\}  \lra \C(D_1(\ast 11))[-1]\{ -2 \} \\
\psi_4 & : & \C(D_1(\ast 10))[-1]\{ -1 \}  \lra \C(D_1(\ast 11))[-1]\{ -2 \} 
\end{eqnarray*} 

We can decompose $\C(D_1),$ considered as a $\Z\oplus\Z$-graded 
$R$-module (see the end of Section~\ref{notations-complexes}), 
 into the following direct sum of $\Z\oplus\Z$-graded 
$R$-modules.  
\begin{eqnarray*} 
\C(D_1) & = & \C(D_1(\ast 00))
\oplus \C(D_1(\ast 01))[-1]\{-1 \}    \\
 & \oplus & \C(D_1(\ast 10))[-1]\{-1 \} 
\oplus \C(D_1(\ast 11)) [-2] \{ -2 \} 
\end{eqnarray*} 
Let's say a few words about this decomposition: $\C(D_1)$ is the direct sum 
of $R$-modules $V_{D_1}(\L)$ which sit in the vertices of the 
 $\I'$-cube $V_{D_1}.$ Since we  presented this 
 cube as a  commutative square of $\I$-cubes 
 $V_{D_1(\ast uv)}\{ -u-v\}$ for $u,v\in \{ 0,1\},$ the above decomposition
 results. Well, almost. 
 Indeed, when we pass from $\J$-cubes to complexes we tensor with 
 the fixed skew cube $E_{\J}.$ 
 To define the left hand side of the above formula we tensor 
 $V_{D_1}$ with the skew $\I'$-cube $E_{\I'},$ while for the 
 right hand side similar tensor products are formed with the 
 skew $\I$-cube $E_{\I}.$ 
 Therefore, we must say how we identify 
 $R$-modules which sit in the vertices of $E_{\I'}$ with $R$-modules 
 sitting in the 
 vertices of $E_{\I}.$ For $D_1(\ast 00):$ we map $E_{\I}(\L)$ where 
 $\L\subset \I$ to $E_{\I'}(\L)$ by sending $z\in o(\L)$ to $z\in o(\L).$ 
 For $D_1(\ast 10)$: map $E_{\I}(\L)$ where $\L\subset \I$ to 
  $E_{\I'}(\L a)$ by sending $z\in o(\L)$ to $za \in o(\L a).$ 
 Similarly for $D_1(\ast 01).$ For $D_1(\ast 11)$: map 
 $E_{\I}(\L)$ to $E_{\I'}(\L ab)$ by sending $z\in o(\L)$ to 
 $zab\in o(\L ab).$ This is not a canonical choice, since we could have 
 sent $z$ to $zba$ and would have gotten minus the original map. 
 So, to define the latter map, we implicitly fix an ordering of $a$ and $b.$ 
  
Note that the above decomposition
 is  not a direct sum of complexes, as the differential 
$d_{uv}^{(-u-v)}$ of $\C(D_1(\ast uv))$ differs from 
$d$ restricted to $\C(D_1(\ast uv))[-u-v]\{ -u-v\}\subset \C(D_1),$  
except when $u=v=1.$ Exactly, we have 
\[
\begin{array}{lllll}  
dx & = & d_{00}x+[-1]\psi_1 x +[-1]\psi_2x & \mbox{ for } & 
x\in \C(D_1(\ast 00)) \\
dx & = & -d_{01}^{(-1)}x - [-1]\psi_3 x & \mbox{ for } & 
x\in \C(D_1(\ast 01))[-1]\{ -1 \}  \\
dx & = & -d_{10}^{(-1)}x+[-1]\psi_4 x  & \mbox{ for } & 
x\in \C(D_1(\ast 10))[-1]\{ -1 \}  \\
dx & = & d_{11}^{(-2)}x & \mbox{ for } & x\in \C(D_1(\ast 11))[-2]\{ -2 \}  
\end{array}
\]
 {\it Some explanation:} applying $\psi_1$ to $x\in \C(D_1(\ast 00))$ 
 we get an element of $\C(D_1(\ast 10))\{ -1\},$ so that we shift 
 $\psi_1x$ by $[-1]$ to land it in $\C(D_1(\ast 10))[-1]\{ -1\}\subset 
 \C(D_1),$ etc. Various signs in the above formulas 
 come from our previous four identifications of the skew cube $E_{\I}$ 
 with codimension $2$ faces of $E_{\I'}.$ 

Let $\alpha$ be the map of complexes 
\begin{equation} 
\alpha: \C(D_1(\ast 01))[-1]\{ -1 \} \lra \C(D_1(\ast 10))[-1]\{-1 \}
\end{equation} 
associated to the surface 

\drawing{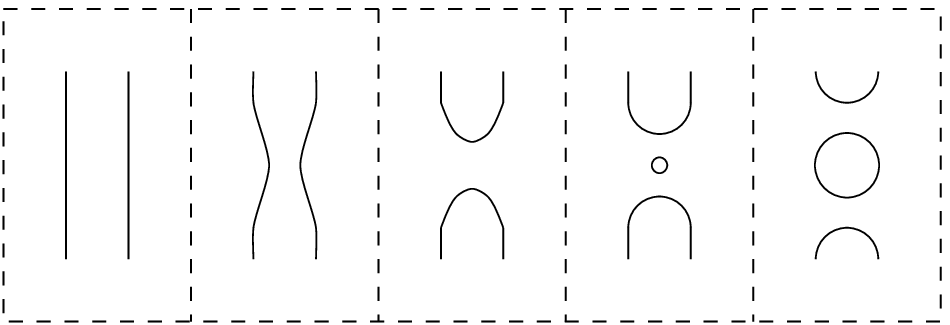} 
Considered as a map of $\Z\oplus \Z$-graded $R$-modules, 
$\alpha$ is grading-preserving. 

Let $\beta$ be the map of complexes 
\begin{equation} 
\beta: \C(D_1(\ast 11))[-2]\{ -2 \} \lra \C(D_1(\ast 10))[-1]\{ -1 \}
\end{equation} 
associated to the surface 

\drawing{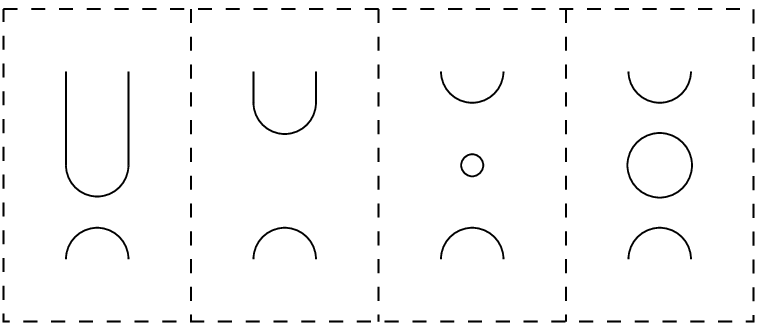} 
Note that $\beta$ is a graded map of degree $(-1,0).$ 

Let $X_1,X_2,X_3$ be $R$-submodules  of $\C(D_1)$ given by 
\begin{eqnarray} 
X_1& = & \{ z+\alpha(z)|z\in \C(D_1(\ast 01))[-1]\{ -1 \}  \} \\
X_2 & = & \{ z + dw | z,w \in \C(D_1(\ast 00)) \}  \\
X_3 & = & \{ z + \beta(w)| z,w \in \C(D_1(\ast 11))[-2 ]\{ -2\}  \} 
\end{eqnarray} 

\begin{prop} \label{prop-proved-first} 
These submodules are stable under  $d$:
\begin{equation} 
dX_i \subset X_i 
\end{equation} 
 and respect the $\Z\oplus \Z$-grading of $\C(D_1).$
\end{prop}

\emph{Proof: } Let us first check that $X_1,X_2$ and $X_3$ are 
direct sums of their graded components. For $X_2$ it follows 
from the fact that $\C(D_1(\ast 00))$ is a direct sum of its 
graded components and $d$ is graded of degree $(1,0).$ 
Submodule $X_3$ is graded because $\C(D_1(\ast 11))[-2]\{ -2\}$ 
is a direct sum of its graded components and $\beta$ is a graded 
map. Finally, $X_1$ is graded since $\alpha$ is grading-preserving. 

We now verify that these three submodules are stable under $d.$ 
For $X_2$ this is obvious.  
 To see it for $X_3$, notice that 
$dz\in \C(D_1(\ast 11))[-2]\{ -2\}$ whenever 
$z\in \C(D_1(\ast 11))[-2]\{ -2\}.$ Moreover, for such a $z,$ 
\begin{equation} 
\label{many-differentials} 
d\beta(z) = -d_{10}^{(-1)}\beta(z) + [-1]\psi_4\beta(z)= 
 - d_{10}^{(-1)}\beta(z)+ z = 
 \beta d_{11}^{(-2)}(z) + z
\end{equation} 
The second equality is implied by $[-1]\psi_4\beta = \mbox{Id}.$ 
Map $\psi_4\beta$ is associated to the surface 

\drawing{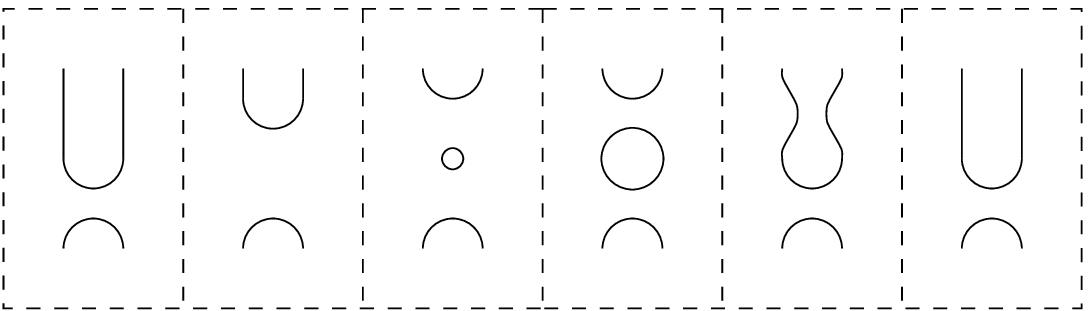}
obtained by composing surfaces to which $\psi_4$ and $\beta$ 
are associated. This surface is isotopic, through an isotopy fixing 
the boundary, to the surface 

\drawing{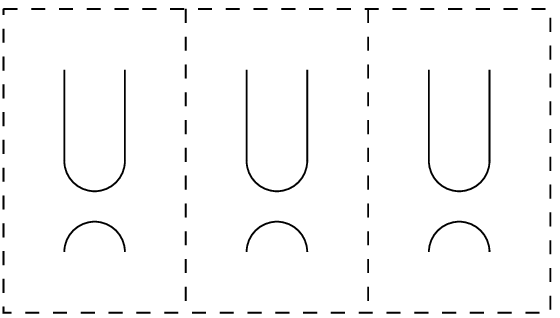} 
representing the identity map. Hence $[-1]\psi_4\beta =\mbox{Id}.$  

Formula (\ref{many-differentials}) implies that $X_3$ is stable 
under $d,$ since the rightmost term $\beta d_{11}^{(-2)}(z)+ z$ 
 lies in $X_3.$ 

Finally, to check the $d$-stability of $X_1,$ we compute, for 
$z\in \C(D_1(\ast 01))[-1]\{ -1\},$ 
\begin{eqnarray*}
 d(z+ \alpha(z)) & = & dz + d\alpha(z) \\
                 & = & 
  - d_{01}^{(-1)}z - [-1]\psi_3z -  d_{10}^{(-1)}\alpha(z) + 
  [-1]\psi_4\alpha(z)   \\
                 & = & - (d_{01}^{(-1)}z + d_{10}^{(-1)}\alpha(z)) + 
                      [-1](-\psi_3 z + \psi_4 \alpha (z)) \\
                 & = & -( d_{01}^{(-1)}z  + d_{10}^{(-1)}\alpha(z)) \\
                 & = & -(d_{01}^{(-1)}z  + \alpha d_{01}^{(-1)}z) \in X_1 
\end{eqnarray*} 
In the fourth equality we used that $\psi_4\alpha=\psi_3,$
in the fifth that $\alpha d_{01}^{(-1)}= d_{10}^{(-1)}\alpha,$ since 
$\alpha$ is a grading-preserving map of complexes.  

 $\square$ 

\begin{corollary} Submodules $X_1,X_2,X_3$ are graded 
subcomplexes of the complex $\C(D_1).$ 
\end{corollary} 

\begin{prop} 
\label{prop-tangency-move} 
\begin{enumerate} 
\item We have a direct sum decomposition 
\begin{equation} 
\label{direct-sum-three} 
\overline{C}(D_1) = X_1 \oplus X_2 \oplus X_3 
\end{equation} 
in the category ${{\mathrm{Kom}}}({\Rm})$ of complexes of graded $R$-modules. 
\item The complexes $X_2$ and $X_3$ 
are acyclic. 
\item The complex $X_1$ is isomorphic to 
the complex $\overline{C}(D)[-1]\{ -1\} .$
\end{enumerate} 
\end{prop}

\emph{Proof: } Since we already know that $X_1,X_2$ and $X_3$ are 
graded subcomplexes of $\C(D_1),$ it suffices to check 
(\ref{direct-sum-three}) on the level of underlying 
abelian groups. We have $\alpha= \beta \psi_3$ and, therefore, 
for $z\in \C(D_1(\ast 01))[-1]\{ -1\} $ 
\begin{equation} 
\alpha z = \beta \psi_3 z \in X_3
\end{equation} 
Subcomplex $X_1$ consists of elements $z+ \alpha z$ and we know that 
$\alpha z\in X_3.$ We are thus reduced to proving the following 
direct sum splitting of abelian groups 
\begin{equation} 
\C(D_1) =  \C(D_1(\ast 01))[-1]\{ -1\}\oplus X_2 \oplus X_3
\end{equation} 
Next recall that $X_2$ consists of elements $z+ dw$ for 
$z,w\in \C(D_1(\ast 00)).$ The differential $dw$ reads 
\begin{equation} 
dw= d_{00}w + [-1]\psi_1 w + [-1]\psi_2 w 
\end{equation} 
Note that $[-1]\psi_2 (w)\in \C(D_1(\ast 01))[-1]\{-1 \}$ and 
$d_{00}w \in \C(D_1(\ast 00)).$ Let $X'_2$ be the subgroup 
of $\C(D_1)$ given by 
\begin{equation}
X'_2= \{ z+ [-1]\psi_1 w | z,w\in \C(D_1(\ast 00))\} 
\end{equation} 
Then it is enough to verify that $\C(D_1)$ is a direct sum of 
its subgroups $ \C(D_1(\ast 01))[-1]\{ -1\}, X'_2$ and $X_3$:
\begin{equation} 
 \C(D_1)=  \C(D_1(\ast 01))[-1]\{ -1\}\oplus X'_2\oplus X_3. 
\end{equation} 
Note that $X_3$ contains $\C(D_1(\ast 11))[-2]\{ -2\}$ and
$X'_2$ contains $\C(D_1(\ast 00)).$   Recall the direct sum 
decomposition 
\begin{eqnarray*} 
\C(D_1) & = & \C(D_1(\ast 00))
\oplus \C(D_1(\ast 01))[-1]\{-1 \}    \\
 & \oplus & \C(D_1(\ast 10))[-1]\{-1 \} 
\oplus \C(D_1(\ast 11)) [-2] \{ -2 \} 
\end{eqnarray*} 
of $\C(D_1).$ Let $X''_2$ and $X'_3$ 
be the following abelian subgroups of $\C(D_1(\ast 10))[-1]\{ -1\}$: 
\begin{eqnarray*} 
X''_2 & = & \{ [-1]\psi_1 (w) | w \in \C(D_1(\ast 00)) \}  \\
X'_3  & = & \{ \beta (w)  | w \in \C(D_1(\ast 11))[-2]\{ -2\} \}
\end{eqnarray*} 
Now we are reduced to proving the direct sum decomposition 
\begin{equation} 
\label{almost-done} 
\C(D_1(\ast 10))[-1]\{ -1\} = X''_2 \oplus X'_3
\end{equation} 
in the category of abelian groups. As an abelian group, 
$\C(D_1(\ast 10))[-1]\{ -1 \} $ is a direct sum of 
$F(D_1(\L a ))$ over all possible resolutions 
 of the $(n-2)$-double points of $D_1.$ 
Similar direct sum splittings can be formed for $X''_2$ and $X'_3$ and 
one sees then 
 that it suffices to check (\ref{almost-done}) when $D_1$ has 
only two double points. There are two such $D_1$'s:  

\drawing{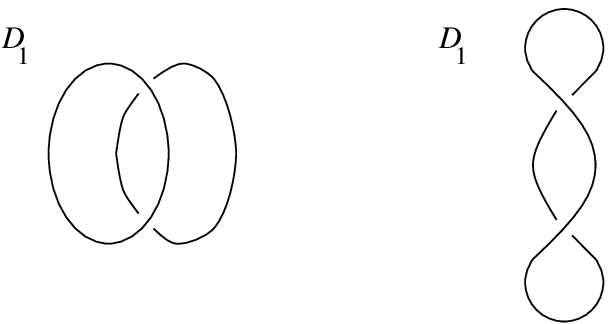}
In each of these two cases decomposition (\ref{almost-done}) follows from  
the splitting (\ref{decomp-one}). 
That proves part 1 of the proposition. 

We next prove part 2. 
 The complex $\C(X_2)$ is isomorphic to the cone of the identity map of 
$\C(D_1(\ast 00))[-1]$ and, therefore, acyclic. 
Similarly, $\C(X_3)$ is acyclic, being isomorphic to 
 the cone of the identity map of $\C(D_1(\ast 11))\{ -2\} [-2].$ 

To prove part 3 of the proposition, notice that the diagrams 
$D$ and $D_1(\ast 01)$ are isomorphic. This induces an isomorphism 
between the complexes 
\begin{equation} 
\C(D) = \C(D_1(\ast 01))
\end{equation} 
An isomorphism 
\begin{equation} 
\gamma: \C(D_1(\ast 01))[-1]\{ -1\} \stackrel{\cong}{\lra} X_1 
\end{equation}
is given by 
\begin{equation} 
\gamma(z) =(-1)^{i} (z + \alpha(z))
\end{equation} 
for $z\in \C^i(D_1(\ast 01))[-1]\{ -1\}.$ 
We need $(-1)^i$ in the above formula to match the 
differentials in these two complexes. 

$\square$ 

\begin{corollary} \label{quasi-quasi} 
The complexes $\C(D)[-1]\{ -1 \}$ and $\C(D_1)$ 
are quasiisomorphic. $\square$ 
\end{corollary} 

Note that $x(D_1)=x(D)+1$ and $y(D_1)=y(D)+1.$ 
From  (\ref{CD-defined}) we get 
\begin{eqnarray*} 
C(D) & = & \C(D)[x(D)]\{ 2y(D)-x(D)\} \\
C(D_1) & = & \C(D_1) [x(D)+1] \{ 2y(D)-x(D)+1\} 
\end{eqnarray*} 
which, together with Corollary~\ref{quasi-quasi},
implies that $C(D)$ is quasiisomorphic to $C(D_1).$ 
Q.E.D.

\subsection{Triple point move} 
\label{section-triple-point-move} 

We are given two diagrams with $n$ double points each, 
$D_1$ and $D_2$, that differ as depicted below. 

\drawing{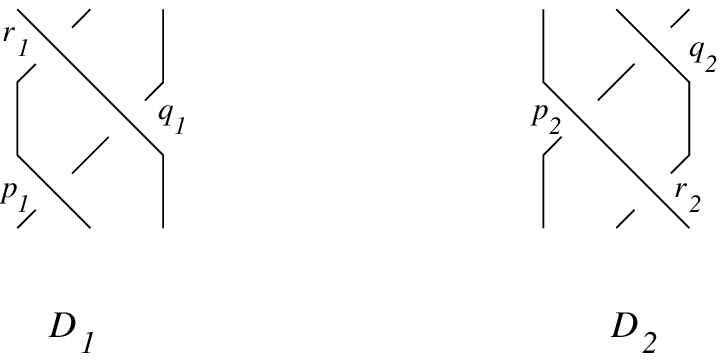} 

In this section we will construct a quasi-isomorphism of complexes 
$C(D_1)$ and $C(D_2).$ 

 Let $\I'$ be the set of double points of $D_1.$ We have
  $\I'= \I \sqcup \{ p_1, q_1, r_1\}$ where $\I$ are all double points 
  not shown on the above picture. In particular, we can identify 
 $\I \sqcup \{ p_2, q_2, r_2\}$ with the set of double points of $D_2.$ 

For starters, consider the diagrams 
$D_1(\ast 0), D_1(\ast 1), D_2(\ast 0), D_2(\ast 1),$ obtained by 
resolving double points $r_1$ of $D_1$ and $r_2$ of $D_2$: 

\drawing{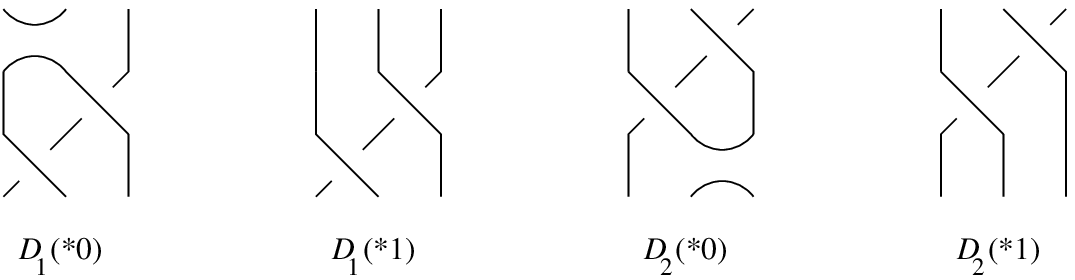} 

Note that diagrams $D_1(\ast 1)$ and $D_2(\ast 1)$ are 
isomorphic and that diagrams $D_1(\ast 0)$ and $D_2(\ast 0)$ 
represent isotopic links. 

We decompose $\C(D_1)$ and $\C(D_2)$ into following 
direct sums: 
\begin{equation} 
\label{a-not-so-long-formula}
\C(D_i)= \C(D_i(\ast 1))[-1]\{ -1\} \oplus 
\oplusop{u,v\in \{ 0,1\}} \C(D_i(\ast u v 0))[-u-v]\{ -u-v\} 
\end{equation} 
These are direct sum decompositions of $\Z\oplus \Z$-graded $R$-modules, 
not complexes. The diagrams $D_i(\ast uv0)$ for $i=1,2$ and $u,v\in \{ 0,1\}$ 
are depicted below
\drawing{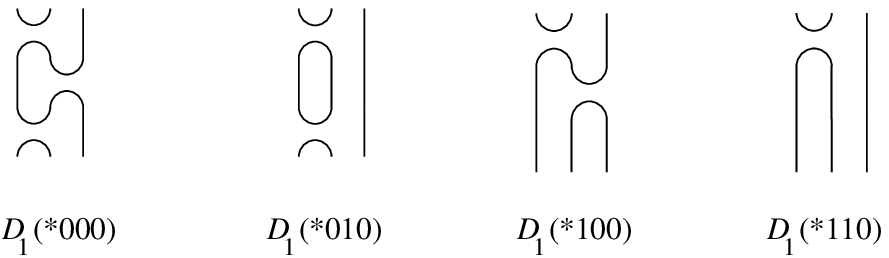} 

\drawing{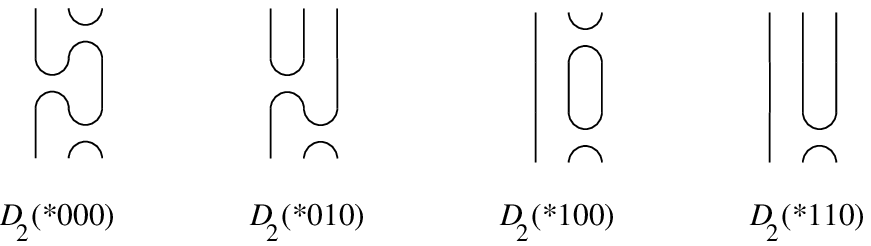} 

For all $i,u,v$ as above, we identify the set of double points of 
 $D_i(\ast uv0)$ with $\I.$ To fix the direct decomposition 
 (\ref{a-not-so-long-formula}) we need identifications between 
 the skew cube $E_{\I}$ and codimension $3$ facets of $E_{\I'}.$ 
 From the discussion in the previous section it should be clear how 
 these identifications are chosen. For instance, for $D_1(\ast 110),$ 
 we map $E_{\I}$ to a codimension $3$ facet of $E_{\I'}$ via 
 maps  $E_{\I}(\L) \to E_{\I'}(\L p_1 q_1)$ given by 
 $o(\L) \ni z \longmapsto z p_1 q_1 \in o(\L \sqcup \{ p_1,q_1\}).$ 

Let $\tau_1$ be the map of complexes 
\begin{equation} 
\tau_1: \C(D_1(\ast 100))[-1]\{ -1 \} \lra 
        \C(D_1(\ast 010))[-1]\{ -1 \} 
\end{equation} 
associated to the surface 
\drawing{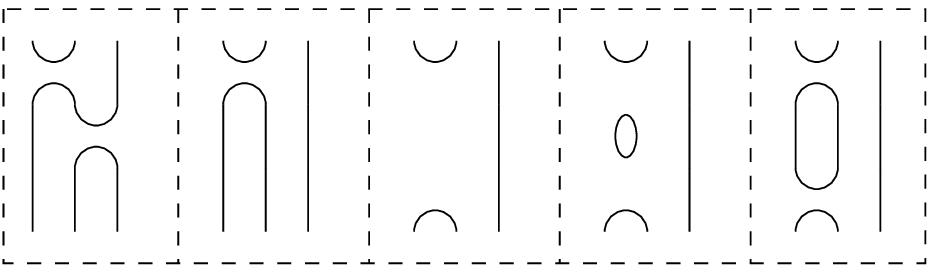} 
Relative to the height function this surface has two critical points, one 
of which is a saddle point and the other -- a local minimum. 
Considered as a map of $\Z\oplus \Z$-graded $R$-modules, $\tau_1$ is 
grading-preserving. 

Let $\delta_1$ be the map of complexes 
\begin{equation} 
\delta_1: \C(D_1(\ast 110))[-2]\{ -2 \} \lra 
        \C(D_1(\ast 010))[-1]\{ -1 \} 
\end{equation} 
associated to the surface 
\drawing{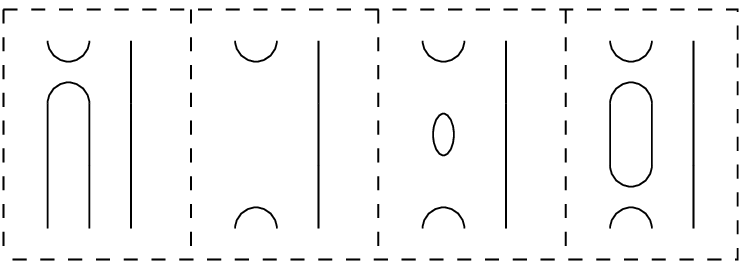}

Let $X_1,X_2,X_3$ be $R$-submodules of $\C(D_1)$ given by 

\begin{equation} 
\label{X-one-two-three} 
\begin{array}{ccl}  
X_1& = & \{ x+\tau_1(x)+y |x\in \C(D_1(\ast 100))[-1]\{ -1 \}, 
y \in \C(D_1(\ast 1))[-1]\{ -1\}  \} \\
X_2 & = & \{ x + d_1 y | x,y \in \C(D_1(\ast 000))\} \\
X_3 & = & \{ \delta_1(x)+ d_1\delta_1(y)| x,y\in \C(D_1(\ast 110))[-2]
\{ -2\}   \} 
\end{array} 
\end{equation} 
where $d_1$ denotes the differential of $\C(D_1).$ 
{\it Warning: } These $X_1,X_2,X_3$ have no relation to the 
 complexes $X_1,X_2,X_3$ considered in 
Section~\ref{section-tangency-move}.

Propositions~\ref{prop-first}-\ref{prop-last} below can be proved 
in the same fashion as 
Propositions~\ref{prop-proved-first} and \ref{prop-tangency-move} of the 
previous section. For this reason and to keep this paper from being 
too lengthy the proofs are omitted. 

\begin{prop} \label{prop-first} 
Submodules $X_1,X_2,X_3$ are stable under $d_1$ and 
respect the $\Z\oplus \Z$-grading of $\C(D_1).$
\end{prop}

\begin{corollary} Submodules $X_1,X_2,X_3$ are graded subcomplexes 
of the complex $\C(D_1).$ 
\end{corollary}

Let $\tau_2$ be the map of complexes 
\begin{equation} 
\tau_2: \C(D_2(\ast 010))[-1]\{ -1 \} \lra 
        \C(D_1(\ast 100))[-1]\{ -1 \} 
\end{equation} 
associated to the surface 
\drawing{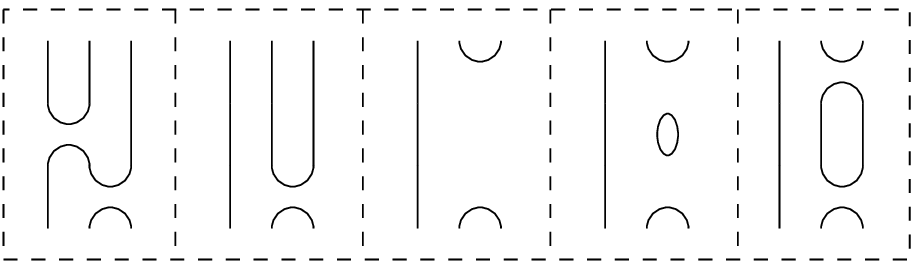} 

Considered as a map of $\Z\oplus \Z$-graded $R$-modules, $\tau_2$ is 
grading-preserving. 

Let $\delta_2$ be the map of complexes 
\begin{equation} 
\delta_2: \C(D_2(\ast 110))[-2]\{ -2 \} \lra 
        \C(D_2(\ast 100))[-1]\{ -1 \} 
\end{equation} 
associated to the surface 
\drawing{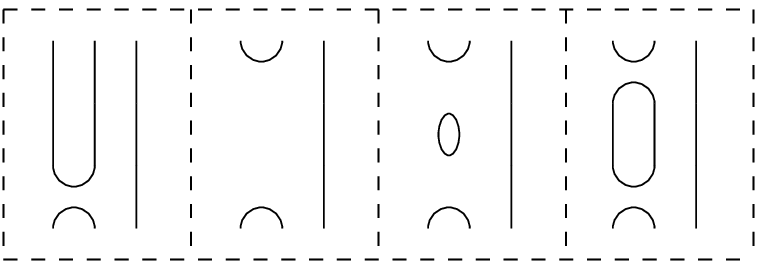} 

Let $Y_1,Y_2,Y_3$ be $R$-submodules of $\C(D_2)$ given by 

\begin{equation}
\label{Y-one-two-three}  
\begin{array}{ccl}  
Y_1& = & \{ x+\tau_2(x)+y |x\in \C(D_2(\ast 010))[-1 ]\{ -1\}, 
y \in \C(D_2(\ast 1))[-1]\{-1 \} \} \\
Y_2 & = & \{ x + d_2 y | x,y \in \C(D_2(\ast 000))\} \\
Y_3 & = & \{ \delta_2(x)+ d_2\delta_2(y)| x,y\in \C(D_2(\ast 110))
[-2]\{ -2\}   \} 
\end{array}
\end{equation}  
where $d_2$ stands for the differential of $\C(D_2).$ 

\begin{prop} 
These submodules are stable under $d_2$  and 
respect the $\Z\oplus \Z$-grading of $\C(D_2).$
\end{prop}

\begin{corollary} Subcomplexes $Y_1,Y_2,Y_3$ are graded subcomplexes 
of the complex $\C(D_2).$ 
\end{corollary} 

\begin{prop}\label{prop-last}  
\begin{enumerate}
\item We have direct sum decompositions  
\begin{eqnarray} 
\C(D_1) & = & X_1 \oplus X_2 \oplus X_3 \\
\C(D_2) & = & Y_1 \oplus Y_2 \oplus Y_3 
\end{eqnarray} 
\item The complexes $X_2, X_3, Y_2$ and $Y_3$ are 
acyclic. 
\item The complexes $X_1$ and $Y_1$ are isomorphic. 
\end{enumerate} 
\end{prop} 

\emph{Proof: } Parts 1 and 2 of this proposition are proved similarly to 
 Proposition~\ref{prop-tangency-move}. The isomorphism $X_1\cong Y_1$ 
 comes from the diagram isomorphisms 
 \begin{equation}
  \label{diagram-isomorphisms}  
 \begin{array}{rcl}  
 D_1(\ast 100)& = & D_2(\ast 010), \\
 D_1(\ast 1)  & = & D_2(\ast 1). 
 \end{array} 
 \end{equation} 
 These diagram isomorphisms induce isomorphisms of complexes 
   \begin{equation}
  \label{complexes-isomorphisms}  
 \begin{array}{rcl}  
 \C(D_1(\ast 100))& = & \C(D_2(\ast 010)) \\
 \C(D_1(\ast 1))  & = & \C(D_2(\ast 1))
 \end{array} 
 \end{equation} 
  which allow us to identify $x$ in the definition 
  (\ref{X-one-two-three}) of $X_1$ with $x$ in the definition 
  (\ref{Y-one-two-three}) of $Y_1$ and, similarly, identify $y$'s. 
  An isomorphism $X_1\cong Y_1$ of complexes is then given by 
 \begin{equation} 
  X_1 \ni x+\tau_1(x) + y \longmapsto x+\tau_2(x) + y \in Y_1. 
 \end{equation} 
 $\square$

\begin{corollary} Compexes $\C(D_1)$ and $\C(D_2)$ are 
quasiisomorphic. 
\end{corollary} 

The above isomorphism of complexes $X_1$ and $Y_1$ 
induces a quasi-isomorphism of  $\C(D_1)$ and $\C(D_2).$ 
Note that $x(D_1)= x(D_2)$ and $y(D_1)= y(D_2).$ Therefore, 
the complexes $C(D_1)$ and $C(D_2)$ are quasi-isomorphic and 
the cohomology groups $H^i(D_1)$ and $H^i(D_2)$ are isomorphic as 
graded $R$-modules. 
Q.E.D. 

This finishes the proof of Theorem~\ref{closed-i}. 

\section{Properties of cohomology groups} 
\label{groups-properties} 

\subsection{Some elementary properties} 

Pick an oriented link $L$ and a component $L'$ of $L.$ Let $L_0$ be 
$L$ with the orientation of $L'$ reversed and let $l$ be the linking number 
of $L'$ and $L\setminus L'.$ Fixing a plane diagram $D$ of $L$, 
we count $l$ as half the number of 
double intersection points in $D$ of $L'$ with $L\setminus L'$ 
with weights $+1$ or $-1$ according to the following convention

\drawing{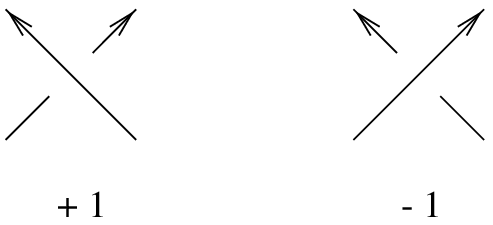} 

Denote by $D_0$ the diagram $D$ with the reversed orientation of $L'.$ 
Since $D_0$ and $D$ are the same as unoriented diagrams, 
$\C(D_0)= \C(D).$ Also
\begin{equation} 
x(D_0) = x(D) - 2l, y(D_0) = y(D) +2l
\end{equation} 
We obtain 
\begin{prop} 
\label{switch-by-2l}
For $L,L_0$ as above, there is an equality 
\begin{equation} 
H^i(L_0) = H^{i+2l}(L)\{ 2l\}
\end{equation} 
of isomorphism classes of graded $R$-modules. 
\end{prop} 

Let $K,K_1$ be oriented knots and $(-K)$ be $K$ with orientation reversed. 
In a similar fashion we deduce 
\begin{prop} There is an equality 
\begin{equation} 
H^i(K \# K_1) = H^i((-K)\# K_1)
\end{equation} 
of isomorphism classes of graded $R$-modules. 
\end{prop}

Let $D$ be a diagram of an oriented link $L$ and denote by 
$\cm(L)$ the number of connected components of $L.$ Then it is easy to 
see that $C^i_j(D)=0$ if parities of $j$ and $\cm(L)$ differ. 
This observation implies 

\begin{prop} 
\label{parity-vanishing} 
For an oriented link $L$ 
\begin{equation} 
H^{i,j}(L) = 0 
\end{equation} 
if $j+1\equiv  \cm(L) (\mbox{mod }2).$ 
\end{prop}

\subsection{ Computational shortcuts and cohomology of $(2,n)$ torus 
 links} 
\label{computations} 

Given a plane diagram $D$, a straighforward computation of 
cohomology groups $H^i(D)$ is daunting. These groups are cohomology 
 groups of the graded 
  complex $C(D)$ and the ranks of the abelian groups $C^i_j(D)$ 
 grow exponentially in the complexity of $D.$ 
Probably there is no fast algorithm for computing 
 $H^i(D)$, since these groups carry full information about 
  the Jones polynomial, computing  which is $\# P$-hard ([JVW]). 

 Yet, one can try to reduce $C(D)$ to a much smaller complex, albeit 
still exponentially  large, but more practical for a computation. 
 In this section we provide an example by simplifying $C(D)$ in 
  the case when $D$ contains a chain of positive half-twists and 
  apply our result by computing cohomology groups of $(2,n)$ torus 
  links. 

Let $D$ be a plane diagram with $n$ crossings and suppose that 
 $D$ contains a subdiagram pictured below 
 
\drawing{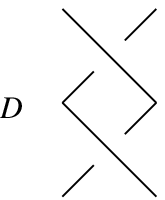} 

Four possible resolutions of these two double points of $D$ produce 
 diagrams $D(\ast 00), D(\ast 01),D(\ast 10), D(\ast 11)$: 

\drawing{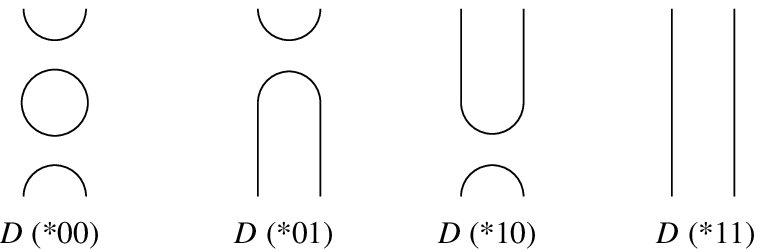} 

Note that diagrams $D(\ast 01)$ and $D(\ast 10)$ are isomorphic and 
 $D(\ast 00)$ is isomorphic to a union of $D(\ast 01)$ and a simple 
 circle. 
 The complex $\C(D)$ is isomorphic to the total complex of the bicomplex 
 \begin{eqnarray*} 
  \label{bicomplex} 
 \cdots \lra  0 & \lra & \C(D(\ast 00))\stackrel{\partial^0}{\lra}
  \C(D(\ast 01))\{ -1\} \oplus \C(D(\ast 10))\{ -1\} 
   \stackrel{\partial^1}{\lra}  \\
   & \stackrel{\partial^1}{\lra} & 
  \C(D(\ast 11))\{ -2\}  \lra 0 \lra \cdots 
  \end{eqnarray*} 
  where the differentials $\partial^0$ and $\partial^1$ are determined by 
  the structure maps of the skew $\I$-cube 
  $V_D\o E_{\I}$ (where $\I$ is the set of crossings of $D$). 
  Denote this bicomplex by $C.$ 
 To simplify notation we denote the diagram $D(\ast 01)$ by $D_0$ and 
 $D(\ast 11)$ by $D_1.$ 

\drawing{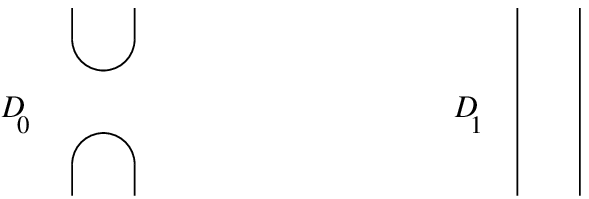} 

Then the bicomplex $C$ becomes 
 
 \begin{eqnarray*} 
  \label{bicomplex-01} 
 \cdots \lra  0 & \lra & \C(D_0)\o A\stackrel{\partial^0}{\lra}
  \C(D_0)\{ -1\} \oplus \C(D_0)\{ -1\}  \stackrel{\partial^1}{\lra} \\
  & \stackrel{\partial^1}{\lra} &  
  \C(D_1)\{ -2\}  \lra 0 \lra \cdots  
  \end{eqnarray*} 

Clearly, the differential $\partial^0,$ if restricted to 
  the subcomplex $\C(D_0)\o \mo$ of $\C(D_0)\o A,$ is injective, and 
  so the total complex of the subbicomplex 
  \begin{equation} 
  0 \lra C(D_0)\o \mo \stackrel{\partial^0}{\lra}  C(D_0) \{ -1 \}  \lra 0
  \end{equation} 
  of $C$ is acyclic. Denote this subbicomplex by $C_s$ and the 
quotient bicomplex by $C/C_s.$ The total complexes $\Tot(C)$ and 
$\Tot(C/C_s)$ of $C$ and $C/C_s$ are quasi-isomorphic, so to compute the 
cohomology of $\C(D)=\Tot(C)$ it suffices to find the cohomology 
of $\Tot(C/C_s).$ 
 
  We next give a precise description of the bicomplex $\Tot(C/C_s).$ 
  Let $u, l , w $ be maps of complexes 
\begin{eqnarray}  
u & : & \C(D(\ast 00))\lra \C(D_0)  \\
l & : & \C(D(\ast 00))\lra \C(D_0)  \\
w & : & \C(D_0) \lra \C(D_1) 
\end{eqnarray} 
induced by surfaces 
\drawing{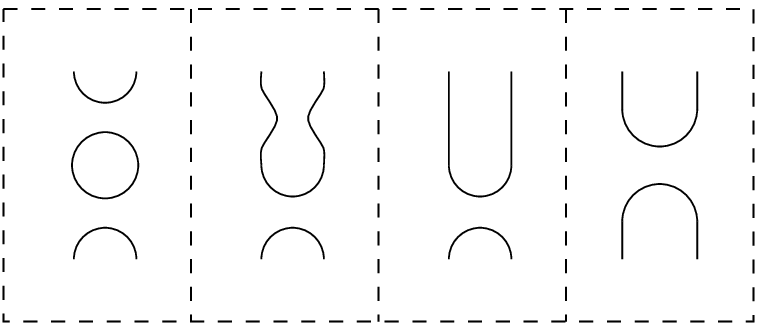} 
\vspace{0.3in} 
\drawing{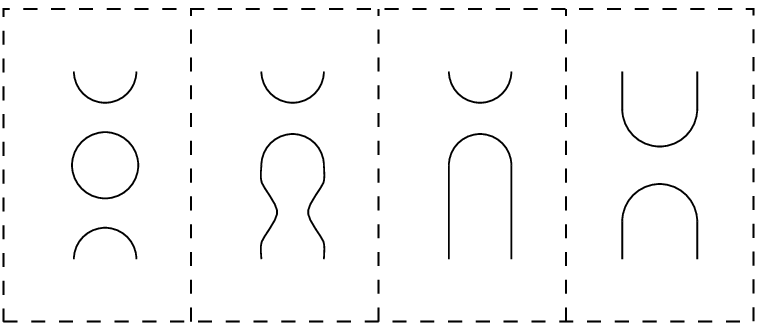} 
and 
\drawing{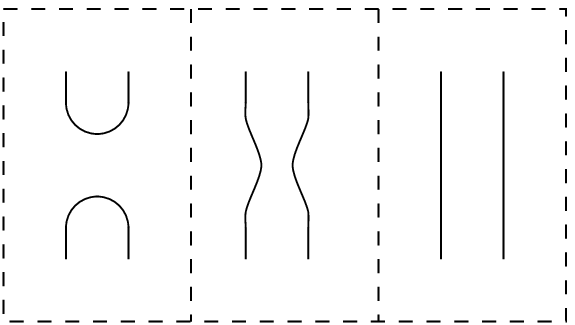} 
 respectively. Note that each of these maps have degree $-1,$ and 
 to make them homogeneous we need to shift gradings 
  of our complexes appropriately. We will use the same notations for shifted 
  maps since it will always be clear what the shifts are. 

  Let
  \begin{equation} 
  v: \C(D_0) \lra \C(D_0)\o A = \C(D(\ast 00))
  \end{equation} 
  be the map of complexes $v(t) = t\o X,   t\in \C(D_0).$ 
  The map $v$ has degree $-1.$ Denote by $u_X$ and $l_X$ the compositions 
  \begin{equation} 
  u_X= u \circ v , \hspace{0.3in} l_X= l \circ v. 
  \end{equation} 
  These are degree $-2$ maps of complexes and for each $i$ they induce 
  degree $0$ maps $\C(D_0)\{ i\} \lra \C(D_0)\{ i-2\},$ also denoted $u_X$ and 
  $l_X.$ 

  \begin{lemma} 
  The bicomplex $C/C_s$ is isomorphic to the bicomplex 
  \begin{equation} 
  \label{small-total}
  0 \lra \C(D_0)\{ 1 \} \stackrel{u_X-l_X}{\lra} \C(D_0)\{ -1\} 
  \stackrel{w}{\lra} \C(D_1)\{ -2\} \lra 0 
  \end{equation} 
  \end{lemma} 
  We skip the proof which is a simple linear algebra. $\square$  
 
  \begin{corollary} 
  \label{induction-base} 
  Cohomology groups $\overline{H}^i(D)$ are isomorphic to 
  the cohomology of the total complex of the bicomplex 
  (\ref{small-total}). 
  \end{corollary} 

  We thus see that the cohomology $\overline{H}^i(D)$ 
  of the diagram $D$ can be computed 
  via the quotient complex $\Tot(C/C_s)$ of $\C(D).$ The quotient complex 
  is smaller than the original one and computing its cohomology requires 
  less work. This reduction is not drastic since ranks of homogeneous 
  components of complexes $\C(D)$ and $\Tot(C/C_s)$ have the same order 
  of magnitude, but a similar reduction (described next,  when $D$ 
  contains a long chain of positive twists)
  leads to an effective computation of $H^i(D)$ for certain diagrams $D.$  

  Suppose that a diagram $D$ contains a chain of $k$ positive
  half-twists 
  \drawing{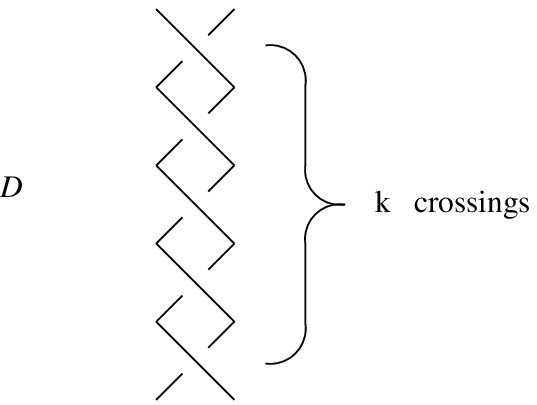} 
  As before, denote by $D_0$ and $D_1$   diagrams that 
  are suitable resolutions of the $k$-chain of $D.$ 
  \drawing{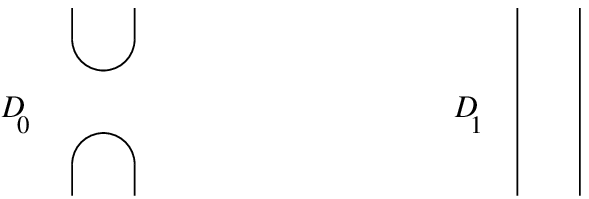} 
  From our previous discussion we retain degree $-2$ maps $u_X,l_X$ 
  and degree $-1$ map $w$ 
  between (appropriately shifted) complexes $\C(D_0)$ and $\C(D_1).$ 
  Let $C'$ be the bicomplex 
  \begin{eqnarray*}  
 0 & \lra & \C(D_0)\{ k-1 \} \stackrel{\partial^0}{\lra}  
         \C(D_0) \{ k-3 \} \stackrel {\partial^1} {\lra} \dots   \\
  & \stackrel{\partial^{k-3}}{\lra} &  \C(D_0) \{ 3-k \}  
  \stackrel{\partial^{k-2}}{\lra} \C(D_0) \{ 1-k \}  
   \stackrel{\partial^{k-1}}{\lra} \C(D_1) \{ -k \} \lra 0  
  \end{eqnarray*} 
  where 
  \begin{eqnarray*} 
   \partial^{k-1}& = & w \\
   \partial^{k-2} & = & u_X-l_X \\
   \partial^{k-3} & = & u_X+l_X \\
   \partial^{k-4} & = & u_X-l_X \\
                  & \dots &     \\
   \partial^{0}   & = & u_X-(-1)^k l_X, 
  \end{eqnarray*} 
  i.e. 
  \begin{equation} 
   \partial^{k-i} = u_X-(-1)^i l_X , \hspace{0.3in}\mbox{ for } 2\le i \le k . 
  \end{equation} 

  \begin{prop} 
  \label{small-complex} 
  The complex $\C(D)$ is quasiisomorphic to the total complex 
  $\Tot(C')$ of the bicomplex $C'.$ Cohomology groups 
  $\overline{H}^i(D)$ are isomorphic 
  to the cohomology groups of $\Tot(C').$ 
  \end{prop} 
  The proof goes by induction on $k,$ induction base $k=2$ being given 
   by Corollary~\ref{induction-base}, and consists of finding a suitable 
   acyclic subcomplex to quotient by. We omit the details. $\square$ 

   We conclude this section by applying this proposition to compute 
  cohomology groups of $(2,k)$ torus links. Fix $k>0$ and denote by $D$ 
  the diagram 

  \drawing{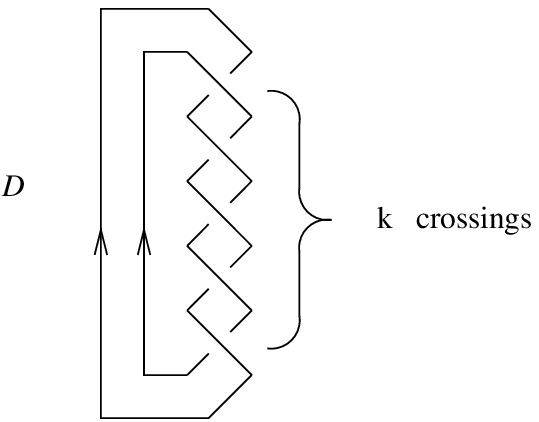}
  of the $(2,k)$ torus link $T_{2,k}.$ 
  
  The diagram $D_0$ is isomorphic to a simple circle and $D_1$ to a disjoint 
  union of two simple circles. Then $u_X=l_X$ is the operator  
  $A\to A$ of multiplication by $X$ 
  and the bicomplex $C'$ becomes a complex 
  \begin{eqnarray*} 
   0 & \lra & A \{ k-1 \} \lra  
         A \{ k-3 \} \lra \dots   \\
  &  \stackrel{0}{\lra} &  A\{ 5-k\}  
   \stackrel{2X }{\lra}  A \{ 3-k \}  
  \stackrel{0 }{\lra} A \{ 1-k \}  
   \stackrel{\Delta}{\lra} A\o A \{ -k \} \lra 0  
  \end{eqnarray*} 
  Recalling that $x(D)=k$ and $y(D)=0,$ we get 
  \begin{prop} The isomorphism classes of the graded $R$-modules 
  $H^i(T_{2,k})$ are given by 
   \begin{equation*} 
   \begin{array}{lll}  
   H^i(T_{2,k}) & = & 0 \hspace{0.2in}\mbox{ for }i<-k \mbox{ and }i> 0, \\
   H^0(T_{2,k}) & = & R\{ k\} \oplus R\{ k-2\}, \\
   H^{-1}(T_{2,k}) & = & 0,  \\
   H^{-2j}(T_{2,k}) & = & (R/2R)\{ 4j+k\} \oplus R\{ 4j-2+k\}
 \hspace{0.2in}\mbox{ for }1\le j \le \frac{k-1}{2},  \\
                          &  & j\in \Z,  \\
    H^{-2j-1}(T_{2,k}) & = & R\{ 4j+2+k\}
 \hspace{0.2in}\mbox{ for }1\le j \le \frac{k-1}{2}, j\in \Z, \\
   H^{-k}(T_{2,k}) & = & R\{ 3k\} \oplus R\{ 3k-2\} \hspace{0.2in} 
  \mbox{ for even } k. 
  \end{array}  
  \end{equation*} 
  \end{prop} 

\subsection{ Link cobordisms and maps of cohomology groups } 
\label{four-d-cobordisms} 

In this section by a surface $S$ in $\R^4$ we mean an oriented, compact
 surface $S,$ possibly with boundary, properly embedded in $\R^3\times 
 [0,1].$ The boundary of $S$ is then a disjoint union
 \begin{equation} 
  \partial S = \partial_0 S \sqcup  - \partial_1 S
 \end{equation} 
 of the intersections of $S$ with two boundary components of 
  $\R^3\times [0,1]$: 
\begin{eqnarray*} 
  \partial_0 S & = & (S\cap \R^3\times \{ 0\}) \\
   - \partial_1 S & = & (S\cap \R^3\times \{ 1\}) 
\end{eqnarray*} 
Note that $\partial_0 S$ and $\partial_1 S$ are oriented links 
 in $\R^3.$ 

The surface $S$ can be represented by a sequence $J$ of plane diagrams 
 of oriented links where every two consecutive diagrams in $J$ are 
 related either by one of the Reidemeister moves I-IV 
(Section~\ref{Reidemeister-moves}) or by one of the four moves 
depicted below (see [CS] where such representations by sequences of 
plane diagrams are studied in detail). 

\drawing{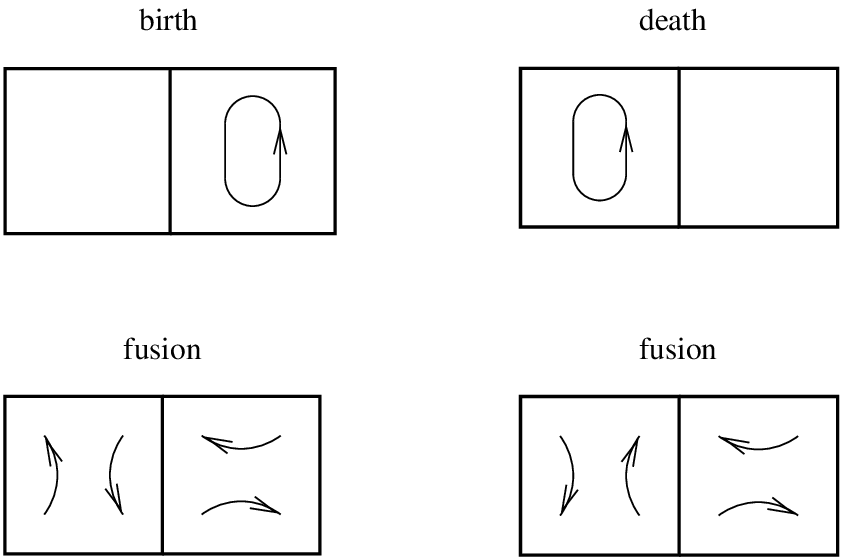} 

Following Carter-Saito [CS] and 
 Fisher [Fs], we call these moves {\it birth, death, fusion}
 ([CS] and [Fs] deal with the non-oriented version of these moves). 
We call $J$ {\it a representation of } $S.$ 
The first diagram in the sequence $J$ is necessarily a diagram 
 of oriented link $\partial_0 S,$ and the last diagram is a 
 diagram of $\partial_1 S.$ 

 The birth move consists 
 of adding a simple closed curve to a diagram $D.$ Denote the 
 new diagram by $D_1.$ Then $C(D_1)= C(D)\o_R A$ and the unit map 
 $\iota: R\to A$ of the algebra $A$ induces a map of complexes 
  $C(D) \to C(D_1).$ This is the map we associate to the birth move. 

 The death move consists of removing a simple circle from a 
 diagram $D_1$ to get a diagram $D.$ In this case the 
 counit $\epsilon: A\to R$ induces a map of complexes $C(D_1) \to C(D).$ 

 Finally, to a fusion move between diagrams $D_0$ and $D_1$ 
 we associate a map $C(D_0) \to C(D_1)$ corresponding to the elementary 
  surface with one saddle point in the manner discussed in 
  Section~\ref{surfaces-comcube-morphisms}. 

 In Section~\ref{section-transformations} to each Reidemeister 
  move between diagrams $D_0$ and $D_1$ we associated a 
  quasi-isomorphism map of complexes $C(D_0)\lra C(D_1).$ 

Given a 
representation $J$ of a surface $S$ by a sequence 
 of diagrams,
  denote the first and last diagrams of $J$ by $J_0$ and 
 $J_1$ respectively. Then to $J$ we can associate a map 
 of complexes 
 \begin{equation} 
  \varphi_J: C(J_0) \to C(J_1) 
 \end{equation} 
  which is the composition of maps associated to elementary transformations 
 between consecutive diagrams of $J.$ The map $\varphi_J$ induces  
 a map of cohomology groups 
  \begin{equation} 
  \theta_J: H^{i,j}(J_0) \to H^{i,j+\chi(S)}(J_1), \hspace{0.3in} i,j\in \Z.  
 \end{equation} 
 We now ready to state our main conjecture. 

  \begin{conjecture} 
  If two representations $J,\widetilde{J}$ of a surface $S$ have the 
  property that  

  (a) diagrams $J_0$ and $\widetilde{J}_0$ are isomorphic, 
 
  (b) diagrams $J_1$ and $\widetilde{J}_1$ and isomorphic, 
 
  then the maps $\theta_J$ and $\theta_{\widetilde{J}}$ are equal, up 
  to an overall minus sign, $\theta_J=\pm \theta_{\widetilde{J}}.$ 
  \end{conjecture} 

  In other words, we conjecture that, after a suitable $\Z_2$ extension
  of the link cobordism category, 
   our construction associates 
  honest cohomology groups $H^{i}(L)$ to oriented links $L$ in $\R^3$ 
  (and not just isomorphism classes of groups) and associates 
  homomorphisms 
  between these groups to isotopy classes of oriented surfaces 
  embedded in $\R^3\times [0,1].$ 
  In the categorical language, we expect to get a functor 
  from the category of ($\Z_2$-extended) oriented link cobordisms
   to the category of bigraded $R$-modules and module homomorphisms. 

  Suppose that the above 
  conjecture is true. Then, in the case of 
   a closed oriented surface $S$
  embedded in $\R^4,$ the map $\theta_S$ 
  of cohomology groups is a homomorphism 
  from $R$ 
  to itself (since $\partial S= \emptyset$ and the cohomology of the 
  empty link is equal to the 
  ground ring $R$). This homomorphism has degree $\chi(S)$ and 
  is automatically $0$ when $\chi(S)<0.$ Thus, the conjectural invariants 
  are zero whenever $S$ has empty boundary and 
  the Euler characteristic of $S$ is negative. 
  If, again, $\partial S = \emptyset$ and 
   the Euler characteristic of $S$ is nonnegative 
  (when $S$ is connected, $S$ is then necessarily a 2-sphere or 
  a 2-torus), the homomorphism $\theta_S:R\to R$ 
  is determined by $\theta_S(1)= k c^{\frac{\chi(S)}{2}} $ and 
  amounts to an integer number $k.$ Hence, we expect to have 
  integer-valued invariants of closed oriented surfaces with non-negative 
  Euler characteristic, embedded in $\R^4.$

\section{Setting $c$ to $0$} 
\label{c-is-zero} 

\subsection{Cohomology groups $\cH^{i,j}$}
\label{small-cohomology}  
Setting $c=0$ and taking $\Z$ instead of $R=\Z[c]$ as the base 
ring, everything 
we did in Sections~\ref{section-prelim}, \ref{section-diagrams} and 
\ref{section-transformations} 
goes through in exactly the same manner. The role of the ring $A$ will 
be played by the free graded abelian group $\cal A$ of rank $2$ 
with generators 
$\mo$ and $X$ in degrees $1$ and $-1$ correspondingly. $\cal A$ has 
commutative algebra and cocommutative coalgebra structures: 
\begin{eqnarray} 
 & & \mo^2 = \mo, \hspace{0.1in} \mo X = X \mo = X, 
 \hspace{0.1in} X^2 =0 \\
 & & \Delta(\mo) = \mo \o X + X\o \mo, \hspace{0.1in}
 \Delta(X) = X \o X 
\end{eqnarray} 
and the identity (\ref{m-and-delta}) holds. 
By abuse of notations, we use $m$ and $\Delta$ to denote multiplication 
and comultiplication in $\cA.$ Earlier $m$ and $\Delta$ were used 
to denote multiplication and comultiplication in $A.$ 
As in Section~\ref{A-and-surfaces}, we construct 
a functor $\cal F$ from the category $\MM$ of closed 
one-manifolds and cobordisms 
between them to the category of graded abelian groups and graded 
homomorphisms. To a disjoint union of $k$ circles functor $\cal F$ 
assigns the group ${\cal A}^{\o k}.$ To elementary 
surfaces $S_2^1,S_1^2,S_0^1,S_1^0, S_2^2$ and $S_1^1$ 
 (see Section~\ref{A-and-surfaces}) 
functor $\cal F$ assigns maps $m,\Delta, \iota, \epsilon,\mbox{Perm}$ and 
$\mbox{Id}$  
between suitable tensor powers of $\cA.$  The maps $\iota: \Z \to \cA$ and 
$\epsilon: \cA \to \Z$ are given by 
\begin{equation} 
\iota(1)= \mo, \hspace{0.2in} \epsilon(\mo)=0, \hspace{0.2in} 
 \epsilon(X) =1,  
\end{equation} 
while $\mbox{Perm}$ is just the permutation map $\cA\o \cA\to \cA \o \cA.$ 

To a diagram $D$ of an 
oriented link $L$ we can then associate a commutative $\I$-cube ${\cal V}_D$ 
of graded abelian groups and grading-preserving 
homomorphism, by the same procedure 
as the one described in Section~\ref{diagrams-to-comcubes}, 
using the functor $\cal F$ instead of $F$. In particular, for $\L\subset \I$ 
we have $\cV_{D}(\L)= {\cal F}(D(\L))\{ -|\L|\}$ 
where $\{ k\}$ shifts the grading down by $k.$ 

Let $\cal{AB}$ be the category of graded abelian groups and grading-preserving 
homomorphisms. Let ${\cal E}_{\I}$ be the skew-commutative 
 $\I$-cube $E_{\I}\o_R \Z $ over $\cal{AB}.$ 

Tensoring $\cV_D$ with ${\cal E}_{\I}$ over $\Z,$ we get a skew 
commutative  $\I$-cube 
$\cV_D\o {\cal E}_{\I}$ over the category $\cal{AB}.$ 
From this skew-commutative 
$\I$-cube we get a complex $\overline{\cC}(\cV_D\o {\cal E}_{\I})$
of graded abelian groups with a grading-preserving differential 
(Section~\ref{section-cubes-complexes}). 
Denote this complex by $\overline{\cal C}(D)$ and  
by ${\cal C}(D)$ the shifted complex 
\begin{equation} 
{\cal C}(D)=\overline{\cal C}(D)[x(D)]\{2 x(D)-y(D)\}
\end{equation} 
If we consider $\Z$ as a graded $R$-module, concentrated in degree $0,$ 
so that $c\Z=0,$ then 
\begin{equation} 
   \overline{\cal C}(D) = \C(D)\o_R\Z \hspace{0.2in} 
  \mbox{ and } \hspace{0.2in} 
  {\cal C}(D) = C(D)\o_R\Z. 
\end{equation} 

To a plane diagram $D$ of 
an oriented link $L$ we thus associate a complex of graded abelian 
groups ${\cal C}(D)$. Denote the $i$-th cohomology group 
of the $j$-th graded summand of ${\cal C}(D)$ by 
${\cal H}^{i,j}(D).$ These cohomology groups are finitely-generated 
abelian groups. For each diagram $D$ as we vary $i$ and $j$ over all integers, 
only a finite number of these groups are non-zero. 

\begin{theorem}
\label{theorem-no-c} 
For an oriented link $L,$ 
 isomorphism classes of abelian groups ${\cal H}^{i,j}(D)$ 
do not depend on the choice of a diagram $D$ of $L$ 
and are invariants of $L.$ 
\end{theorem}

\emph{Proof: } Set $c=0$ in the proof of Theorem~\ref{closed-i}. 
$\square$ 

For a diagram $D$ of the link $L,$ denote the isomorphism 
 classes of ${\cal H}^{i,j}(D)$ 
by ${\cal H}^{i,j}(L).$ 

Denote by $\overline{\cC}^i(D)$ (respectively, 
$\cC^i(D)$) the $i$-th group of the complex 
$\overline{\cC}(D)$ (respectively, $\cC(D)$) 
and by $\overline{\cC}^i_j(D)$ (respectively, by $\cC^i_j(D)$) 
the $j$-th graded 
component of $\overline{\cC}^i(D)$ (respectively, $\cC^i(D)$), 
so that $\overline{\cC}^i(D)= \oplusop{j\in \Z}\overline{\cC}^i_j(D), $  
respectively, $\cC^i(D)=\oplusop{j\in \Z}{\cC}^i_j(D).$
For an   diagram $D$ denote by $\cH^i(D)$ the graded 
abelian group $\oplusop{j\in \Z}\cH^{i,j}(D).$ In other words, 
$\cH^i(D)$ is the $i$-th cohomology group of $\cC(D).$ 
Denote by 
$\overline{\cH}^{i}(D)$ the $i$-th 
cohomology group of the complex $\overline{\cal C}(D)$ and by 
$\overline{\cH}^{i,j}(D)$ the $j$-th graded component of 
$\overline{\cH}^i(D),$ so that 
$\overline{\cH}^i(D)= \oplusop{j\in \Z} \overline{\cH}^{i,j}(D).$ 

\subsection{Properties of $\cH^{i,j}$: 
Euler characteristic, change of orientation} 

The Kauffman bracket of an oriented link $L$ is equal to the graded  
Euler characteristic of the cohomology groups $\cH^{i,j}(L),$ as stated in 
the following proposition. 

\begin{prop}
\label{no-c-Kauffman-bracket} 
 For an oriented link $L,$ 
\begin{equation} 
K(L) = \sum_{i,j\in \Z} (-1)^i q^j 
\ddim_{\Q}({\cal H}^{i,j}(L)\o \Q),
\end{equation} 
where $K(L)$ is the scaled 
Kauffman bracket (see Section~\ref{Kauffman-bracket}). 
\end{prop} 

The proof is completely analogous to that of 
formula (\ref{first-characteristic}). $\square$ 

The statements and proofs of 
Propositions~\ref{switch-by-2l}-\ref{parity-vanishing} 
transfer without change 
to the case of cohomology groups $\cH^{i,j},$ as indicated below.  

Let $L$ be an oriented link and $L'$ a component of $L.$ Denote by 
$l$ the linking number of $L'$ with its complement $L\setminus L'$ in $L.$ 
Let $L_0$ be the link $L$ with the orientation of $L'$ reversed. 

\begin{prop} For $i,j\in \Z$ there is an equality of isomorphism classes 
 of abelian groups
\begin{equation} 
\cH^{i,j}(L_0) = \cH^{i+2l,j+2l}(L). 
\end{equation} 
\end{prop} 

\begin{prop} 
\label{no-change} 
Let $K$ and $K_1$ be oriented knots and $(-K)$ be $K$ with the reversed
orientation. Then 
\begin{equation} 
{\cal H}^{i,j}(K \# K_1)= {\cal H}^{i,j}((-K)\# K_1)
\end{equation} 
\end{prop} 
Similarly to Proposition~\ref{parity-vanishing} we can prove

\begin{prop} For an oriented link $L$ 
\begin{equation} 
\cH^{i,j}(L) = 0 
\end{equation} 
if $j+1\equiv \cm(L) ({\mathrm{mod }}2).$ 
\end{prop}

\subsection{ Cohomology of the mirror image} 

Let $L$ be an oriented link and denote by $L^!$ the mirror image of $L$. 
Let $D$ be a diagram of $L$ with $n$ crossings, $\I$ the set of these 
crossings, and $D^!$ the corresponding diagram of $L^!$: 

\drawing{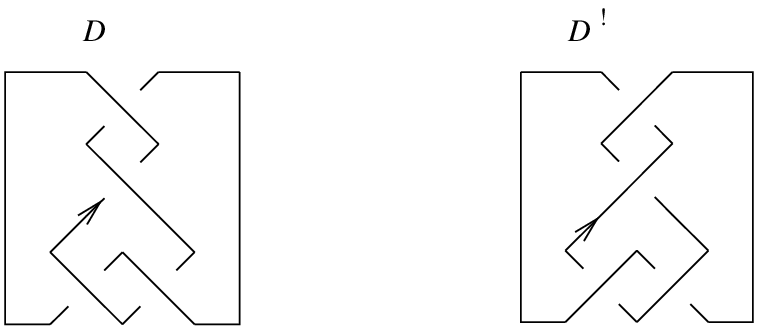} 

If $M$ is a graded abelian group, $M= \oplusop{j\in \Z}M_j,$ define the 
dual graded abelian group $M^{\ast}$ by 
$(M^{\ast})_j=\mbox{Hom}(M_{-j},\Z).$ 
The dual map $f^{\ast}: N^{\ast}\to M^{\ast}$ 
of a map $f: M\to N$ is defined 
 as the dual of $f$ in the sense of linear algebra. 

For $\L\subset \I$ denote by $\widetilde{\L}$ the complement 
 $\I \setminus \L.$ 

Let $\cV$ be a commutative $\I$-cube over the category $\cal{AB}$ of 
graded abelian groups and grading-preserving homomorphisms. Define the 
dual cube $\cV^{\ast}$ by 
\begin{equation} 
 \cV^{\ast}(\L) = ( \cV(\wta))^{\ast}  
\end{equation} 
 and the structure map
\begin{equation}  
\xi_a^{\cV^{\ast}}: \cV^{\ast}(\L)\lra \cV^{\ast}(\L  a)
\end{equation}
  of $\cV^{\ast}$ being the dual of the structure
 map
 \begin{equation} 
 \xi_a^{\cV}: \cV( \wta\setminus a  ) \lra \cV(\wta ) 
 \end{equation} 
 of $\cV.$  

 Denote by $\{ s\}$ the automorphism of the category $\cal{AB}$ that 
 shifts the grading down by $s.$ If $\cV$ is a commutative $\I$-cube over 
 $\cal{AB},$ denote by $\cV\{ s\}$ the commutative $\I$-cube $\cV$ with 
  the grading of each group $\cV(\L)$ shifted down by $s.$ 

\begin{prop} Let $D$ be a diagram with $n$ crossings and 
$D^!$ the dual diagram (see above). Then the 
commutative $\I$-cube $\cV_{D^!}\{ -n\} $ is isomorphic to the 
 dual $(\cV_D)^{\ast}$ of the $\I$-cube $\cV_D.$ 
\end{prop} 

\emph{Proof: } 
 Introduce a basis $\{ \mo^{\ast}, X^{\ast}\}$ 
  in the abelian group ${\cal A}^{\ast}= \mbox{Hom}({\cal A}, \Z)$ by 
  \begin{equation} 
  \mo^{\ast}(\mo) =0, \mo^{\ast}(X)=1, X^{\ast}(\mo)=1, X^{\ast}(X) = 0
  \end{equation} 
  Denote by $m^{\ast}, \Delta^{\ast}$ maps dual to $\Delta$ and $m$ 
  respectively: 
  \begin{eqnarray*} 
   m^{\ast}& : & {\cal A}^{\ast}\o {\cal A}^{\ast} \lra {\cal A}^{\ast} \\
   \Delta^{\ast}& : & {\cal A}^{\ast} \lra {\cal A}^{\ast} \o {\cal A}^{\ast} 
  \end{eqnarray*} 
  Then, in the basis $\{ \mo^{\ast}, X^{\ast}\} $ these maps are 
  \begin{eqnarray*} 
  & & m^{\ast}(\mo^{\ast}\o X^{\ast}) = 
        m^{\ast}(X^{\ast}\o \mo^{\ast})= X^{\ast} \\
  & &  m^{\ast}(\mo^{\ast}\o \mo^{\ast})= \mo^{\ast} \\
  & &  m^{\ast}(X^{\ast}\o X^{\ast}) = 0 \\
  & & \Delta^{\ast}(\mo^{\ast}) = \mo^{\ast}\o X^{\ast}+ X^{\ast}\o 
  \mo^{\ast} \\
  & & \Delta^{\ast}(X^{\ast}) = X^{\ast}\o X^{\ast} 
  \end{eqnarray*}
  Hence, under 
  the isomorphism $\mu: {\cal A} \to {\cal A}^{\ast}$ of graded 
  abelian groups,  
  given by $\mu(\mo)=\mo^{\ast}$ and $\mu(X) = X^{\ast},$ 
  maps $m,\Delta$ become $m^{\ast},\Delta^{\ast}.$

   Note that the $\wta$-resolution $D(\wta)$ 
  of diagram $D$ and the $\L$-resolution $D^!(\L)$ of $D^!$ are isomorphic. 
  Let $k$ be the number of circles in $D(\wta).$ Then 
  \begin{equation}  
  {\cal F}(D(\wta))={\cal F}(D^!(\L))= \cA^{\o k}
  \end{equation} 
 and, via $\mu,$ we can identify 
  \begin{equation} 
  \label{A-dual} 
  {\cal F}(D(\wta)) = (\cA^{\ast})^{\o k}= ({\cal F}(D^!(\L)))^{\ast}   
  \end{equation}
  Since $\mu$ maps $m,\Delta$ to $m^{\ast},\Delta^{\ast},$   we see that 
  after suitable shifts (recall that $\cV_D(\L)$ is equal to 
  ${\cal F}(D(\L))$ shifted up by $|\L|$) the identification (\ref{A-dual})  
  extends to an 
  isomorphism of $\I$-cubes $\cV_{D^!}\{ -n\} $  and $(\cV_D)^{\ast}.$ 

   $\square$ 

  Given a complex $C$ of graded abelian groups and grading-preserving 
  homomorphisms 
  \begin{equation}
 \cdots \lra C^i \stackrel{d^i}{\lra} C^{i+1} \lra \cdots 
  \end{equation} 
  define the dual complex $C^{\ast}$ by 
  $(C^{\ast})^i = (C^{-i})^{\ast}$ and the differential $(d^{\ast})^i$ 
  being the dual of the differential $d^{-i-1}$ of $C.$ 
   
  From the last proposition we easily obtain 
  \begin{prop} The complex $\cC(D^!)$ is isomorphic to the 
   dual of the complex $\cC(D).$ 
  \end{prop}   

   This implies 
    \begin{corollary} For an oriented link $L$ and integers $i,j$ there are 
   equalities of isomorphism classes of abelian groups 
  \begin{eqnarray} 
  \cH^{i,j}(L^!) \o \Q & = & \cH^{-i,-j}(L)\o \Q \\
   \Tor(\cH^{i,j}(L^!)) & = & \Tor(\cH^{1-i,-j}(L))
  \end{eqnarray} 
  where $\Tor$ stands for the torsion subgroup. 
  \end{corollary} 
  
Note that this corollary provides a necessary condition for a link to be 
amphicheiral.

\subsection{Cohomology of the disjoint union and connected sum of knots} 

Pick diagrams $D_1,D_2$ of oriented links $L_1,L_2$ and consider 
a diagram $D_1\sqcup D_2$ of the disjoint union $L_1\sqcup L_2.$ 
We then have an isomorphism of cochain complexes 
\begin{equation} 
{\cal C}(D_1\sqcup D_2) = {\cal C}(D_1) \o {\cal C}(D_2)
\end{equation} 
of free graded abelian groups. From the K\"unneth formula we derive 
\begin{prop} There is a short split exact sequence of cohomology 
groups 
\begin{eqnarray*} 
0 & \to  &  \oplusop{i,j\in \Z}( {\cal H}^{i,j}(D_1)\o {\cal H}^{k-i,m-j}(D_2))
\to {\cal H}^{k,m} (D_1\sqcup D_2) \to  \\
 & \to  & \oplusop{i,j\in \Z} 
 \Tor_1^{\Z}({\cal H}^{i,j}(D_1), {\cal H}^{k-i+1,m-j}(D_2))\to 0
\end{eqnarray*} 
\end{prop} 

\begin{corollary} 
 For each $k,m \in \Z$ there 
  is an equality of isomorphism classes of abelian groups 
  \begin{eqnarray*} 
   {\cal H}^{k,m} (L_1\sqcup L_2) = & & 
   \oplusop{i,j\in \Z}({\cal H}^{i,j}(L_1)\o {\cal H}^{k-i,m-j}(L_2)) \oplus \\
   & &  \oplusop{i,j\in \Z} 
 \Tor_1^{\Z}({\cal H}^{i,j}(L_1), {\cal H}^{k-i+1,m-j}(L_2)) 
 \end{eqnarray*} 
\end{corollary}

  Let $D_1,D_2$ be diagrams of oriented knots $K_1,K_2.$ We 
  picture $D_1$ and $D_2$ schematically as 

  \drawing{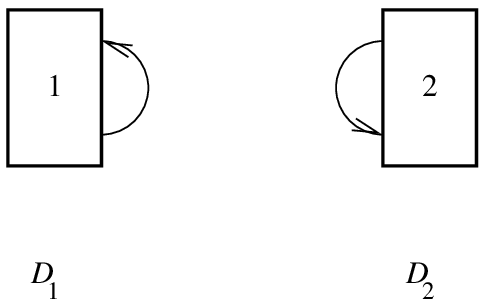} 
  Consider diagrams $D_3,D_4$ and $D_5,$ depicted below, of oriented links 
  $K_1\sqcup K_2,$  $K_1\# K_2,$ and  $K_1\# (-K_2)$: 

  \drawing{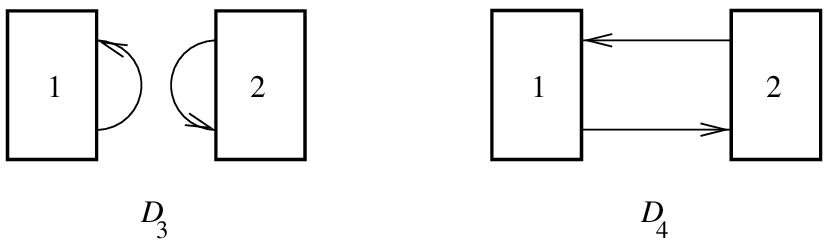} 
  \vspace{0.3in} 
  \drawing{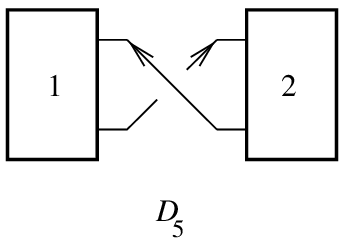} 

  By  resolving the central double point of $D_5,$ 
  we get a short exact sequence of complexes of graded abelian groups
  \begin{equation} 
  0 \lra \overline{\cC}(D_3)[-1]\{ -1\} \lra 
         \overline{\cC}(D_5) \lra \overline{\cC}(D_4) \lra 0 
  \end{equation} 
  After shifts, we obtain an exact sequence 
    \begin{equation} 
  0 \lra \cC(D_3)[-1]\{ -1\} \lra 
         \cC(D_5)[-1]\{ -2\}  \lra \cC(D_4) \lra 0 
  \end{equation} 
  which induces, for every integer $j,$ a long exact sequence of 
  cohomology groups 
  \begin{equation} 
  \begin{array}{ccccccccc} 
  \cdots & \to  & \cH^{i-1,j-1}(D_3) & \to &  \cH^{i-1,j-2}(D_5) 
  & \to &  \cH^{i,j}(D_4) & \to &  \\
  & \to  & \cH^{i,j-1}(D_3) & \to &  \cH^{i,j-2}(D_5) 
  & \to  & \cH^{i+1,j}(D_4) & \to &  \cdots 
  \end{array} 
  \end{equation} 
  Since the diagrams $D_3,D_4$ and $D_5$ represent oriented links
   $K_1\sqcup K_2, K_1\# K_2$ and  $K_1\# (-K_2),$ respectively, then  
    in view of Proposition~\ref{no-change}, we obtain
  
  \begin{prop} For oriented knots $K_1,K_2,$ 
  the isomorphism classes of the abelian groups $\cH^{i,j}(K_1\sqcup K_2), 
  \cH^{i,j}(K_1\# K_2)$ can be arranged into 
  long exact sequences 
  \begin{equation} 
   \begin{array}{ccccccc}  
  \to  & \cH^{i-1,j-1}(K_1\sqcup K_2) & \to &  \cH^{i-1,j-2}(K_1\# 
  K_2) 
  & \to &  \cH^{i,j}(K_1\# K_2) & \to   \\
  \to  & \cH^{i,j-1}(K_1\sqcup K_2 ) & \to &  \cH^{i,j-2}(K_1\# K_2) 
  & \to &  \cH^{i+1,j}(K_1\# K_2) & \to  
  \end{array}
  \end{equation} 
  \end{prop} 

\subsection{A spectral sequence} 
\label{spectral-sequence} 

Let $D$ be an plane diagram of a link. In this section we 
 construct a spectral sequence whose $E_1$-term is made of groups 
 $\cH^{i,j}(D)$ and which converges to cohomology groups 
 $H^{i,j}(D).$ 

  Due to the direct sum decomposition $R=\oplusop{k\ge 0} 
  c^k \Z$ of abelian groups, we have an abelian group decomposition 
  \begin{equation} 
  \label{two-complexes} 
  C^{i}_j(D) = \oplusop{k\ge 0} \cC^{i}_{j-2k}(D)
  \end{equation} 
  where, recall, $C^i_{j}(D)$ and $\cC^{i}_{j}(D)$ 
  were defined in Sections~\ref{diagrams-to-comcubes} and 
  \ref{small-cohomology} respectively.  
  Let us fix a $j\in \Z.$ Denote by $d$ the differential 
  in the weight $j$ subcomplex of the complex $C(D)$:  
  \begin{equation} 
  \cdots  \stackrel{d}{\lra} C^{i-1}_{j}(D) \stackrel{d}{\lra} 
  C^{i}_{j}(D) \stackrel{d}{\lra} C^{i+1}_{j}(D) 
  \stackrel{d}{\lra} \cdots 
  \end{equation} 
  and by $\partial$ the differential in the complex 
  \begin{equation} 
    \cdots \stackrel{\partial}{\lra} 
   \cC^{i-1}_{j-2k}(D) \stackrel{\partial}{\lra} 
    \cC^{i}_{j-2k}(D) \stackrel{\partial}{\lra} 
   \cC^{i+1}_{j-2k}(D) \stackrel{\partial}{\lra} \cdots 
  \end{equation} 
  where we suppress the dependence of $\partial$ on $k$.  
  Under the identification (\ref{two-complexes}) differential 
  $d$ becomes a differential of the complex 
    \begin{equation} 
    \cdots \stackrel{d}{\lra}  \oplusop{k\ge 0} 
   \cC^{i}_{j-2k}(D) \stackrel{d}{\lra} \oplusop{k\ge 0} \cC^{i+1}_{j-2k}(D) 
   \stackrel{d}{\lra} \cdots 
  \end{equation}
  Consider a bigraded abelian group 
  \begin{equation}   
  C=\oplusop{k\ge 0, i\in\Z} \cC^{i}_{j-2k}(D) 
  \end{equation} 
  where we set the grading of  $\cC^{i}_{j-2k}(D)$ to  $(i,-k).$ 
  We thus have a bigraded abelian group $C$ and 
  two maps, $d$ and $\partial$, from  $C$ to $C.$  
  Map $\partial$ is bigraded of degree $(1,0)$  while 
  $d$ is only graded relative to the first grading. However, we 
  can decompose 
  \begin{equation} 
  d= \partial + \widetilde{\partial} 
  \end{equation} 
  where $\widetilde{\partial}$ has grading $(1,-1)$  and satisfies 
   \begin{eqnarray} 
   \widetilde{\partial} \widetilde{\partial} & = & 0, \\
   \partial\widetilde{\partial} +\widetilde{\partial}\partial & = & 0
   \end{eqnarray} 
   Besides, since $\partial$ is a differential, $\partial\partial=0.$ 
   Therefore, $H^{i,j}(D)$ is equal to the $i$-th cohomology group 
  of the total complex of the bicomplex $(C,\partial, \widetilde{\partial}).$ 
   Group  $\cH^{s,j-2k}(D)$ is equal to the $s$-th cohomology group 
   of the subcomplex (relative to the differential $\partial$)  
   \begin{equation} 
  \cdots  \stackrel{\partial}{\lra} 
   \cC^{i-1,j-2k}(D) \stackrel{\partial}{\lra} 
  \cC^{i,j-2k}(D) \stackrel{\partial}{\lra} C^{i+1,j-2k}(D) 
  \stackrel{\partial}{\lra} \cdots 
  \end{equation}     
   of $C.$ Therefore, for each $j\in \Z$ we get a spectral sequence 
   whose $E_1$-term  is given by cohomology groups $\cH^{s,j-2k}(D)$, 
   $s\in \Z, k\ge 0$ and which converges to cohomology groups 
   $H^{i,j}(D).$ In the few cases where we managed to compute 
   cohomology groups, we have $H^{i,j}(D)= \oplusop{k\ge 0}\cH^{i,j-2k}(D)$ 
   and,consequently, the spectral sequence degenerates at $E_1.$ 
   We have no idea whether this is true for any diagram $D.$ 

\subsection{Examples} 
\label{examples} 

 Perhaps the graded 
 groups $\cH^i(D)$ are easier to compute than  $H^i(D).$ The latter 
  are computed via the complex $C(D)$ of free graded $R$-modules, 
  and a complex for $\cH^i(D)$ is obtained by tensoring $C(D)$ with 
  $\Z$ over $R,$ so that free graded $R$-modules become free abelian 
  groups of the same rank. In practice, the computation of $\cH^i(D)$ 
  faces the problem of effectively simplifying complexes of abelian 
  groups of exponentially high rank. The shortcut for the computation of 
  $H^i(D)$ described in Section~\ref{computations} works equally 
  well for groups $\cH^i(D),$ with Proposition~\ref{small-complex} 
  generalized to complexes $\overline{\cC}(D).$ This
  easily leads to a computation of 
  cohomology groups $\cH^{i,j}(T_{2,k})$ of the $(2,k)$ torus link
  $T_{2,k},$ oriented as in  Section~\ref{computations}.
    
 \begin{prop} Cohomology groups $\cH^{i,j}(T_{2,k}), k> 1$ are
  isomorphic to 
   \begin{equation*} 
   \begin{array}{lll}  
   \cH^{0,-k}(T_{2,k}) & = & \Z, \\
    \cH^{0,2-k}(T_{2,k}) & = & \Z, \\
    \cH^{-2j-1,-4j-2-k}(T_{2,k}) & = & \Z 
 \hspace{0.2in}\mbox{ for }1\le j \le \frac{k-1}{2}, j\in \Z, \\
    \cH^{-2j,-4j-k}(T_{2,k}) & = & \Z_2 
 \hspace{0.2in}\mbox{ for }1\le j \le \frac{k-1}{2}, j\in \Z, \\
    \cH^{-2j,-4j+2-k}(T_{2,k}) & = & \Z 
 \hspace{0.2in}\mbox{ for }1\le j \le \frac{k-1}{2},  j\in \Z, \\
   \cH^{-k,-3k}(T_{2,k}) & = & \Z \hspace{0.2in}  \mbox{ for even } k, \\
   \cH^{-k,2-3k}(T_{2,k}) & = & \Z \hspace{0.2in}  \mbox{ for even } k, \\
   \cH^{i,j}(T_{2,k}) & = & 0 \hspace{0.2in} \mbox{ for all other values of } 
   i \mbox{ and }j 
  \end{array}  
  \end{equation*} 
  \end{prop} 

\subsection{An application to the crossing number} 
\label{crossing-number} 

\begin{definition} A plane diagram $D$ with the set $\I$ of double points 
is called $+$adequate if for each double point $a$ the diagram 
 $D(\I\setminus\{ a\})$ has one circle less than $D(\I).$ 
\end{definition} 

\begin{definition} A plane diagram $D$ with the set $\I$ of double points 
is called $-$adequate if 
for each double point $a$ the diagram $D(\{ a\})$ has one circle less 
 than $D(\emptyset).$ 
\end{definition} 

\begin{definition} A plane diagram $D$ is called adequate if it is both  
$+$ and $-$adequate. 
\end{definition} 

These definitions are from [LT] and [T]. 

\begin{prop} 
\label{homology-sides} 
Let $D$ be a diagram with $n$ crossings. 
Then $\overline{\cH}^0(D)\not= 0$ if and only if $D$ is $-$adequate and 
     $\overline{\cH}^n(D)\not= 0$ if and only if $D$ is $+$adequate. 
\end{prop} 

\emph{Proof: } The differential 
$\partial^0 : \overline{\cal C}^0(D) \to \overline{\cal C}^1(D)$ 
is not injective and, hence, $\overline{\cH}^0(D)\not= 0$ 
 if and only if $D$ is $-$adequate. Similarly for $\overline{\cH}^n(D)$ and 
 $+$adequate diagrams (groups $\overline{\cH}^i(D)$ were defined 
 at the end of Section~\ref{small-cohomology}). $\square$

\begin{definition} Homological length $\hl(L)$ of an oriented 
 link $L$ is the 
difference between the maximal $i$ such that $\cH^i(L)\not= 0$ and 
the minimal $i$ such that $\cH^i(L)\not= 0.$ 
\end{definition} 

Denote by $c(L)$ the crossing number of $L$. It is the minimal number 
of crossings in a plane diagram of $L.$ 

\begin{prop} For an oriented link $L$
\begin{equation} 
c(L) \ge \hl(L)
\end{equation} 
\end{prop} 

\emph{Proof: } Let $D$ be a   diagram of $L$ with $c(L)$ 
crossings. Then $\overline{\cal C}^i(D)=0$ for $i<0$ and for $i> \hl(L).$ 
Consequently, $\overline{\cH}^i(D)=0$ for $i< 0$ and $i> \hl(L).$ 

$\square$ 

\begin{corollary}
\label{due-to-Th} Let $D$ be an adequate diagram with $n$ crossings 
of a link $L.$ Then $c(L)= n.$  
\end{corollary} 

\emph{Proof: } By Proposition~\ref{homology-sides} 
$\overline{\cH}^0(D)\not= 0$ and $\overline{\cH}^n(D)\not= 0.$ Therefore, 
$c(L)\ge \hl(L)\ge n.$ 
But since $D$ is an $n$-crossing diagram of $L,$ the 
crossing number of $L$ is $n.$ $\square$ 

Corollary~\ref{due-to-Th} was originally obtained by Thistlethwaite 
(Corollary 3.4 of [T]) through the analysis  of the $2$-variable Kauffman 
 polynomial (not to be confused with the Kauffman bracket). 

\section{Invariants of $(1,1)$-tangles}
\label{section-invariants-from-modules} 

\subsection{Graded $A$-modules} 
\label{graded-modules} 

Recall that in Section~\ref{here-is-A} we defined algebra $A$
as a free module of rank $2$ over the ring $R=\Z[c],$  
generated by $\mo$ and $X,$ with the multiplication rules  
 \begin{equation} 
 \mo\mo = \mo , \hspace{0.2in} \mo X = X\mo = X, \hspace{0.2in} X^2=0. 
 \end{equation}  
Gradings of $\mo$ and $X$ are equal to $1$ and $-1$, respectively, 
 so that the multiplication in $A$ is a graded map of degree $-1.$ 

\begin{definition} 
A graded $A$-module $M$ is a $\Z$-graded abelian group  
$M=\oplusop{i\in \Z} M_i,$ together with group homomorphisms
\begin{eqnarray}
 X : & M_i \lra M_{i-2}, & \hspace{0.2in}i\in \Z, \\
 c : & M_i \lra M_{i+2}, & \hspace{0.2in} i\in \Z, 
\end{eqnarray} 
 that satisfy relations 
\begin{equation} 
Xc = cX \hspace{0.3in} \mbox{and} \hspace{0.3in} X^2=0.  
\end{equation} 
\end{definition} 

\begin{definition} A homomorphism  of graded $A$-modules $M$ and $N$ is 
 a grading-preserving homomorphism of abelian groups $f: M \to N$ 
that intertwines the action of $X$ and $c$ in $M$ and $N$: 
\begin{equation} 
 Xf = fX, \hspace{0.5in} cf= fc . 
 \end{equation} 
\end{definition}

Denote by $A\mbox{-mod}_0$ the category with objects being 
graded $A$-modules and 
morphisms being grading-preserving homomorphisms of graded $A$-modules. 
 Note that $A\mbox{-mod}_0$ 
 is an abelian category.  Denote by $\{ n\}$ 
 the automorphism of $A$-mod that shifts the grading down 
 by $n.$ Let $A$-mod be the category of graded $A$-modules and graded 
maps. $A$-mod has the same objects as $A\mbox{-mod}_0,$ but more morphisms.

 Given a graded $A$-module $M,$ define the multiplication map 
  \begin{equation} 
  m_M:A\o_R M \lra M
  \end{equation} 
  by 
  \begin{equation} 
  m_M(\mo \o t) = t, \hspace{0.2in}
  m_M(X \o t) = Xt,  \hspace{0.2in}  t\in M. 
  \end{equation} 
  The multiplication map $m_M$ is a degree $-1$ map with the 
  grading on $A\o_R M$ defined as the product grading of gradings 
  of $A$ and $M$.

 To a graded $A$-module $M$ associate a map 
 \begin{equation} 
 \Delta_M: M \lra A\o M 
 \end{equation}  
 (recall that all tensor products are over $R=\Z[c]$) by 
 \begin{equation} 
 \Delta_M(t) = X \o t + \mo \o X t  + cX \o X t, \hspace{0.3in} t\in M.   
 \end{equation} 

Then $\Delta_M$ equips $M$ with the structure of a cocommutative 
comodule over $A.$ Map $\Delta_M$ has degree $-1.$ 
The following relation between $\Delta_M$ and $m_M$ is straightforward 
to check 
\begin{equation}
\Delta_M m_M= (Id_A\o m_M)(\Delta\o Id_M) = 
(m\o Id_M)(Id_A\o \Delta_M)  
\end{equation}
Denote by 
$\iota_M$ the map $M \lra A\o M$ given by $\iota_M(t)= \mo \o t$ for 
 $t\in M.$ 

Introduce an $A$-module structure on $A^{\o n}\o M$ by 
\begin{equation} 
a ( x \o y) =  x\o ay\mbox{ for } a\in A, x\in A^{\o n}, y\in M, 
\end{equation} 
where $ay$ means $m_M(a\o y).$ 
Then $m_M,\Delta_M$ and $\iota_M,$ after appropriate shifts by 
 $\{ 1\}$ or $\{ -1 \},$ are maps of graded $A$-modules. 

\begin{prop} 
\label{graded-module-splittings} 
For any graded $A$-module $M$ we have direct sum 
decompositions of $A\o M,$ considered as a graded $A$-module 
\begin{eqnarray}
A\o M & = & \Delta_M M \oplus \iota_M M \label{equality-one} \\
 \label{middle-equality} 
A\o M & = & \iota_M M \oplus (\Delta_M - \iota_M m_M \Delta_M) M  \\
A\o M & = & \Delta_M M \oplus (\iota_M - c \iota_M m_M \Delta_M) M 
\label{equality-three}  
\end{eqnarray}
\end{prop} 

\emph{Proof:} Let us check (\ref{middle-equality}), for instance. 
We have a decomposition 
\begin{equation} 
A\o M= (\mo \o M) \oplus (X \o M).
\end{equation} 
 Denote 
by $p$ the projection $A\o M \to X\o M,$ orthogonal to $\mo \o M.$ 
Since $\iota_M M = \mo \o M,$ it suffices to check that 
\begin{equation} 
p(\Delta_M -\iota_M m_M \Delta_M): M \to X\o M
\end{equation} is an $A$-module 
isomorphism. This map is given by 
\begin{equation} 
t\longmapsto X \o (1 +cX)t, \hspace{0.1in}\mbox{ where }t\in M. 
\end{equation} 
The inverse map is 
\begin{equation} 
X\o t \longmapsto (1-cX)t
\end{equation} 
Decompositions (\ref{equality-one}) and (\ref{equality-three}) 
can be verified analogously. $\square.$ 

\subsection{Non-closed (1+1)-cobordisms} 
\label{module-cobordisms} 

Let $\M1$ be the category whose objects are one-dimensional manifolds 
that are unions of a finite number of circles and one interval.
An ordering of the ends of this interval is fixed. 
Morphisms between objects $\alpha$ and $\beta$ of $\M1$ are oriented 
surfaces whose boundary is the union of $\alpha, \beta$ and 2 intervals 
that join corresponding ends of the intervals of $\alpha$ and $\beta.$ 
An example is depicted below. 
\drawing{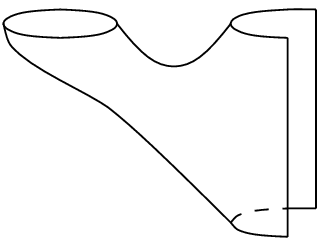}
This surface represents a morphism from an interval to a union of a
circle and an interval. 
We require that a surface can be presented as a composition of 
disjoint unions of surfaces $S_2^1,S_1^2,S_0^1,S_1^0,S_2^2,S_1^1,$  
defined in Section~\ref{A-and-surfaces}, and surfaces $T_1,T_2,$ depicted 
below

\drawing{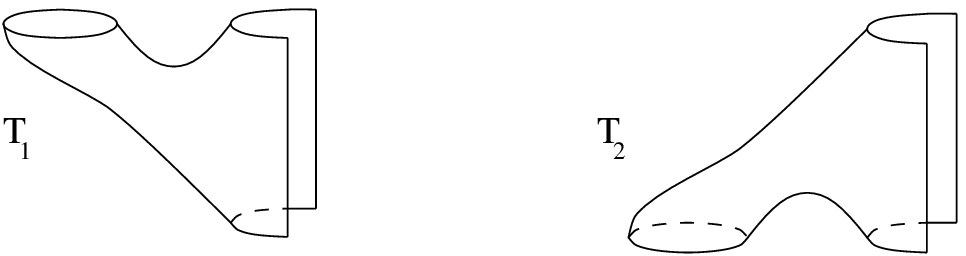}

We compose morphisms in this category by concatenating surfaces.  

Category $\M1$ is a module-category over $\MM,$ 
defined in Section~\ref{A-and-surfaces}.
The bifunctor 
$\MM \o \M1 \to \M1$ is defined on objects and morphisms by taking 
disjoint unions.  

Recall from the previous section that 
  $A$-mod denotes the category of graded $A$-modules and graded 
 homomorphisms. 
Given a graded $A$-module $M$, define a monoidal functor 
$$F_M: \M1 \lra A\mbox{-mod} $$
by assigning the graded $A$-module $A^{\o n} \o M$ 
to a union of $n$ circles and one interval, maps 
$\Delta_M$ and $m_M$ (defined in the previous section) 
to the elementary surfaces $T_1$ and $T_2$: 
\begin{equation} 
F_M(T_1) = \Delta_M, \hspace{0.1in} F_M(T_2)= m_M. 
\end{equation} 

To the other six elementary surfaces $S_2^1,S_1^2,S_0^1,S_1^0,S_2^2,S_1^1$ 
(Section~\ref{A-and-surfaces}) associate 
the same maps as for the functor $F$: 
\begin{equation}
\label{functor-FM-defined}
\begin{array}{lll} 
 F_M(S_2^1) =m, &  
  F_M(S_1^2) = \Delta, & 
  F_M(S_0^1)= \iota,  \\
  F_M(S_1^0)=\epsilon, & 
  F_M(S_2^2)= \mbox{Perm}, & 
  F_M(S_1^1)= \mbox{Id}. 
\end{array}  
\end{equation}

\subsection{$(1,1)$-tangles} 
\label{section-loong-links} 

A $(1,1)$-tangle is a proper smooth embedding 
$e: T \hookrightarrow \R^2\times [0,1]$ 
 of a finite collection $T$ of circles and one interval $[0,1]$ 
 into $\R^2 \times [0,1]$ such that the boundary points of $[0,1]$ go 
 to the corresponding boundary component of $\R^2\times [0,1]:$ 
  \begin{equation} 
  e(0) \in \R^2\times \{ 0\}, \hspace{0.2in} 
  e(1) \in \R^2\times \{ 1\}. 
  \end{equation} 
 
Two $(1,1)$-tangles are called equivalent if they are 
isotopic via an isotopy that fixes the boundary. 
 Oriented $(1,1)$-tangles are 
  $(1,1)$-tangles with a chosen orientation of 
  each component,  orientation of $[0,1]$ always chosen 
  in the direction from $0$ to $1.$ 

 Define a marked oriented link (in $\R^3$) as an oriented link with a 
 marked component. Obviously, there is a natural one-to-one correspondence 
  between marked oriented links and oriented $(1,1)$-tangles: the closure 
  of a $(1,1)$-tangle is a marked oriented link. 
 Denote this map from oriented $(1,1)$-tangles to oriented 
  marked links by $\mbox{cl}.$

\subsection{Invariants} 
\label{inv-long-links}

 A plane diagram $D$ of an oriented $(1,1)$-tangle $L$ is a generic projection 
of $L$ onto $\R\times [0,1].$

 If $D$ is a plane diagram of an oriented $(1,1)$-tangle, define $x(D)$ and 
 $y(D)$ in the same way as for plane diagrams of oriented links 
 (see Section~\ref{Kauffman-bracket}). 

Fix a graded $A$-module $M.$ Let $n$ be the number of double points of $D$, 
so that $n=x(D)+y(D)$ and $\I$ the set of double points of $D.$ 
 To $M$ and $D$ associate 
a commutative 
 $\I$-cube $V_D^M$ over the category $A\mbox{-mod}_0$ of graded $A$-modules 
 as follows. 

For $\L \subset \I$ the $\L$-resolution $D(\L)$ of 
$D$ consists of a disjoint union of circles and an interval. 
The functor $F_M$ (see Section~\ref{module-cobordisms})
assigns a graded $A$-module to $D(\L).$ Define 
\begin{equation} 
V_D^M(\L) = F_M(D(\L))\{ -|\L|\} .
\end{equation} 
Maps between $V_D^M(\L)$ for various subsets $\L$ 
are defined by the procedure completely 
analogous to the one described in Section~\ref{diagrams-to-comcubes}. 
Due to shifts $\{ -|\L|\},$ 
these maps of graded $A$-modules are grading-preserving, rather than  
 just graded maps, so that $V_D^M$ is a commutative cube over 
$A\mbox{-mod}_0.$ 

\emph{Example:}  For a diagram $D,$ depicted below, 

\drawing{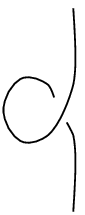} 

resolutions of $D$ are 

\drawing{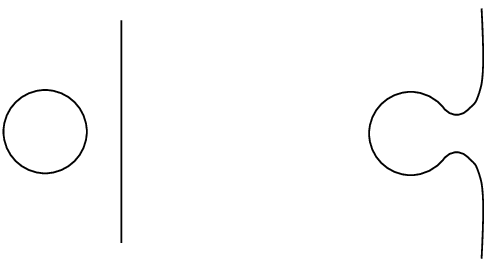} 

so that 
\begin{eqnarray*} 
V_D^M(\emptyset) & = & A\o M  \\
 V_D^M(\I)& = & M\{ -1\} 
\end{eqnarray*} 
 and the structure map $V_D^M(\emptyset) \lra V_D^M(\I)$ 
is the multiplication map $m_M: A\o M \to M\{ -1 \}.$ 

Next we transform the commutative $\I$-cube $V_D^M$ into a skew 
 commutative $\I$-cube 
 by putting minus signs in front of some structure maps 
 of $V_D^M,$ or, equivalently, by tensoring it with $E_{\I}.$ 
 Denote by $\C_M(D)$ the complex 
 $\C(V_D^M\o E_{\I})$ of graded $A$-modules. Define 

\begin{equation} 
C_M(D) = \overline{C}_M(D) [ x(D) ]\{ y(D)-2x(D)\}
\end{equation} 

Denote the $i$-th cohomology group of the complex $C_M(D)$ by 
 $H^i(D,M).$ These cohomology groups are graded $A$-modules. 
Denote the $j$-th graded component of $H^i(D,M)$ by $H^{i,j}(D,M)$  

\begin{theorem}
\label{theorem-long-i} 
 For a graded $A$-module $M,$ an oriented $(1,1)$-tangle $L$ and a
diagram $D$ of $L,$ isomorphism classes of graded $A$-modules $H^{i}(D,M)$ 
do not depend on the choice of $D$ and are invariants of $L.$ 
\end{theorem} 

Our proof of Theorem~\ref{closed-i} immediately generalizes without 
  essential modifications to a proof of Theorem~\ref{theorem-long-i}. 
  Proposition~\ref{graded-module-splittings} is used to establish 
 direct sum decompositions of $C_M(D)$, analogous to decompositions 
 of  $C(D)$, for suitable $D,$ given by 
  Propositions~\ref{comcube-splitting-iota-jmath},
  \ref{comcube-splitting-aleph-delta}, 
  \ref{prop-tangency-move}.1, 
  \ref{prop-last}.1. 

$\square$ 

  Cohomology groups $H^i(D),$ 
  defined in Section~\ref{diagrams-to-comcubes}, is a special 
  case of groups $H^i(D,M),$ as the next proposition explains. 

\begin{prop} Let $D$ be a diagram of an oriented $(1,1)$-tangle 
 $L$ and denote by $\cl (D)$ the associated   diagram 
 of the marked oriented link $\cl(L).$ 
 Considering $A$ as a graded $A$-module,  
 we have a canonical  isomorphism of cohomology groups 
 (as graded $R$-modules)  
 \begin{equation}
  H^i(D,A) \cong H^i(\cl(D)), \hspace{0.5in} i\in \Z. 
 \end{equation} 
 \end{prop} 

$\square$ 

Given a finitely-generated graded $A$-module $M,$ define the graded Euler 
 characteristic $\widehat{\chi}(M)$ by 
\begin{equation} 
  \widehat{\chi}(M) = \sum_{j\in \Z} \ddim_{\Q}(M_j\o_{\Z}\Q)
 \end{equation} 

\begin{prop} 
Let $M$ be a finitely-generated graded $A$-module, $L$ an oriented
 $(1,1)$-tangle and $D$ a diagram of $L$. Then 
\begin{equation} 
\frac{K(\cl(L))\widehat{\chi}(M)}{q+q^{-1}} = 
\sum_{i,j\in \Z} (-1)^i q^j \ddim_{\Q}(H^{i,j}(D,M)\o_{\Z}\Q).
\end{equation} 
that is, the Kauffman bracket of $\cl(L)$ is proportional to the Euler 
characteristic of groups $H^{i,j}(D,M).$  
\end{prop} 

$\square$ 

Given two graded $A$-modules $M,N$ and a grading-preserving homomorphism 
 $f: M\to N$, it induces a map of commutative cubes $V_D^M\to V_D^N,$ which, 
 in turn, induces a map of complexes $C_M(D)\to C_N(D)$ and a 
 map of cohomology groups $H^i(D,M) \to H^i(D,N).$  
 So, in fact, each diagram $D$ of an oriented long link 
  defines functors 
 $H^i_D$ from the category $A\mbox{-mod}_0$ of graded $A$-modules to itself, 
  $H^i_D(M)= H^i(D,M).$ If two diagrams $D_1,D_2$ 
  are related by a Reidemeister move, constructions of 
  Section~\ref{section-transformations} extend to functor 
  isomorphism $H^i_{D_1}\stackrel{\cong}{\lra}H^i_{D_2}.$ Let us frame 
  this observation into a proposition. 

 \begin{prop} 
  \label{no-more-propositions} 
 For an oriented (1,1)-tangle $L$ and a diagram $D$ of 
 $L$ isomorphism classes of functors 
 \begin{equation}   
  H^i_D: A\mbox{\rm -mod}_0 \lra A\mbox{\rm -mod}_0
  \end{equation}   
  do not depend on the choice of $D$ and are invariants of $L.$ 
 \end{prop} 

  Oriented long links with one component correspond one-to-one 
  to oriented knots in $\R^3.$ Thus, 
  Proposition~\ref{no-more-propositions} gives invariants of 
  oriented knots in $\R^3.$ Moreover, if a diagram $D$ represents 
  an oriented knot 
  and $D'$ is the diagram obtained from $D$ by reversing the 
  orientation of the underlying curve, there is 
  a natural in $M$ isomorphism $H^i(D,M) = H^i(D',M).$ 
  Consequently, for knots, isomorphism classes of functors $H^i_D$ do not 
  depend on the orientation, and $H^i_D$ provide ``functor-valued'' 
  invariants of non-oriented knots. Of course, these invariants depend 
  on how the ambient 3-space is oriented. 
  
Let $D$ be the 3-crossing diagram of  
the left-hand trefoil (knot $T_{2,3}$ in the notations 
of Section~\ref{computations}). The functors $H^i_D$ are written below

  \begin{eqnarray*} 
   H^{-3}_D(M) & = &  \kker \hspace{0.05in} 2X(M)\{ 8\}, \\
   H^{-2}_D(M) & = &  (M/(2XM))\{ 6 \}, \\
   H^{0}_D(M) & = &  M\{ 2\},\\
   H^{i}_D(M) & = &  0 \mbox{ for all other values of }i,  
  \end{eqnarray*} 
  where $\kker \hspace{0.05in} 2X(M)= \{ t\in M| 2Xt=0\}.$

\textsc{Institute for Advanced Study, Princeton} 

\textit{E-mail address:} 
\textsf{mikhailk@math.ias.edu}

\end{document}